\documentclass[singlespace]{article} % default 10pt
\usepackage{cite}
\usepackage{hyperref}
\usepackage{url}
\usepackage{amsmath}
\usepackage{bm}
\usepackage{amsfonts}
\usepackage{algorithm}
\usepackage{algorithmic}
\usepackage{amsthm,paralist}
\usepackage{graphicx}
\usepackage{color}
\usepackage{subfigure}   % For two subfigures put together

\usepackage{diagbox} % draw slash in the table    \usepackage{slashbox}
\usepackage{stfloats}  % put table at bottom
\usepackage{footnote}
\usepackage{multirow}  % multi-row use
\usepackage[scriptsize,bf]{caption}
\usepackage{amssymb}
\usepackage{diagbox} % draw slash in the table    \usepackage{slashbox}
\usepackage{stfloats}  % put table at bottom

\def\argmin{\mathop{\rm argmin}}

\def\R{\mathbb{R}}
\def\C{\mathbb{C}}

\def\avg{\mathrm{avg}}

\newcommand{\beq}{\begin{equation}}
\newcommand{\eeq}{\end{equation}}

\newcommand{\trace}{{\mbox{\textrm{\rm Tr}}}}

\newcommand{\cO}{{\mbox{$\mathcal{O}$}}}

\newtheorem{theorem}{Theorem}[section]
\newtheorem{lemma}{Lemma}[section]

\newtheorem{proposition}{Proposition}[section]

\newtheorem{claim}{Claim}[section]

\newtheorem{conjecture}{Conjecture}[section]

\def\smskip{\par\vskip 5 pt}

\def\QED{\hfill{\bf Q.E.D.}\smskip}

\def\O{\mathcal{O}}
\def\tO{\tilde{\mathcal{O}} }

\title{ Worst-case Complexity of Cyclic Coordinate Descent: $O(n^2)$ Gap with Randomized Version }
\author{{Ruoyu Sun} \thanks{Department of Industrial and Enterprise Systems Engineering,  Univeristy of Illinois at Urbana-Champaign, Urbana, IL. \texttt{ruoyus@illinois.edu}.}  { \quad \quad  \quad   Yinyu Ye}\thanks{Department of Management Science and Engineering,  Stanford University, Stanford, CA. \texttt{yyye@stanford.edu}. } \
	\date{ }
}

\begin{document}

\maketitle

\vspace{-0.5cm}
\begin{abstract}

 This paper concerns the worst-case complexity of cyclic coordinate descent (C-CD) for minimizing a convex quadratic function, which is equivalent to Gauss-Seidel method and can be transformed to Kaczmarz method and projection onto convex sets (POCS).
We observe that the known provable complexity of C-CD can be $\O(n^2)$ times slower than randomized coordinate descent (R-CD), but no example was rigorously proven to exhibit such a large gap.
In this paper we show that the gap indeed exists. We prove that there exists an example for which C-CD takes at least $\O(n^4 \kappa_{\text{CD}} \log\frac{1}{\epsilon})$ operations,  where $\kappa_{\text{CD}}$ is related to Demmel's condition number and it determines the convergence rate of R-CD.
It implies that in the worst case C-CD can indeed be $\O(n^2)$ times slower than R-CD, which has complexity  $\O( n^2 \kappa_{\text{CD}} \log\frac{1}{\epsilon})$.
Note that for this example, the gap exists for any fixed update order, not just a particular order.
Based on the example, we establish several almost tight complexity bounds of C-CD for quadratic problems. 
One difficulty with the analysis is that the spectral radius of a non-symmetric iteration matrix does not necessarily constitute a \textit{lower bound} for the convergence rate.

An immediate consequence is that for Gauss-Seidel method, Kaczmarz method and POCS, there is also an $\O(n^2) $ gap between the cyclic versions and randomized versions (for solving linear systems). 
We also show that the classical convergence rate of POCS by
Smith, Solmon and Wager [1] is always worse and sometimes can be infinitely times worse than our bound.

\end{abstract}

\section{Introduction}

Coordinate descent (CD) algorithms have been very popular recently due to their efficiency for solving large-scale optimization problems (see, e.g., \cite{wright2015coordinate} for a recent survey).
In the most basic form, cyclic CD (C-CD) optimizes over one variable at a time with other variables fixed, and the variables are chosen according to a fixed order.
Due to the simplicity, CD methods are one of the most widely used class of optimization methods in science and engineering. Its applications include tensor decomposition \cite{kolda2009tensor}, libsvm package for SVM in machine learning \cite{chang2011libsvm, hsieh2008dual}, glmnet package for Lasso in statistics \cite{friedman2010regularization, bradley2011parallel,mazumder2012sparsenet}, resource allocation in wireless communications \cite{razaviyayn2013unified,baligh2014cross,sun2013long,hong2013joint}, to name a few;
see some other applications in \cite{canutescu2003cyclic,bouman1996unified,,wen2012block,
	sun2015guaranteed}.

For the theoretical analysis, most early works focused on the exact conditions for the convergence (e.g., Powell \cite{powell1973search}, Bertsekas \cite{bertsekas99}, Tseng\cite{tseng01}, Grippo and Sciandrone \cite{Grippo00}) and the quality of convergence  (e.g. Luo and Tseng \cite{Luo92CD}).
A landmark in the history of CD algorithms is the establishment of the explicit convergence rate of randomized CD (R-CD) \cite{leventhal2010randomized,nestrov12}, a variant which  updates variables randomly. In particular, R-CD was shown to be $\O(1)$ to $\O(n)$ times faster than GD, where $n$ is the number of variables. Note that the introduction of randomized update order is crucial since it makes the analysis of CD methods quite simple.
Ever since then, randomized update rule has been a new standard for theoretical analysis of CD-type methods
\cite{shalev2013stochastic,richtarik12,lu2015complexity,qu2014randomized,
	lin2014accelerated,zhang2015stochastic,fercoq2015accelerated,liu2015asynchronousJMLR,
	patrascu2015efficient,hsieh2015passcode,wright2015coordinate}.
%For example, the results on the convergence speed of randomized versions have been extended to composite objective functions % \cite{richtarik12,lu2015complexity}.
Furthermore,  accelerated R-CD  was shown to have better  complexity than conjugate gradient method 
(in some parameters) when  solving symmetric PD (positive definite) linear systems  % $\mathcal{O}(\frac{ n L }{ \sum_i L_i })$  over conjugate gradient method
\cite{lee2013efficient}, and can improve the complexity of solving packing and covering LP (Linear Programming) \cite{allen2015nearly}.

With all these nice theoretical results on R-CD, one may wonder whether the same results can be achieved for C-CD. There are several reasons for studying cyclic methods.
% the above ``cut and load'' process may happen several times.
 % If the communication between the cache, the memory and the hard disk is fast enough, applying randomized coordinate selection can be feasible;
 % otherwise it might be beneficial to process all coordinates in a fixed order.
  %loaded into CPU one by one, and all coordinates loaded into the first-level cache one by one, etc.
  % While it is not easy to design an algorithm tailored for a specific hardware architecture (and new architectures are emerging), our point is that
 % In short, cyclic CD methods might be more hardware-friendly than randomized versions in some scenarios.
(1) The complexity of deterministic algorithms is theoretically important (partly because generating random bits is highly non-trivial). For example, the first polynomial-time deterministic algorithm for PRIME was regarded as a great achievement \cite{agrawal2004primes}.
% , even though polynomial time randomized algorithms had been known for decades.
Another example is the interesting open question whether there exists a version of deterministic simplex method that can solve LP in polynomial time.
 % It also remains open whether 
% deterministic algorithms can achieve the same complexity as the randomized methods for symmetric PD linear systems and positive LP analyzed in \cite{lee2013efficient, allen2015nearly}.
% In our case, as R-CD is the crucial components of fast algorithms for solving PD linear system and positive LP, it is an interesting question whether the same complexity can be achieved by deterministic algorithms. Understanding the complexity of C-CD will serve as an important first step towards answering this question.
(2) The study of the cyclic order may help us understand other update orders.
% In fact, for many orders such as cyclic, random permutation and double sweep, it was unknown what the worst-case convergence rates are. 
For instance, the random permutation order was observed to perform very well in practice for CD, SGD and ADMM  \cite{zhang2015stochastic,recht2013parallel, sun2015expected}, but the best known convergence rate bounds of randomly permuted CD are almost the same as that of cyclic versions \cite{sun2015improved} (except for some special cases \cite{lee2016random}).
%  tight convergence rate analysis remains  open. 
%  (except for some special cases \cite{lee2016random})).
% the best known convergence rate bounds of randomly permuted CD are almost the same as that of cyclic versions  (better rates are given for some special cases \cite{lee2016random}).
(3) In practice, it is not always easy or desirable to randomly pick coordinates.
% A computer system usually divides the whole set of variables in the hard disk into many chunks and load one chunk of variables into the memory. 
The typical computer architecture consists of multiple layers including caches, memory and hard disk, and fully randomized coordinate selection might be time-consuming when the communication between components of the system is not very fast. For example, it was pointed out  in \cite{recht2013parallel} that the sampling time of randomized order is not negligible. 
In certain distributed optimization algorithm \cite{xiao2017dscovr}, independently randomized order was deliberately avoided due to specific  design requirement. % , and a mixture of cyclic and permutation order was proposed. 
% , but the analysis was only performed for independently randomized order that is infeasible in practice. 
(4) Many practitioners are still using cyclic versions of CD; one example is that statisticians are still using
cyclic CD to solve Lasso \cite{yang2014coordinate}. If cyclic CD performs well in practice and already implemented in software packages, why would one change to randomized versions?
A more comprehensive understanding of different update orders may help practitioners choose
an appropriate update order.

% (4) The study of C-CD may shed light on other related methods, such as POCS and ADMM. 
%  which are much worse than randomized versions. 

There have been some recent efforts to understand the convergence speed of C-CD
 \footnote{	In fact, the analysis applies to cyclic BCGD (Block Coordinate Gradient Descent) for solving convex problems. 
 	 % BCGD for quadratic problems is equivalent to block successive over-relaxation method.
	%I n BCGD, a typical stepsize for the $i$-th block is $1/L_i$, and for randomized BCGD this stepsize can lead to an optimal rate.
	For minimizing convex quadratic functions, cyclic CGD with a special stepsize is the same as cyclic CD (i.e. exactly minimizing each subproblem).  }
 \cite{Beck13, Beck13b, Saha10, hong13complexity, sun2015improved}.
% but there are two strange phenomenons based on the theoretical bounds.
 For simplicity, we will discuss these bounds for applying C-CD to strongly convex quadratic functions $x^T A x - 2b^T x$,  which is equivalent to Gauss-Seidel method, Kaczmarz method and POCS in this special setting
 (see discussions later).
  We further assume the coefficient matrix $A$ has equal diagonal entries. 
  Suppose the maximum eigenvalue, minimum eigenvalue and average eigenvalue of $A$ are
  $\lambda_{\max},\lambda_{\min},\lambda_{\avg}$ respectively,
  the condition number $\kappa \triangleq \lambda_{\max}/\lambda_{\min} $,  and
  $\kappa_{\mathrm CD} \triangleq \lambda_{\avg}/\lambda_{\min} $. 
  % $ \lambda_{\max} = \lambda_1 \geq \dots \geq \lambda_n = \lambda_{\min}$, and define $\lambda_{\text{avg}} = \sum_i \lambda_i /n $.
 It is well-known that the  complexity of GD is $\tO( n^2 \kappa ) $,
  and the complexity of R-CD is $\tO( n^2  \kappa_{\mathrm CD} ) $,
  in which we ignore an $\O(\log 1/\epsilon)$ factor.
  This implies that R-CD is $\tau \triangleq \lambda_{\max}/\lambda_{\text{avg}} \in [1,n ]$ times better than GD; here, 
%  The quantity $\lambda_{\avg}/\lambda_{\min} $ is the condition number,
%  and $ n \lambda_{\avg}/\lambda_{\min}  $ is  Demmel condition number.
  note that the gap $ \tau $ can be as large as $n$.
  % and for the Wishart random matrix $\tau$ is usually of order $ \O(\sqrt{n})$ \cite{edelman1988eigenvalues}.
%  \begin{itemize} \item  \tO(n^2 \frac{ \lambda_{\max}^2  }{ \lambda_{\avg} \lambda_{\min} } \log^2 n  ) = 
The best known complexity of C-CD for quadratic problems is approximately  
  $\tO(n^2 \tau \kappa \log^2 n  ) $   % \cite{oswald1994convergence} % \cite{sun2015improved}
 %  \footnote{The reference  \cite{oswald1994convergence} discussed Gauss-Seidel method and seems to be unknown to optimization researchers. 
 %  	We first derived this bound by applying the analysis of the paper \cite{sun2015improved} (which studies non-strongly cnovex problems) to strongly convex quadratic problems. When preparing this paper we became aware of 	the early reference \cite{oswald1994convergence} which obtained a similar bound. }
   % See details in Proposition \ref{prop: upper bound} and Section \ref{sec:}. 
 %  	 and the proved complexity is $O(n^2 \tau L \log^2 n /\epsilon )$ for quadratic problems, which is $\tau$-times worse than GD; but the same analysis in \cite{sun2015improved} will lead to a  bound of $\tO(n^2 \tau \kappa \log^2 n ) $ for strongly convex quadratic functions (see the proof of Proposition \ref{prop: upper bound}).  },
  which is at least $\tau$-times worse than GD and $\tau^2$ times worse than R-CD.
 This theoretical bound does not match the numerical experiments which almost always show that C-CD converges much faster than GD. The existing results seem so weak that they even make a wrong prediction on whether C-CD is faster than GD.
  The potential $\O(n^2)$ gap between C-CD and R-CD also seems quite strange, as such a huge gap has not been reported by practitioners.  % observed in real world applications. 

It is very tempting to think that we might be able to prove C-CD is faster than GD, or even comparable to R-CD.
The discrepancy between the theory and the practice might just be because of the weakness of the proof techniques.
% might just attribute to the weakneess of the theory.
This impression may be enhanced when we reflect on the existing proofs of upper bounds.
The proof idea of \cite{Beck13,sun2015improved} is to view C-CD as an inexact version of GD,
and the major effort is spent on bounding the difference between C-CD and GD.
One obvious drawback of such a proof framework is that it cannot show a better convergence rate than GD; even if the difference is
zero, only the same rate would be established.
%--------------------------------------------
%This raises a question related to deeper understanding of CD-type methods: should we really view CD-type methods as inexact GD?
%The fact that R-CD is faster than GD implies that randomized versions should \emph{not} be viewed as inexact GD.
%--------------------------------------------
It seems possible that there exists a different proof framework for C-CD that leads to better convergence rates. To understand this issue is the main purpose of this paper. 
% To understand the ``best convergence rate'' (or tight convergence rates)  of C-CD is the main purpose of this paper.
%Understanding whether a better convergence rate exists or not is the main motivation of this paper.
% However, in this paper we will show that the large gap between C-CD and GD/R-CD really exists, at least when the convergence rate is expressed in terms of the current parameters. 

%To summarize, the known theory explains neither why C-CD is faster than GD in practice, nor
%why larger stepsize within range $[1/L, 1/L_i]$ leads to faster convergence of cyclic CGD in practice. In addition, the existing proof framework seems to be not strong enough to provide the tightest bound. However, we will show that the known results are almost the best we can hope for.

% If C-CD were really faster than GD, there should be a different proof framework based on new perspectives on CD methods.
% the ``true'' worst-case complexity of C-CD puzzling. It is possible that from the beginning we should have chosen
% when the block size is greater than $1$ or the objective function is not quadratic, cyclic BCGD with stepsize $1/L_i$
% is different from cyclic BCD.

\subsection{CD, Gauss-Seidel Method, POCS, Kaczmarz method}\label{subsec: review these algorithms}
% The earliest form of CD is perhaps  , a parallel version of G-S method 
In this subsection, we review several closely related methods: Gauss-Seidel method, Kaczmarz method and POCS (Projection Onto Convex Sets, a.k.a., alternating projection method). 
We will see that they are equivalent in the simple yet important setting of solving linear systems, thus understanding convergence speed is a common issue for all these methods. 
% also applied to these methods. 

Gauss-Seidel (G-S) method, first proposed by Gauss and Seidel in 19th centry, is one of the oldest iterative algorithms.
It can be used to solve any  system of linear equations, though the convergence is only guaranteed when the coefficient matrix
satisfies some diagonally dominant properties or is symmetric PSD (Positive SemiDefinite).
Regarding the convergence speed, it is well-known that for some special matrices, asymptotically G-S method converges twice as fast as Jacobi method (see, e.g., \cite{greenbaum1997iterative}).

% G-S method is closely related to 
POCS is a method to solve the convex feasibility problem, i.e., find a point in the intersection of closed convex sets. 
POCS has found many applications in applied mathematics and engineering; see, e.g., a survey of ten applications of POCS by Deutsch \cite{deutsch1992method}. 
The convergence of POCS   was proved by Von Neumann for two sets  in 1933 \cite{von1950functional} and Halperin for more than two sets \cite{halperin1962product}. 
The convergence rate of POCS was given by Smith, Solmon and Wagner \cite{smith1977practical},
and improved by a few works (e.g. \cite{kayalar1988error, deutsch1997rate}). For a detailed review of numerous works in this field, we refer the readers to Bauschke, Borwein and Lewis \cite{bauschke1997method},  Escalante and Raydan \cite{escalante2011alternating} and Galantai \cite{galantai2013projectors}.

Kaczmarz method is an old method for solving linear systems of equations proposed in 1937 \cite{kaczmarz1937angenaherte}.  A recent work \cite{strohmer2009randomized}
proved explicit convergence rate of randomized Kaczmarz method, which motivated works on R-CD.
Note that Kaczmarz method is a special case of POCS when
when the sets are hyperplanes.

The basic versions of G-S method, POCS, Kaczmarz method and CD are equivalent.
As mentioned above, Kaczmarz method is a special case of POCS. 
  Under a basis transformation, Kaczmarz method  is equivalent to G-S method for solving a symmetric PSD linear system, which is equivalent to cyclic CD (C-CD) for minimizing convex quadratic functions (see Appendix \ref{appen: Proof of equivalence of GS and Kaczmarz}). 
Note that G-S, POCS and CD are not equivalent in more general settings; in fact, G-S can be used to solve non-symmetric linear systems, POCS can be used to find intersection of any closed convex sets, and CD can be used to solve non-quadratic non-smooth problems. 
It seems not easy to obtain a unified convergence analysis for all of them.
Nevertheless, to understand the worst-case complexity, we need to first study the simplest setting, in which these methods are equivalent and thus can be analyzed altogether. 
In particular, the major question we want to answer is:
\begin{equation*}
\begin{split}
& \quad  \quad  \text{ For coordinate descent, G-S method, POCS and Kaczmarz method, 
is there an $\O(n^2)$ gap between} \\
& \text{ the worst-case convergence rate of their cyclic versions and randomized versions? }
\end{split}
\end{equation*}
% The equivalence of these methods for solving linear systems indi

\subsection{Summary of Contributions}
%Accelerated randomized BCD methods can be efficiently implemented for certain problems
%including the quadratic minimization \cite{lee2013efficient}, leading to better complexity for solving PSD linear systems than the Conjugate Gradient method.
% The focus has shifted from characterizing conditions for convergence to proving explicit rate of convergence.
% (choosing a stepsize can be viewed as solving a regularized version of the subproblem).

We will focus on the worst-case complexity of C-CD for minimizing convex quadratic functions
$ \min_{x \in \R^{n}} \frac{1}{2}x^T A x - b^T x $.
As discussed in Section \ref{subsec: review these algorithms}, in this simple setting, C-CD is equivalent to 
G-S method, Kaczmarz method and POCS. 
% This algorithm is equivalent to Gauss-Seidel method for solving a symmetric PSD linear system $Ax = b$, and Kaczmarz method for solving a linear system $U x = b$ where $A = U^T U$, and POCS for finding a  point in the intersection of hyperplanes defined by the equations in the system $U x = b$. 
% The extension to the case with block size more than $1$ and to inexact update is left as future work. 
%  It is not hard to extend our results to cyclc BCD (the block size can be larger than $1$) and cyclic BCGD with different stepsizes, but for simplicity we will not discuss these extensions.
In the following, we will say an algorithm has complexity $\tO( g(n, \theta) )$,
if it takes $\O( g(n, \theta)  \log(1/\epsilon) )$ unit operations % (plus, minus, multiplication and division)
to achieve relative error $ \frac{f(x) - f^*}{f(x^0) - f^* } \leq \epsilon $.
It is well-known that GD has complexity $\tO( n^2 \kappa ) $ and R-CD has complexity $ \tO( n^2 \kappa_{\mathrm CD} ) $, where $\kappa = \lambda_{\max}/\lambda_{\min}$ is the condition number, $ \kappa_{\mathrm CD} =  \lambda_{\avg} / \lambda_{\min} $. Denote $\tau \triangleq \lambda_{\max}/\lambda_{\text{avg}} \in [1,n ]$. 
%Assumption: All diagonal entries of A are equal to $L_1$.
% This assumption is mild since it be achieved by double-sided diagonal preconditioning.

We summarize our results, when specialized to the equal-diagonal case (i.e. all diagonal entries are the same) in the following table, ignoring a factor of $\O(\log 1/\epsilon)$. The non-equal-diagonal case is quite subtle and related to conjectures on Jacobi-preconditioning; see Section \ref{sec: Jacobi Precond}. Our main contribution is to establish several lower bounds by analyzing the convergence rate of a simple class of examples. Our discovery is that the upper bounds are ``almost'' tight (up to $\O(\log^2 n)$ factor) in the equal-diagonal case. More specifically, the table shows the following results:

\begin{table}[!htbp]
	\centering
	% \caption{Complexity of C-CD, GD and R-CD for equal-diagonal case (ignoring a $\log 1/\epsilon$ factor)} \label{table of complexity}
	\begin{tabular}{|c|c|c|c|}
		\hline 
		Parameters	&   $\kappa$ &  $\kappa$ and $\tau$ &  $\kappa_{\mathrm CD}  ( = \kappa/\tau) $    \\ 
		\hline 
		C-CD Upper bound (Proposition \ref{prop: upper bound})	&   $   n^3 \kappa  $  &
		$   \frac{1}{10} n^2  \kappa \tau \log^2 n  $
		&  $   n^4 \kappa_{\mathrm CD} $    \\ 
		\hline 
		C-CD Lower bound (Theorem \ref{thm: lower bound})	&   $ \frac{1}{40 }   n^3 \kappa $  &
		$ \frac{1}{40 } n^2  \kappa \tau   $
		&  $  \frac{1}{40} n^4 \kappa_{\mathrm CD} $    \\ 
		\hline
	\end{tabular} 
\end{table}

% We discuss two types of bounds, one does not involve $\tau = \lambda_{\max}/\lambda_{\text{avg}} $ and the other does.

% (note that $\lambda_{\max}(A) = L, \lambda_{\text{avg}(A)} = L_1 = A_{11}$) 
% , neither of which dominates the other.
% in addition, later we will discuss why a stronger bound dominating both does not exist.
\begin{itemize}
	\item %For a simple comparison with GD (or randomized CD), we want to express the complexity in terms of $n$ and $\kappa$ (or $\kappa_{\mathrm CD}$)
	In terms of $\kappa$ or $\kappa_{\mathrm CD}$, the worst-case complexity of C-CD is
	\begin{equation}\label{lower bound, 1st type}
	\tO( n^3 \kappa ) \text{ or }  \tO( n^4 \kappa_{\mathrm CD} )  .
	\end{equation}
	Both bounds are tight up to constant factors. This implies that C-CD can be $\O(n)$ times slower than GD and $\O(n^2)$ times slower than R-CD.
	
	\item It is more precise to characterize the complexity using an extra parameter $\tau  $ together with $\kappa$ or $\kappa_{\mathrm CD}$. The lower bound for the complexity of C-CD is
	\begin{equation}\label{lower bound, 2nd type}
	\tO\left( n^2  \kappa \tau   \right) \text{ or } \tO\left( n^2  \kappa_{\mathrm CD} \tau^2   \right),
	\end{equation}
	which is $\tau$ times worse than GD  or  $\tau^2$ times worse than R-CD.
	The range of the gap $\tau$ is  $ [1,n]$ and can be large in most cases.
	% and for Wishart random matrices   $\tau = O(\sqrt{n} )$       \cite{edelman1988eigenvalues},
	% thus the gap $\tau$ or $\tau^2$ can be very large.
	These two bounds are ``almost'' tight as they are only $\O(\log^2 n)$-times smaller than the upper bounds.
	% There is an $\O(\log^2 n)$ gap between these lower bounds and the corresponding upper bounds, but the $\log^2(n)$ factor was previously known
	% to be unavoidable in some cases. See more detailed discussion in Section \ref{sec: discuss two types of bounds}.
\end{itemize}

To prove the lower bounds, we only need to estimate the convergence rate of our specific examples, and there are at least two difficulties. Firstly, there is no closed form expression of the spectral radius of the iteration matrix and we need to consider the limiting behavior of a class of examples (still with fixed $n$).
Secondly, the spectral radius does not directly lead to a lower bound of the convergence rate
when the iteration matrix is \emph{non-symmetric}, and we need to explore some special structure of the examples.

%We observe that when the off-diagonal entries are large and thus $\tau = \lambda_{\max}/\lambda_{\avg} $ is large, C-CD is indeed slow, which is consistent with the theoretical bound of $\tilde{O}(n^2 \kappa \tau)$. However, in almost all scenarios  C-CD converges much \emph{faster} than GD, and the gap with R-CD is far smaller than $\O(\tau^2)$, which cannot be explained by the theoretical bounds. 
% \subsubsection{Numerical Findings}
Simulation shows that our worst-case bound is \textit{partially} consistent with the numerical experiments.
We perform  numerical experiments for dozens of random distributions of matrix $A$, and the relation between the numerical findings and the theory are summarized below.
%matches our theory. What matches worst-case bounds:  
% We summarize our numerical findings below (again, only consider the equal-diagonal case).
\begin{itemize}
	\vspace{-0.2cm}
	\item Our theoretical bound of $\tO(n^2 \kappa \tau)$ indicates that C-CD  converges slowly when $\tau = \lambda_{\max}/\lambda_{\avg} $ is large. Interestingly, we do observe that when the off-diagonal entries are large and thus $\lambda_{\max}/\lambda_{\avg} $ is large, C-CD is indeed slow.  This shows that the theory is partially consistent with the simulations. 
	 % What does not match worst-case bounds:
	% this phenomenon is indeed observed in our numerical experiments.
	%C-CD converges faster than R-CD when off-diagonal entries are small but slower than R-CD when off-diagonal entries are large.
	% This phenomenon is rather universal as it is observed for a wide variety of generative models of the coefficient matrix.
	
	\item  In almost all scenarios (except random perturbations of our example) C-CD converges much \emph{faster} than GD, which is \emph{opposite} to the theory.  % (or more precisely, about $\O( L/L_1)$ times)which says that C-CD converges $\O(1)$ to $\O(n)$ times \emph{slower} than GD.
	The gap between C-CD and R-CD in the experiments is far from the theoretical gap $\O(\tau^2)$. 	This discrepancy reveals the weakness of the worst-case analysis.
	% The theory also says that C-CD converges up to $\O(n^2)$ times slower than R-CD, but seems to be much smaller than $\tau$, and
	% Therefore, our experiments show that C-CD is still a very competitive algorithm
	%   even though we prove it converges slowly in the worst case.

	%   nevertheless, this issue has long been recognized and requires much more effort to resolve. 
	
	% \item % theory correctly predicts that large off-diagonal entries lead to slower convergence of C-CD, but cannot explain why C-CD is faster than GD.% recognized in numerical linear algebra and
\end{itemize}

%The  gap between the practice and the theory partially explains why the $\O(n^2)$ gap has never been reported.
% This gap has motivated our research, but our research does not fill in the gap.

%It was known that randomized algorithms are typically faster than cyclic versions, but no such theoretical results have been established
%(except that cyclic ADMM can diverge while randomly permuted ADMM converge in expectation for equally constrained quadratic minimization).
\subsection{Discussions }
We further discuss a few interesting issues related to this work. 

% present a few implications of our theoretical results and numerical findings.  % , as well open theoretical questions.
 \emph{ Gap Between Cyclic and Randomized Algorithms.}  We prove for the first time that C-CD, Gauss-Seidel method, Kaczmarz method and POCS can be $\O(n^2)$ times slower than their randomized counterparts. 
 Despite the long history of these algorithms, this $\O(n^2)$ gap was not rigorously established before.
 This is one of the few examples in continuous optimization that a large gap between a certain deterministic algorithm and its randomized counterpart is established.
    
    \emph{ Robustness of Worst-case Examples.} Our worst-case example appears to be quite robust. 
    % One may wonder whether a small perturbation of the algorithm or problem can fix the worst-case. 
    % A possible solution is to perturb the update order. 
    A common belief is that C-CD can be slow because  one particular order can be very bad, and randomly pick an order and fix it will be good.  Indeed, this is the case for Example 2 in Section \ref{sec: role of examples}.
    However, for our example, \textit{any fixed order} out of  all $n!$ possible orders is equally slow.
   % most orders out of all $n!$ possible orders might be good. 
    %    Intuitively, why is cyclic order possibly bad? One intuition is that there may exist some bad orders that the algorithm proceeds slowly (see Example 2 in Section \ref{sec: role of examples}). 
    % One may wonder   Indeed for Example 2, a random fixed order makes C-CD as fast as R-CD.  
    Another possible way to fix the worst-case example is to perturb the problem input. 
    In a different scenario,  a small perturbation of the problem input makes the complexity of the simplex method much better \cite{spielman2004smoothed}. However, under a small perturbation of our example, C-CD is still $\O(n^2)$ times slower than R-CD. 
    
    %\emph{Classical POCS bounds.} The convergence rate of C-CD was rarely studied before 2010, but the convergence rate of POCS has been studied since 1970s. What is the relationship between the 
    
    \emph{Role of Examples in Convergence Analysis.}  Our contribution is not just to provide an example that C-CD is much slower than R-CD. A single example itself says little, because there might exist another example that C-CD is much faster than R-CD. What is more interesting is how the example interacts with the theoretical bounds. There is an $\O(n^2)$ gap between existing bounds of C-CD and R-CD, and our contribution is to prove that our example matches both bounds of C-CD and R-CD, thus validating the $\O(n^2)$  gap. Not all examples can make the same ``achievement''. 
    See more discussions in Section \ref{sec: role of examples}.

\emph{Fundamental Gap Between Deterministic and Randomized CD?} Our results only establish a large gap between a single deterministic version of CD and R-CD.
 A natural question arises:
is there a fundamental gap between deterministic CD and randomized CD?
%is there a deterministic variant of CD that achieves the same complexity of R-CD? %  (which means $L/L_{\min}$ times better than GD)?
% be used to test various candidates.
There has to be an answer: either we can prove a large lower bound for \emph{all} deterministic CD methods, or we can find one deterministic CD that performs close to R-CD.
Both possibilities are very interesting.
For the latter possibility, % any candidate should be tested for our example. % , and it should match R-CD and be $n$ times faster than GD.
there are a few  candidates such as CGD with stepsize other than $1/L_i$ (equivalent to SOR, i.e. successive over relaxation) and double sweep method (a.k.a. symmetric SOR),
but they are far worse than R-CD for our example.

%----RESERVE-----
%CGD with stepsize $1/L$ does improve over C-CD for our example; in fact, it can match GD in almost all tests, but
%this essentially means it is much worse than C-CD in most practical scenarios.
%% One interesting phenomenon, related to the stepsize dilemma, is
% Not to mention matching R-CD, even a CD variant that theoretically outperforms GD while practically matches C-CD is unknown to us.
 %--------------
% While it may not be easy to derive the worst-case complexity of all variants,
% If the complexity is measured in $\kappa$, then the worst-case complexity of R-CD is actually
% (the theoretical upper bound \cite{sun2015improved} is just $\O(\log^2 n)$ times worse than GD),

 \emph{Deterministic Complexity.} Recent progress on the complexity of some important classes of problems (e.g. PD linear systems, positive LP) is based on randomized versions of CD methods.
      As we have established a large gap between C-CD and R-CD, it is unclear whether the same complexity can be achieved for deterministic algorithms.
   % The theoretically fastest iterative algorithm for solving linear systems is accelerated randomized BCD \cite{lee2013efficient,allen2015even}.
   % However, since our result establishes a gap between cyclic BCD and randomized BCD,
    For example, CG (conjugate gradient) is still the fastest \emph{deterministic} iterative algorithm for solving PSD linear systems,
   even though accelerated R-CD is faster in a probabilistic sense. % One natural question is: can we improve the \emph{deterministic complexity} for solving linear systems?

  \emph{Bridging the Gap Between Theory and Practice.} It is an interesting question how to explain the large discrepancy between the theory and the practical performance of C-CD. This kind of discrepancy may lead to novel theoretical advances. % , such as the smoothed analysis \cite{spielman2004smoothed}. 
 One famous example is the smoothed analysis developed by Spielman and Teng \cite{spielman2004smoothed} that aims to explain such a gap for the simplex method.
 %  why the simplex method has exponential time complexity but performs well in practice.
 What type of analysis is suitable for explaining the practical performance of C-CD (e.g. why is it usually much faster than GD)? Smoothed analysis is not enough as a small perturbation of our example still exhibits the large gap.
 This seems to be a difficult question that is currently beyond our reach.   
 We think one possibility is to introduce a new metric that measures the convergence speed.
 %  Maybe it can be explained by using different quantities to measure the convergence speed, but also some other frameworks. 
 % One possibility is to prove C-CD is faster than GD for some types of random matrices; but we suspect that   % Another example is compressive sensing
 
  \emph{How to Compare Algorithms?} It is widely accepted that Lanczos method is faster than power method, and conjugate gradient method is faster than GD, both theoretically and empirically.
 In particular, one theoretical justification is that in both cases the former achieves a rate dependent on $\sqrt{\kappa}$ while the latter achieves a rate dependent on $\kappa$.
 When it comes to the comparison of cyclic algorithms and randomized algorithms, the conclusion is far less clear. One issue is that there is no longer a proper metric like $\kappa$ to quantify the convergence rate of both algorithms. While $\kappa_{\mathrm CD}$ is a natural choice for 
 R-CD, the choices for C-CD are more abundant. In POCS literature, the rate is quantified by complicated functions of the angles between subspaces; in optimization literature, the rate is quantified by both $\kappa$ and $\kappa_{\mathrm CD}$, and sometimes complicated functions of the Hessian \cite{sun2015improved}. We also argue in Section \ref{sec: Jacobi Precond} that for non-equal-diagonal case, a natural metric should depend on eigenvalues of a Jacobi-preconditioned matrix, not the original matrix. It is for the comparison purpose that we express the convergence rate of C-CD in terms of the metric for R-CD. 
 Therefore, our work cannot provide a complete answer to the worst-case complexity of C-CD and Kaczmarz method, and the investigation on other quantities and the influence on the comparison is left as future work. 
 
 %Is randomization always an acceleration technique? 
 %The answer to this question is less clear than 
 %We heard some researchers say cyclic order is usually faster, and some say randomized order is usually faster. 
 %We want to understand whether the theoretical bounds are tight, and whether they make sense at all. 
 %A-GD is faster than GD since they have the same parameter.
 %C-CD is analyzed by kappa, but is that right? We need to compare with randomized versions.
 %There are many other bounds, but without benchmarks, we don't know the bounds are going to.
 %Build a theoretical framework to compare. Not just an example. Not sometimes A is faster, sometimes B is faster, but guidance. 

 \emph{Related Algorithms.} We hope this research will shed light on the study of related algorithms, such as POCS, SGD (Stochastic Gradient Descent) and ADMM  (Alternating Direction of Multiplier Method).
For ADMM, % the gap between cyclic version and randomly permuted version (not the independently randomized version) was recently discovered
it was recently found that the cyclic version with at least $3$ blocks can be divergent \cite{chen2016direct}, while randomly permuted version converges in expectation for solving linear systems \cite{sun2015expected}, so a fundamental gap between cyclic versions and randomly permuted versions exsits.
Nevertheless, it was also known that for certain problems (e.g. strongly convex) the small-stepsize versions of cyclic ADMM can be convergent \cite{hong2012linear,lin2014convergence,cai2014direct}.
Based on the results of the current work,  it is reasonable to conjecture that in these cases cyclic ADMM still achieve worse convergence rate than randomized versions of ADMM.
% ALM (Augmented Lagrangian Method, the ``batch version'' of ADMM) and 

\subsection{Notations and Organization}
Most notations in this paper are standard.
Throughout the paper, $A \in \R^{n \times n}$ is a symmetric positive semi-definite matrix.
Let $ L =\lambda_{\max}(A), \lambda_{\min}(A), \lambda_{\avg}(A)$ denote the maximum eigenvalue, minimum \emph{non-zero} eigenvalue and and average  eigenvalue of $A$ respectively;
sometimes we omit the argument $A$ and just use $\lambda_{\max}$, $\lambda_{\min}$ and $ \lambda_{\avg}$.
% It might be more clear to use $\lambda_{\min, \mathrm{nz}} $ to denote the minimum non-zero eigenvalue for the PSD case;
% but to save the notations we will just use the same notation $\lambda_{\min}$ as the positive definite case.
The condition number of $A$ is defined as $\kappa = \frac{\lambda_{\max}(A)}{\lambda_{\min}(A)}$.
Denote $A_{ij}$ as the $(i,j)$-th entry of $A$ and
 $L_i = A_{ii}$ as the $i$-th diagonal entry of $A$.
We use redundant notations $L $ and $L_i$ to be consistent with the optimization literature: $L$ represents the global Lipschitz constant
and $L_i$ represents the $i$-th coordinate Lipschitz constant of the gradient of the function $\frac{1}{2} x^T A x$.
We denote $\mathcal{R}(A) = \{  A x\mid  x \in \R^n \}$ as the range space of $A$.
Denote $ A^{\dag} $ as the pseudo-inverse of $A$, which can be defined
as $ V \text{diag}\{1/\lambda_1, \dots, 1/\lambda_r, 0, \dots, 0  \} V^{-1} $
when the eigen-decomposition of $A$ is $ V \text{diag}\{\lambda_1, \dots, \lambda_r, 0, \dots, 0  \} V^{-1} $,
where $\lambda_1, \dots, \lambda_r$ are all the non-zero eigenvalues of $A$.

The less widely used notations are summarized below.
We denote $L_{\max} = \max_i L_i$ and $L_{\min} = \min_i L_i$ as the maximum/minimum per-coordinate Lipschitz constant (i.e. maixmum/minimum diagonal entry of $A$),
and $L_{\avg} = (\sum_{i=1}^n L_i)/n$ as the average of the diagonal entries of $A$ (which is also the average
of the eigenvalues of $A$).
Denote $\kappa_{\mathrm CD} = \frac{L_{\avg} }{  \lambda_{\min} } = \frac{\lambda_{\avg}}{\lambda_{\min}}$ which  is a well-studied quantity that characterizes the convergence rate of R-CD.
% (to be more precise, the rate is determined by ).
We usually use $\Gamma$ to denote the lower triangular part of matrix $A$ with diagonal entries, i.e. $\Gamma_{ij} = A_{ij}$ iff $i \leq j$.
We also use $D_{A}$ to denote the diagonal matrix consisting of diagonal entries of $A$.
Finally, an important quantity $ \tau \triangleq \frac{L}{ L_{\min} } = \frac{\lambda_{\max}}{\lambda_{\avg}} $, a crucial ratio that characterizes the difference between GD, C-CD and R-CD.

The rest of the paper is organized as follows.
In Section \ref{sec: prelim}, we review the algorithms discussed in the paper.
In Section \ref{sec:BCPG}, we present our theoretical results on the complexity of C-CD as well as the comparison
of C-CD with other algorithms.
In Section \ref{sec: overview of proofs}, we provide an overview of the proof techniques and main steps. 
Section \ref{sec: proof of Main Result} is devoted to the proof of the main result Theorem \ref{thm: lower bound}.
In Section \ref{sec: simulation}, we present some numerical experiments.
In Section \ref{sec: conclulsion}, we summarize our findings and discuss some future directions. 
The proofs of results other than Theorem \ref{thm: lower bound} are provided in the appendix.

%  \emph{Lower Bound. } The complexity lower bound of general CD-type methods remains largely open (we only show the lower bound of one particular algorithm).
%  A recent lower bound on minimizing a finite sum (including regularized least square problems,
%  but not the least square problem itself) \cite{agarwal2014lower} only applies to deterministic algorithms,
%   while no deterministic algorithms are known to achieve that lower bound.
%     Our work indicates that there is probably a large theoretical gap between cyclic algorithms and randomized algorithms.
   % Randomized algorithms can achieve this lower bound, thus an interesting question is whether there exists deterministic algorithms that can achieve the same complexity as randomized algorithms, closing the gap between the upper bound and the lower bound.

%   Many works on cyclic multi-block ADMM (e.g. \cite{hong2012linear}) are based on the analysis of cyclic BCD.
%   Our work implies that if the analysis for cyclic multi-block ADMM is based on the analysis of cyclic BCD, then
%   the obtained complexity bounds are worse than one might expect.
%  to work, it is highly likely that the stepsize should be of the order $O(1/n)$.

%Some recent fast solvers for Laplacian systems and LP rely on randomized coordinate descent methods.
%   A similar question is whether the deterministic computational complexity of these problems are the same as randomized complexity.
   %  We expect BCD to be used as sub-routines in many other algorithms.

\section{Preliminaries: Several Algorithms }\label{sec: prelim}

In this section we will review several variants of CD, G-S method, Kaczmarz method and POCS. 
% For simplicity, we present the versions of CD for minimizing convex quadratic functions instead of  the most general form.
We mainly consider the quadratic minimization problem
$$
  \min_{x \in \R^n }  f(x) \triangleq \frac{1}{2} x^T A x - b^T x,
$$
where $A \in \R^{n \times n}$ is a symmetric PSD (positive semi-definite) matrix, $b \in \mathcal{R}(A)$ and $A_{ii} \neq 0, \ \forall i$.
All the optimal solutions of the problem satisfy the system of linear equations
$$
  A x = b.
$$
When $A$ is non-singular (thus positive definite), the unique minimizer $x = A^{-1} b $ is the unique solution to the linear system.
When $A$ is singular, there are infinitely many optimal solutions. 
%  of the form $A^{\dag} b + x^{\bot} $ where $A x^{\bot} = 0$.
% Most results in this paper (except the convergence of iterates) can be easily extended to the non-strongly convex case (i.e. PSD linear systems) by restricting to the range space of $A$,  but we will not discuss such extensions for simplicity.
% While the direct method for computing $A^{-1} $ takes $\O(n^3)$ time, iterative algorithms are definitely more favorable
% for many large-scale problems.

\textbf{Gradient descent.} 
GD (gradient descent) is one of the most basic iterative algorithms. Starting at $x^0 \in \R^n$, GD proceeds as follows:
$$
  x^{k+1} = x^k - \frac{1}{L} \nabla f(x)
      = x^k - \frac{1}{L}(Ax - b).
$$
There are many other choices of stepsizes, but we use a constant stepsize $1/L$ in the paper because it is simple and already leads to the standard complexity $\O(n^2 \kappa \log 1/\epsilon  )$ for quadratic problems.
 % it is somewhat standard in optimization literature; 2) this stepsize already leads to the well known comple

\textbf{Cyclic Coordinate Descent and Gauss-Seidel Method. } 
The C-CD algorithm updates the variables cyclically by minimizing the objective function over one variable with other variables fixed. Each cycle of C-CD consists of the update of all variables: 
 $$ x_i \leftarrow  \argmin_{x_i} f(x_i; x_{-i}), \;, i=1,\dots, n ,$$ 
 where $x_{-i}$ denotes  the collection of all variables except $x_i$.
The update order in each cycle is fixed, such as $(12\dots n)$.
For the quadratic problem, the subproblems are single-variable quadratic problems with closed-form solutions. Thus it can be written in the following way, assuming the initial point is $x^0 = x^{0,0}$:
\begin{equation}\label{cyclic CD update}
\begin{split}
 & x^{k, j} = x^{k,j-1} - \frac{ A(j,:)x^{k,j-1} - b_j }{ A_{jj} } e_j,  \quad j=1,2, \dots, n ; \\
 & x^{k+1} = x^{k,n}, \  x^{k+1, 0} = x^{k+1}.
\end{split}
\end{equation}
where $e_j$ is the $j$-th standard unit vector with only one nonzero entry $1$ in the $j$-th position,  $A(j,:)$ denotes the $j$-th row of $A$, and $A_{jj}$ denotes the $j$-th diagonal entry of $A$. 

The algorithm \eqref{cyclic CD update} is also the Gauss-Seidel method for solving the linear system $Ax = b$.
Note that even if $A$ is not symmetric, one can still apply G-S method (the update equations are exactly the same as above), but it only converges under certain assumptions on $A$. 

% For our example, $A(j,j)=1, \forall j$.
We can write the above update equation as a simple matrix recursion  % When using C-CD to solve the problem,
$$
   x^{k+1} - x^* = (I - \Gamma^{-1}A ) ( x^k - x^*),
$$
where $x^*$ is one optimal solution, and $\Gamma$ is the lower triangular part of $A$ with diagonal entries, i.e.,
 $\Gamma_{ij} = A_{ij}, 1 \leq j \leq i \leq n $.
%\begin{equation}\label{L def}
%   \Gamma = \begin{bmatrix}
%    1          &  0      &  \dots  &  0   \\
%    c          &  1      &  \dots  &  0   \\
%    \vdots     &  \vdots  &   \ddots &  \vdots \\
%    c          &   c     &   \dots &   1
% \end{bmatrix}
%\end{equation}
 We denote the iteration matrix as  %The convergence speed of CD is related to the spectral radius of the update matrix
$$
  M = I - \Gamma^{-1}A.
$$

\textbf{Randomized Coordinate Descent.}
R-CD (randomized coordinate descent) algorithm starts at $ z^0 $ and proceeds as follows:
\begin{equation}\label{R-CD}
\begin{split}
  & \text{FOR } j=1,2, \dots  \\
 & \quad \quad \text{Randomly pick } t \in \{1,\dots, n \} \text{ uniformly at random}, \\
 &  \quad \quad z^{ j+1 } = z^{j} - \frac{ A(t,:)z^{ j } - b_j }{ A(t,t) } e_t .  \\
% & \text{END }  \\
\end{split}
\end{equation}
The output of R-CD is a sequence $(z^1, z^2, \dots)$. We further define
$
  x^{k } = z^{kn}, \; \forall \; k,
$
to be comparable with GD and C-CD. Here $k$ can be viewed as the index of ``epoch'', where each epoch consists of $n$ iterations.

\textbf{Randomly Permuted Coordinate Descent.}
We also consider RP-CD (randomly permuted coordinate descent). Starting from $x^0 = x^{0,0}$, the algorithm proceeds as follows.
\begin{equation}\label{RCD update equation}
\begin{split}
  & \text{At epoch } k, \text{ pick a permutation } \sigma_k \text{ uniformly at random from the set of all permutations}.  \\
  & \text{FOR } j= 1,2, \dots, n,  \\
  &  \quad \quad  t = \sigma_k(j) ; \\
  &  \quad \quad x^{ k, j } = x^{k,j-1} - \frac{ A(t,:) x^{k,j-1} - b_j }{ A_{tt} } e_t ;  \\
  &   x^{k+1} = x^{k,n}, \  x^{k+1, 0} = x^{k+1}.
 % &  \text{END }  \\
\end{split}
\end{equation}

 According to \cite[Section II.A]{sun2015expected}, the recursion formula of RP-CD is
 \begin{equation}\label{RP CD update equation}
   x^{k} = (I - \Gamma_{\sigma_k}^{-1} A ) x^{k-1},
 \end{equation}
 where $\sigma_k$ is the permutation used in the $k$-th iteration, and $\Gamma_{\sigma} \in \R^{n \times n}$ is defined by
 \begin{equation}\label{Lsigma def}
   \Gamma_{\sigma}( \sigma(i) , \sigma(j) ) \triangleq \begin{cases}
  A_{\sigma(i), \sigma(j)} &  j \leq i .    \\
0  &  j > i ,
\end{cases}
 \end{equation}
 For example, when $n=3$ and $\sigma = (\sigma(1), \sigma(2), \sigma(3)) = (231)$,
 $$
  \Gamma_{(231)} =
   \begin{bmatrix}
   A_{11} &    A_{12}   &     A_{13}        \\
   0 &  A_{22}  &      0      \\
   0 &  A_{32} &   A_{33}     \\
   \end{bmatrix}.
$$

\textbf{POCS \cite{von1949rings,halperin1962product, bauschke1997method,escalante2011alternating,galantai2013projectors}.} 
POCS
 is a general method to find a common point of $m $ closed convex sets $\mathcal{M}_1 , \dots, \mathcal{M}_m$. Starting from any point $x^0$, the algorithm proceeds by performing projection onto these sets one by one:
$$
x^{k+1} =  P_m P_{m-1} \dots P_1 x^k,
$$
where $ P_j z = \text{Proj}_{\mathcal{M}_j} (z) $ is the projection of $z$ onto the set $\mathcal{M}_j$. 
%   \text{Proj}_{\mathcal{M}_{m-1} } \dots    \text{Proj}_{\mathcal{M}_{1} }

\textbf{Kaczmarz Method \cite{kaczmarz1937angenaherte}.}
Consider a  linear system of equations $ U y =b $, where $U \in \R^{ n \times m }, y \in \R^{m \times 1}, b \in \R^{n \times 1}$ and $n \geq m$. 
Suppose $U^T = (u_1, u_2, \dots, u_n)  $ and $b^T = (b_1, \dots, b_n)$, 
then a solution of $ U y = b $ is a point in the intersection of $n $ hyperplanes
$\mathcal{H}_k = \{ y \mid \langle u_k,  y \rangle = b_k \},k=1,\dots, n. $ 
Kaczmarz method is a special case of POCS for finding the intersection of hyperplanes. More specifically, starting from an arbitrary initial point $y^{0,0}$, the algorithm 
proceeds as
\begin{equation}\label{Kaczmarz update}
\begin{split}
& y^{k, j} = \text{Proj}_{\mathcal{H}_j} (y^{k,j-1}) = y^{k,j-1} + \frac{  b_j - \langle u_j, y^{k,j-1}  \rangle }{\|u_j \|^2 }  u_j,     \quad j=1,2, \dots, n ; \\
& y^{k+1} = y^{k,n}, \  y^{k+1, 0} = y^{k+1}.
\end{split}
\end{equation}

\textbf{Connections between Different Methods.}
As mentioned above, G-S method for solving a symmetric PSD linear system of equations is a special case of C-CD. Kaczmarz method is a special case of POCS. The following claim shows that G-S method for solving a symmetric PD linear system $UU^T x = b$  is equivalent to Kaczmarz method for a full-rank square system $Uy =b$.

\begin{claim}\label{equivalence of Kaczmarz and G-S}
	Suppose $b \in  \R^{ n \times 1}$, $A = UU^T \in \R^{n \times n} $, where $U \in \R^{n \times n}$ is full rank. Then  Gauss-Seidel method for solving $A x = b$ is equivalent to Kaczmarz method
	for solving $U y = b$; here, the equivalence means that there is a one-to-one mapping between
	the iterates of the two algorithms. 
\end{claim}

The proof is given in Appendix \ref{appen: Proof of equivalence of GS and Kaczmarz}.
Intuitively, under a coordinate transformation Kaczmarz method is equivalent to G-S method. More specifically,
any vector $y $ can be expressed as $y = x_1 u_1 + \dots + x_n u_n = U^T x$, i.e., $x_j$'s are the coordinates of $y$ under the basis $u_1, \dots, u_n$, where $u_j$'s are columns of $U^T$. Thus, updating one coordinate $x_j$ is equivalent to updating $y$ according to one equation $\langle u_j, y \rangle = b_j$.

When $U$ is not square and/or not full rank, as long as the initial point of Kaczmarz method lies in the row space of $U$, we can still show the almost ``equivalence'' of Kaczmarz method and G-S method, though there is no one-to-one mapping but a one-to-many mapping; see Appendix \ref{appen: Proof of equivalence of GS and Kaczmarz}. 
Therefore, in the basic setting, all four methods C-CD, G-S, Kaczmarz and POCS are equivalent. 

\section{  Main Results }\label{sec:BCPG}
%The linear convergence of CD method for solving quadratic minimization problem (even non-strongly convex case) was established
%in \cite{luo1992convergence}. The rate of convergence depends on the Hoffman bound, and this rate was explicitly given in Section 3.4 of \cite{wang2014iteration} (after eq. (28)). We restate this result below.

Consider the quadratic minimization problem  $$ \min_x f(x) \triangleq x^TA x - 2b^T x,  $$ where $A \in \R^{n \times n} $ is symmetric positive semi-definite, $b \in \mathcal{R}(A)$  and $ A_{ii} \neq 0, \ \forall i$.
We can assume $b \in \mathcal{R}(A) $ since otherwise the minimum value of $\min_x x^TA x - 2b^T x$ will be $- \infty$. 
We can assume $A_{ii} \neq 0, \ \forall i$, since when some $A_{ii} = 0 $ all entries in the $i$-th row and the $i$-th column of $A$ should be zero, which means that the $i$-th variable does not affect the objective and thus can be deleted. 
Recall that the maximum eigenvalue, minimum eigenvalue and average eigenvalue of $A$ are
$\lambda_{\max},\lambda_{\min},\lambda_{\avg}$ respectively,
the condition number $\kappa = \lambda_{\max}/\lambda_{\min} $,  and
$\kappa_{\mathrm CD} = \lambda_{\avg}/\lambda_{\min} $. 

To help the readers understand our main results, we first summarize the main results in the following Table \ref{table of complexity} for the equal-diagonal case (i.e. all diagonal entries of $A$ are equal). The upper bounds will be given in Proposition \ref{prop: upper bound},  and the lower bounds will be given in Theorem  \ref{thm: lower bound}. In this table, we ignore the  $\log 1/\epsilon$ factor, which is always necessary for an iterative algorithm to achieve error $\epsilon$. 
\begin{table}[!htbp]
	\centering
	\caption{Complexity of C-CD, GD and R-CD for equal-diagonal case (ignoring a $\log 1/\epsilon$ factor)} \label{table of complexity}
 \begin{tabular}{|c|c|c|c|}
 	\hline 
 	Parameters	&   $\kappa$ &  $\kappa$ and $\tau$ &  $\kappa_{\mathrm CD}  ( = \kappa/\tau) $    \\ 
 	\hline 
 	C-CD Upper bound	&   $   n^3 \kappa  $  &
 	$   \frac{1}{10} n^2  \kappa \tau \log^2 n  $
 	&  $   n^4 \kappa_{\mathrm CD} $    \\ 
 	\hline 
 	C-CD Lower bound  (Theorem \ref{thm: lower bound})	&   $ \frac{1}{40 }   n^3 \kappa $  &
 	$ \frac{1}{40 } n^2  \kappa \tau   $
 	&  $  \frac{1}{40} n^4 \kappa_{\mathrm CD} $    \\ 
 	\hline 
 	GD  &  \multicolumn{2}{|c|}{$  n^2 \kappa  $ }  & --   \\
 	\hline 
 		R-CD   &  -- &  $n^2 \kappa /\tau$  &  $ n^2 \kappa_{\mathrm CD}  $ \\
 		\hline
 \end{tabular} 
\end{table}

This table shows that the lower bounds match the upper bounds, up to constant and $\log^2 n$ factors.  In addition, the table reveals the relations between the worst-case complexity of C-CD, GD and R-CD. To make the relationships easy to read, we extract the  results on parameters $\kappa $ and $\tau$ (i.e. the middle column)  and normalize them by the complexity of GD to create Table \ref{table simple compare}.

\begin{table}
	\centering	
	\caption{ Complexity for equal-diagonal case (divided by $n^2 \kappa \log \frac{1}{\epsilon}$ and ignoring constants. $\tau = \lambda_{\max}/\lambda_{\avg} \in [1, n]$ ) }	\label{table simple compare}
	\begin{tabular}{|c|c|l|c|}
		\hline
		&       C-CD          &              GD         &   R-CD \\
		\hline
		Lower bound       &     $  \tau $             &           $  1 $       &    --  \\
		\hline
		Upper bound       &     $ \min\{   \tau  \log^2 n ,  n \} $    &           $ 1  $       &   $  \frac{1}{\tau}  $   \\
		\hline
	\end{tabular}
\end{table}
% $   \frac{1}{\tau}  $     

According to the tables, the main implications of our results are the following:
\begin{itemize}
 \item  C-CD is roughly $\O(\tau)$ times slower than GD, and R-CD is $\O(\tau)$ times faster than GD.	
  \item When $\tau $ achieves the maximum $\O(n)$, C-CD is $\O(n)$ times slower than GD and R-CD is $\O(n)$ times faster than GD.	This implies C-CD can be $\O(n^2)$ times slower than R-CD. 
  \end{itemize}
Note that in the above statement ``method 1 is X-times slower than method 2'' does not mean that method 1 is always slower than method 2 (of course
rarely can one make such a strong statement), but that ``the worst-case complexity of method 1 is X-times worse than that of method 2, and both complexity bounds can be simultaneously achieved''.

%We say ``roughly'' when comparing C-CD with GD because of the $\log^2 n$ gap between the upper bound and the lower bound; see more detailed discussion
%on this gap in the following subsection.

% \end{itemize}

% & $\kappa_{\mathrm CD} $ and $\tau$  
%  $  n^2  \kappa_{\mathrm CD} \tau^2 \log^2 n   $
% & $ n^2  \kappa_{\mathrm CD} \tau^2    $
%  \multicolumn{2}{|c|}{ $ n^2 \kappa_{\mathrm CD}  $ }  

%%--------------------------------------------------------------------
%%-------------   a table with O 
%%--------------------------------------------------------------------
%\begin{table}[!htbp]
%	\centering
%	\caption{Complexity of C-CD}
%	\begin{tabular}{|c|c|c|c|c|}
%		\hline 
%		Parameters	&   $\kappa$ &  $\kappa$ and $\tau$ &  $\kappa_{\mathrm CD} $ & $\kappa_{\mathrm CD} $ and $\tau$  \\ 
%		\hline 
%		Upper bound	&   $  \tO\left( n^3 \kappa \right) $  &
%		$\tO\left( n^2  \kappa \tau \log^2 n  \right)$
%		&  $   \tO\left( n^4 \kappa_{\mathrm CD} \right) $  & $	\tO\left( n^2  \kappa_{\mathrm CD} \tau^2 \log^2 n   \right) $ \\ 
%		\hline 
%		Lower bound	&   $  \tO\left( n^3 \kappa \right) $  &
%		$\tO\left( n^2  \kappa \tau  \right)$
%		&  $   \tO\left( n^4 \kappa_{\mathrm CD} \right) $  & $	\tO\left( n^2  \kappa_{\mathrm CD} \tau^2   \right) $  \\ 
%		\hline 
%		GD  &  \multicolumn{2}{|c|}{$ \tO\left( n^2 \kappa \right)  $ }  &  -- & --  \\
%		\hline 
%		R-CD   & - & - &  \multicolumn{2}{|c|}{$ \tO\left( n^2 \kappa_{\mathrm CD} \right)  $ }   \\
%		\hline
%	\end{tabular} 
%\end{table}
%--------------------------------------------------------------------
 
Now we formally state the upper bounds and lower bounds on the convergence rate of C-CD. 
\begin{proposition}\label{prop: upper bound}(Upper bound of C-CD)
Consider the quadratic minimization problem  $\min_x f(x) \triangleq x^TA x - 2b^T x $ where $A \in \R^{n \times n} $ is positive semi-definite, $b \in \mathcal{R}(A)$  and $ A_{ii} \neq 0, \ \forall i$.
   For any $x^0 \in \R^n $, let $x^k$ denotes the output of C-CD after $k$ cycles, then
\begin{subequations}
\begin{align}
     f(x^{k+1}) - f^* \leq \min \left\{ 1 - \frac{1}{ n \kappa  } \frac{L_{\min}}{L_{\avg}} ,  1 - \frac{L_{\min} }{ L (2 + \log n/ \pi)^2 } \frac{1}{ \kappa  } \right\}  (f(x^k) - f^*).  \label{upper bound in kappa} \\
     f(x^{k+1}) - f^* \leq \min \left\{ 1 - \frac{1}{ n^2 \kappa_{\mathrm CD}  } \frac{L_{\min}}{L_{\avg}} ,
   1 - \frac{ L_{\min} L_{\avg} }{ L^2 (2 + \log n/ \pi)^2 } \frac{1}{ \kappa_{\mathrm CD}  } \right\}  (f(x^k) - f^*). \label{upper bound in kappaCD}
\end{align}
 \end{subequations}
%We also have
%\begin{equation}
%
%\end{equation}
Here, $f^*$ is the minimum value of the function $f$,   %$x^*$ is the minimum of function $f$, $f^* = f(x^*)$,

\end{proposition}

%It is very tempting to think that the huge gap between C-CD and GD/R-CD is just due to the weakness of the proof techniques.
%In particular, it seems very strange that C-CD converges much slower than GD, since people expect that it should be faster than GD.
%Rather surprisingly, we will show that the bounds given in Proposition \ref{prop: upper bound} is actually tight
%when all diagonal entries are the same (possibly up to $\log^2 n$ factor).
% up to a constant factor,

\begin{theorem}\label{thm: lower bound}(Lower bound of C-CD)
	For any initial point $x^0 \in \R^n $, any $\delta \in (0,1]$, there exists a quadratic function $f(x) = x^TA x - 2b^T x $ such that % the sequence $\{x^k\}$ generated by cyclic coordinate descent algorithm
	%  satisfies
	\begin{subequations} \label{desired bound}
		\begin{align}
		f(x^k) - f^* \geq (1 - \delta) \left( 1 - \frac{ 2\pi^2 }{ n \kappa }  \right)^{2k + 2 } ( f(x^0) - f^*) , \;  \forall k,        \label{desired bound: kappa} \\
		f(x^k) - f^* \geq (1 - \delta) \left( 1 - \frac{ 2\pi^2 }{ n^2 \kappa_{\mathrm CD} }  \right)^{2k + 2 } ( f(x^0) - f^*)       , \;  \forall k,   \label{desired bound: kappaCD}
		\end{align}
	\end{subequations}
	% Corollary:  f(x^k) - f^* \geq (1 - \delta) \mu \left( 1 - \frac{ 2\pi^2 }{ n \kappa }  \right)^{2k + 2} \| x^0 - x^* \|^2,   \\ \\
	where  $x^k$ denotes the output of C-CD after $k$ cycles,
	$f^* $  is the minimum of the objective function $f$.
\end{theorem}
% , $\kappa = \frac{\lambda_{\max}(A)}{\lambda_{\min}(A)}$ is the condition number,
% $\kappa_{\mathrm CD} =  \frac{ L_{\avg} }{ \lambda_{\min}(A) } = \frac{\sum_i A_{ii}}{n \lambda_{\min}(A)} . $

The overview of the proofs will be given in Section \ref{sec: overview of proofs}.
The formal proof of Proposition \ref{prop: upper bound} will be given in Appendix \ref{appen: proof of upper bound}, and the formal proof of Theorem \ref{thm: lower bound}  will be given in Section \ref{sec: proof of Main Result}.

Remark 1: The example we construct is simple: all diagonal entries of $A$ are $1$ and all off-diagonal entries are $c  $, where $c$ is a constant close to $1$.
It is known that the SDD (symmetric diagonally dominant) system can be solved in almost linear time (see, e.g., \cite{spielman2004nearly,spielman2014nearly,Peng14} and the references therein).
While for SDD system the off-diagonal entries have very small magnitude,
the hard instance we construct can be viewed as the opposite of the SDD system: the off-diagonal entries are chosen as large as possible
so that the matrix remains positive definite.

Remark 2: Throughout the paper, our discussion focuses the comparison of the \textit{total time complexity},
instead of the \textit{iteration complexity}. 
For quadratic problems, the two are closely related because each epoch of C-CD, GD and R-CD (under the proper definition of ``epoch'') takes approximately the same time. 
For general convex problems, our lower bound result Theorem \ref{thm: lower bound} should be viewed as a lower bound on the iteration complexity of C-CD.

%\begin{corollary}\label{prop: upper bound, in kappaCD}(Another upper bound of C-CD)
%\end{corollary}  L_{\max} = \max_i A_{ii},
%Remark:
% It is easy to prove $\| \Gamma^T A^{-1} \Gamma  \|$
% f(x^k) - f^* \leq  \left( 1 - \frac{ 1 }{ 4(6 + 5 \sqrt{n} + n )\kappa + (7 + 2 \sqrt{n}) }  \right)^{k } ( f(x^0) - f^*)
% f(x^k) - f^* \leq  \left( 1 - \frac{ 1 }{ 4(6 + 5 \sqrt{n} + n )n \kappa_{\mathrm CD} + (7 + 2 \sqrt{n}) }  \right)^{k } ( f(x^0) - f^*)

We then describe how to obtain Table \ref{table of complexity} from the two results. 
 As mentioned in the introduction,  we will say an algorithm has complexity $\tO( g(n, \theta) )$,
if it takes $\O( g(n, \theta)  \log(1/\epsilon) )$ unit operations 
% (plus, minus, multiplication and division)
to achieve relative error $\epsilon$.
Each iteration of GD, each epoch (i.e. $n$ iterations) of C-CD and R-CD  all take $\O(n^2)$ operations \footnote{When the matrix is sparse, the time is actually $O(\text{nnz}(A))$, but to simplify the discussions, we do not consider the sparsity in this work.}. 
 Using the fact $ -\ln(1 - z) \geq - z, z \in (0,1) $ one can immediately show that to achieve
$ (1 - 1/u )^k \leq \epsilon $ one only needs $ k \geq u \log(1/\epsilon) $ epochs. 
Thus we can transform the convergence rate to the   number of epochs, then the complexity \footnote{To be precise, the upper bounds on the convergence rate can be transformed to upper bounds of the complexity, but the lower bounds require a bit of more work. We can make it precise, but let us ignore this minor issue, and just assume both upper bounds and lower bounds of convergence rate can be transformed to corresponding complexity bounds.  }

Consider the equal-diagonal case (i.e. $L_i = L_1, \ \forall i$) for now and we will discuss the general case later in Section \ref{sec: Jacobi Precond}.
In this case, $ L_{\avg} = L_{\min}   $, thus greatly simplifying the bounds; further,
$ \frac{L}{ L_{\min} } $ is just the quantity  $\tau = \frac{\lambda_{\max}}{\lambda_{\avg}} $.
The upper bounds on convergence rate \eqref{upper bound in kappa}  can be transformed to the following upper bound of complexity
\begin{equation}\label{upper bound of complexity, kappa}
 \min \left\{ \tO\left( n^3 \kappa \right),  \tO\left( n^2  \kappa  \tau \log^2 n \right) \right\}.
  \end{equation}
  These two quantities are those in the first two entries of C-CD upper bound in Table \ref{table of complexity}. Similarly, the other bounds on convergence rate in Proposition \ref{prop: upper bound} and Theorem \ref{thm: lower bound} can be transformed to corresponding bounds
  on the complexity, and they form the rest of Table \ref{table of complexity}.

 \section{Extensions and Discussions}
 \subsection{Comparison with Known Convergence Rate of POCS}
 
The convergence rate of POCS for finding the intersection of closed subspaces of a real Hilbert space has been stuided since 1970s. One of the first convergence rate results is given by Smith, Solmon and Wagner \cite{smith1977practical}, and cited as a major convergence rate result of POCS in
\cite{galantai2013projectors, escalante2011alternating}.
Further results are given in Kayalar and Weinert \cite{kayalar1988error} and Deutsch and Hundal \cite{deutsch1997rate}, but these rates are very complicated.
% : the former bound is the minimum of four terms, and the latter is the minimum of $n$ terms. 
% It was shown in \cite{deutsch1997rate} that none of the three bounds are sharp for general convex sets. 

% UNNECESSARY Paragraph
%It is interesting to compare the convergence rate of POCS with that of ours. 
%A major difference between the results for POCS and the results for CD lies on the parameters. The rates for POCS depend on the angles between intersections of the desired subspaces, while the rates for CD depend on the Lipschitz constants of the partial gradients and sometimes strong convexity paramters. In general, it is not clear how to compare these two kinds of convergence rate. 
%Neverthelss, for the simple setting of solving linear systems studied in this work, 
%such a comparison is possible. 

Due to the complication of the bounds of \cite{kayalar1988error} and \cite{deutsch1997rate},
we will only consider the  classical convergence rate in \cite{smith1977practical}. 
The original result characterizes the rate by the angles between subspaces; interestingly,
for the simple case of solving a  linear system of equations $Uy = b$,
the convergence rate can be charaterized by the determinant of the matrix $UU^T$.
For simplicity, we present the result for a full-rank square linear system. 

\begin{proposition}\cite{smith1977practical, galantai2005rate}\label{prop: POCS rate}
	Consider the linear sytem of equations $U y = b$,
	where $U^T = (u_1, \dots, u_n)$ is an $n \times n $ matrix with full rank and $\|u_j\| = 1, \forall j$. Suppose the  sequence generated by Kaczmarz method is 	$\{  y^k \}$, then  
	$$
	  \| y^k - y^* \| \leq (1 - \det(UU^T)  )^{k/2} \| y^0 - y^* \|. 
	$$  
\end{proposition} 

This rate can be transformed to a rate dependent on the eigenvalues by using the fact $\det(UU^T) = \lambda_1 \dots \lambda_n $, where  $\lambda_1 \geq  \dots \geq \lambda_n$ are  the eigenvalues of $A = UU^T$.  The number of epochs to achieve a relative error $\epsilon$ predicted by the above result is
$$
 2 \frac{1}{ \lambda_1 \lambda_2 \dots \lambda_n } \log \frac{1}{\epsilon}.
$$
The number of epochs predicted by our result is approximately
$$
 \frac{ \lambda_{\max}^2 }{ \lambda_{\min}} \log \frac{1}{\epsilon} = \frac{\lambda_1^2}{ \lambda_n } \log \frac{1}{\epsilon},
$$
in which we ignore the constant factor and $\log n $ factor. 
To simplify the comparison, let us denote 
$$ T_{\text{POCS}} \triangleq \frac{1}{ \lambda_1 \lambda_2 \dots \lambda_n } , \quad 
 T_{\text{C-CD}} \triangleq \frac{\lambda_1^2}{ \lambda_n } . $$
 
For the example that achieves the lower bound (see Example 2 of Section \ref{sec: role of examples}), the eigenvalues are 
$$
 \lambda_1 = 1 - c + cn, \; \lambda_2 = \lambda_3 = \dots  = \lambda_n = 1 -c,
$$
where $0 < c < 1$. 
Then $ T_{\text{POCS}} = \frac{1}{ (1 - c + cn) (1 - c)^{n-1} }$, and $ T_{\text{C-CD}} 
=  \frac{ (1 - c + cn)^2}{ 1 -c }  . $ 
The ratio of the two quantities are
$$
    \frac{ T_{\text{C-CD}}  } {  T_{\text{POCS}}  }=   (1-c)^{n - 2}  (1 -c + cn) ^3  \rightarrow 0, \quad \text{ as }  c \rightarrow 1.
$$
Thus the POCS bound, as given by Proposition \ref{prop: POCS rate}, is very loose for our example, and can be infinitely times worse than our bound. It is easy to show that as long as $ c > 1 - K^{1/(n-2)} n^{-3/(n-2)} \approx 0 $, the above ratio is less than $1/K$, meaning that the POCS bound is $K$ times worse than our bound (up to a $\log(n)$ factor).

In general, we can show that $  T_{\text{C-CD}} \leq 10 T_{\text{POCS}}  $. 
We  need the condition $\| u_j\|^2 = 1  $, which means $A_{jj} = 1$ and thus
$\sum_{i=1}^n \lambda_i = \trace(A) = n $. By algebraic-mean-geometric-mean inequality, we have
% where $K$ captures the missing constant factor and $\log n$ factor. I
$$
   \frac{  T_{\text{C-CD}}  }{T_{\text{POCS}} } =
  \lambda_1^3 \lambda_2 \dots \lambda_{n-1} 
  = 27 (\frac{ 1 }{3} \lambda_1 )^3    \lambda_2 \dots \lambda_{n-1}
  \leq 27 ( \frac{  \lambda_1 + \lambda_2 + \dots + \lambda_{n - 1}   }{ n + 1 } )^{n+1}
   \leq 27  (\frac{ n}{ n+1} )^{n+1}  \leq 27/e  \leq 10, %\frac{1}{3}. 
$$
This relation means that our bound is at least as good as the POCS bound (up to constant and $\log(n)$ factors). 

This comparison has a few implications. First, an interesting question is whether  the bound in this paper can be improved by using other metrics like the angles between subspaces. From the comparison we know that at least the classical bound of POCS does not provide the improvement. Second, as the classical bound of POCS can be infinitely times worse than our bound, there is large room of improvement for general POCS. 
%  as our complexity bound is already quite large, the classical POCS bound in \cite{smith1977practical}
% is even worse than our bound, thus
% % : we compare our bound of cyclic POCS with the classical bound of cyclic POCS in \cite{smith1977practical}. 

In this subsection we only consider the deterministic cyclic version of POCS.
Randomized versions of POCS (see, e.g., \cite{strohmer2009randomized,necoara2018randomized}) can have much faster convergence rates.  In fact, since randomized POCS has the same rate as R-CD for solving linear systems \cite{strohmer2009randomized}, we infer that the bound of
randomized POCS can be infinitely times better than the classical rate of cyclic POCS in \cite{smith1977practical} \footnote{As a historical remark, the paper \cite{strohmer2009randomized} could have used the comparison of their rate with the bound of cyclic POCS in \cite{strohmer2009randomized} to justify their proposal of randomized Kaczmarz method.}. By using our bound,
cyclic POCS is up to $O(n^2)$ times worse than randomized POCS.
% , and can be infinitely times better. nd R-CD can be $O(n^2)$-times faster than C-CD,

% Recall that one main result of our paper is that there is an $O(n^2)$ gap between our bound of cyclic POCS and the bound of randomized POCS for solving linear systems. Combining this result with the comparison done in this subsection, we infer that there is at least an $O(n^2)$ gap between the classical bound of cyclic POCS \cite{smith1977practical} and the bound of randomized POCS. 
% Due to the lack of communication between the researchers in POCS and the researchers in optimization, this comparison is probably unknown to both fields before.  
%n fact, the convergence rate of randomized Kaczmarz method for linear systems \cite{strohmer2009randomized} 
% Randomized versions of POCS (see, e.g., \cite{strohmer2009randomized,necoara2018randomized}) can have much faster convergence rates, similar to randomized versions of CD.  	

% (a special case is randomized Kaczmarz method \cite{strohmer2009randomized}) 
%-------------------------------------------------------------------------------------------
%Estimating the spectral norm of each order and then taking the average is not going to work.
%This is because each fixed order is slow. We have to analyze the spectral norm of the average,
%not the average of the spectral norm. 
%-------------------------------------------------------------------------------------------

\subsection{Role of Examples in Convergence Rate Analysis}\label{sec: role of examples}
We propose the following example to show the lower bound of the convergence rate of C-CD. 
The rigorous analysis of this example is long and technical, and will be provided in a later section.

\textbf{Example 1:}
For any constant $c \in (0,1)$, consider minimizing the following quadratic function
\begin{equation}\label{quadratic min}
\min_{x \in \R^n } f(x) \triangleq x^T A_c x ,
\end{equation}
where $ A_c \in \R^{n \times n }$ is defined as
% , we pick $A = A_c$ as follows:
\begin{equation}\label{Ac def}
A_c = \begin{bmatrix}
1          &  c      &  \dots  &  c   \\
c          &  1      &  \dots  &  c   \\
\vdots      &  \vdots  &   \ddots &  \vdots \\
c          &   c     &   \dots &   1
\end{bmatrix}
\end{equation}
Simple calculation shows that $A_c$ is a positive definite matrix.
% with one eigenvalue $1-c$ with multiplicity $n-1$ and one
%eigenvalue $1-c + cn $ with multiplicity $1$. 
%The optimum of the problem is $x = (0; 0; \dots ; 0) $.
%Solving this problem is also equivalent to solving a linear system of equations $Ax = 0$.

After posting the first version of the paper in April 2016, Steven Wright pointed out to us that he proposed the matrix that we analyzed in this paper in a talk in Paris in July 2015  and in a talk at NYU in December 2015.  He also noticed the large gap between C-CD and randomized CD for this example, although no theoretical analysis was available on public. 

Another example was brought to our attention independently by Strohmer and Richtarik after posting the first version of this paper. This example shows that cyclic Kaczmarz method  can be much slower than the randomized version. 

\textbf{Example 2:} Consider solving a linear system of equations $U y = 0$ where
$U^T = (u_1, \dots, u_n) \in \R^{2 \times n} $, and $u_j = ( \cos(\theta_k), \sin(\theta_k) ), k=1,...,n,$ where $\theta_k = 2 k \pi/n. $ The hyperplanes $\mathcal{H}_j = \{y \mid \langle u_j , y \rangle = 0 \}, j=1,\dots, n$
are $n$ lines crossing the origin with angles $2 \pi/n$ between two adjacent lines.
Cyclic projection to the lines $H_1, H_2, \dots, H_n$ one by one can be very slow, 
and randomized projection is much faster. 

We have checked this example by simulations. 
Since Kaczmarz method can be transformed to Gauss-Seidel method or equivalently C-CD,
instead of solving $Uy = 0$ by Kaczmarz method, we consider solving $UU^T x = 0$ by C-CD,GD and R-CD.  We have some interesting findings:

1) C-CD is slow if we update the coordinates in the order $(12\dots n)$. If we pick a random order and use this order throughout, C-CD is actually very fast. 
In contrast, for Example 1, any fixed order is slow (similarly, for the divergent examples of cyclic ADMM in \cite{chen2016direct} and \cite{sun2015expected}, any fixed order is divergent). In this sense, Example 2 is a ``weak'' bad examlpe for C-CD, and Example 1 is a ``strong'' bad example. 

2) In Example 2, the condition number $\kappa = 1$, and the spectral radius of the update matrix  of C-CD is approximately $1 - 20/n$ (we check it numerically for $n$ from 10 to 1000).
% \textbf{Observation 1}: 
The gap between C-CD and GD for Example 2 is at least $n/20$, similar to Example 1. Even the constant $20$ is the same. The difference with Example 1 is that for Example 1 there is an $\tO(n)$ gap between GD and R-CD, leading to $\tO(n^2)$ gap between C-CD and R-CD; for Example 2, GD and R-CD converge at the same speed, thus Example 2 does not show the $\tO(n^2)$ gap between C-CD and R-CD. 

In short, this example does \textit{not} provide numerical evidence that the complexity of C-CD is at least $\O(n^4 \kappa_{\mathrm CD} ) $, but only that the complexity of C-CD is at least $\O(n^3 \kappa ) $. 
% Note that the latter claim requires a rigorous proof which is probably nontrivial. The reasons are the same as our proof: the spectral radius may not be easy to estimate, and the spectral radius itself does not constitute a lower bound. 

3) The actual gap between GD/R-CD and C-CD in this example is larger than $n/20$ because the rates of the former do not depend on $\log(1/ \epsilon)$. In fact, GD/R-CD both take 2 epochs to converge while C-CD takes $ \log(1/ \epsilon) n/20$ epochs to converge, thus the "true" gap between R-CD and C-CD is $ \frac{n}{40}  log(1/\epsilon) $. 
For instance, when $n=100$ and $\epsilon= 10^{-15}$,  GD/R-CD both take 2 epochs and C-CD takes about $80$ epochs, and the true gap $40=80/2$ is close to $ \log(1/\epsilon) n/40 = 15 \times 2.5 = 37.5$.
Note that in this computation, $\log(1/\epsilon) =15 $ has a much larger contribution
than $n / 40 = 2.5$, hence by looking at the $37.5$-times gap itself, it is not easy to tell where
this number $37.5$ comes from.
To see the effect of $\tO(n)$, one may need to choose $n > 600$ or even larger. 
%In fact, 
If the gap were completely due to the contribution of $ \log(1/\epsilon) $, this large gap between GD/R-CD and C-CD should be considered a constant gap, since it is impossible to get a theoretical bound of GD/R-CD independent of $\log(1/\epsilon)$ in general. 
In this sense,  the gap of $\log(1/\epsilon)$  is very special to the example, and should not be considered an evidence of GD/R-CD being faster than C-CD. 
% does not match the theoretical prediction.

 % The second observation shows the subtlety of using examples to illustrate the gaps between algorithms. 
 We emphasize that our contribution is not only the proposal of an example (independent of Steven Wright), but also the theoretical analysis related to the example. 
  Just one example empirically showing algorithm A being much faster than algorithm B is not very meaningful for a theoretical  understanding, for at least two reasons.
 First, it is possible that someone comes up with another example showing that algorithm B  is much faster than algorithm A. In fact, there are many numerical examples to show C-CD is faster than R-CD; even though the gap is not as large as $\tO(n)$, one could not claim that such an example does not exist. 
 Second, the gap may be a ``fake'' gap that cannot be explained by any existing theoretical bounds, just like the
 $\log(1/\epsilon)$-factor gap analyzed above.
 Therefore,  in addition to proposing an example,  it is important to prove that the example exhibits the behavior of the theoertical bounds, thus validating the tightness of the established bounds as well as the gap between the bounds.

% The $15$-times factor contributed by $ \log(1/\epsilon) $ is significant in practice. Nevertheeless,
% theoretical bounds. 
%  R-CD/GD perform better than the typical behavior (almost independent of $\log(1/\epsilon)$), 

%
%Note that $ da$ 
%As R-CD performs better than the typical behavior (independent of $\epsilon$), we believe cannot claim the extra gap of log(1/epsilon) in general.
%This observation 

%\subsection{ Different Levels of Lower Bound }\label{sec: possible improvement}
% $1/L_{\min}$ factor,
%Lower bounds are tricky when there are more than one problem parameter.
%We will discuss in more detail what the implications of our results are, and what possible improvement can be made.

\subsection{Non-equal-diagonal Case,  Jacobi Preconditioning and Open Questions}\label{sec: Jacobi Precond}
We will discuss the complexity bounds when the diagonal entries $L_i $'s are not equal. It turns out the ``true'' complexity in this general case is  more subtle
than the equal-diagonal case (i.e. the case where all $L_i$'s are equal) and related to an old problem in numerical linear algebra.
% and there is even a large gap % between our lower bound and the upper bound.

In the previous discussions we often assume $L_i = L_1, \ \forall i$ since one can always
scale the coefficient matrix $ A $ to get a new matrix $ D_A^{ - 1/2} A D_A^{-1/2} $ and modify the algorithm correspondingly.
Such a preprocessing procedure is called Jacobi preconditioning in numerical linear algebra, and is a common data preprocessing trick in machine learning.
 It is very simple to implement and only slightly increases the total complexity of the algorithm.

Nevertheless, one may still wonder what the complexity in the non-equal-diagonal case is.
Our Proposition \eqref{prop: upper bound} implies an upper bound which is more general than
 \eqref{upper bound of complexity, kappa}:
\begin{equation}\label{upper bound of complexity, Lmin}
  \min \left\{ \tO\left( n^3 \kappa \frac{L_{\avg}}{L_{\min}} \right),  \tO\left( n^2  \kappa \log^2 n \frac{L}{L_{\min}} \right) \right\} .
  \end{equation}
 Notice that $L_{\min }$ appears in the denominator of both bounds, thus as $L_{\min} \rightarrow 0 $ both bounds approach infinity.
  Intuitively, this implies that when one coordinate has very little contribution to the whole function C-CD will converge very slowly.
 However, this phenomenon will not happen in practice and the dependency in $1/L_{\min}$ is somewhat artificial.
 In fact, theoretically we can prove a stronger upper bound of C-CD that does not depend on $1/L_{\min}$,
 but instead depends on a new condition number.

 \begin{proposition}\label{prop: upper bound, JacobiPre}(Stronger Upper Bounds)
Consider the same setting as Proposition \ref{prop: upper bound}. We have
% quadratic minimization problem  $\min_x f(x) \triangleq x^TA x + 2b^T x $ where $A \in \R^{n \times n} $
% is positive definite.
 %  For any $x^0 \in \R^n $, let $x^k$ denotes the output of C-CD after $k$ cycles, then
\begin{subequations}
\begin{align}
     f(x^{k+1}) - f^* \leq \min \left\{ 1 - \frac{1}{ n \hat{ \kappa }  }  ,  1 - \frac{ 1 }{ \hat{L} (2 + \log n/ \pi)^2 } \frac{1}{
     \hat{ \kappa }  } \right\}  (f(x^k) - f^*).  \label{Jacobi Pre,upper bound in kappa}  \\
       f(x^{k+1}) - f^* \leq \min \left\{ 1 - \frac{1}{ n^2 \hat{\kappa}_{\mathrm CD}  } ,
   1 - \frac{ 1 }{ \hat{L}^2 (2 + \log n/ \pi)^2 } \frac{1}{ \hat{\kappa}_{\mathrm CD}  } \right\}  (f(x^k) - f^*). \label{Jacobi Pre,upper bound in kappaCD}
\end{align}
\end{subequations}
%We also have
%\begin{equation}
%\end{equation}
Here, the parameters $  \hat{L} = \lambda_{\max}( \hat{A})$, $ \hat{\kappa} = \lambda_{\max}(\hat{A})/\lambda_{\min}(\hat{A}) $,
and $ \hat{\kappa}_{\mathrm CD}  = 1/\lambda_{\min}(\hat{A}) $,
where $\hat{A} = D_{A}^{-1/2} A D_A^{-1/2} $ is the Jacobi-preconditioned matrix, and $D_A$ is a diagonal matrix consisting of all diagonal entries of $A$.
% $f^*$ is the minimum value of the function $f$,   %$x^*$ is the minimum of function $f$, $f^* = f(x^*)$,
% $\hat{\kappa} = \frac{\lambda_{\max}(A)}{\lambda_{\min}(A)}$ is the condition number, $L = \lambda_{\max}(A)$,
% $L_{\max} = \max_i A_{ii}, L_{\min} = \min_i A_{ii}  $ and  $\kappa_{\mathrm CD} =  \frac{ L_{\max} }{ \lambda_{\min}(A) } $.
\end{proposition}

The proof of Proposition \ref{prop: upper bound, JacobiPre} will be given in Appendix \ref{appen: proof of Jacobi Pre upper bound}.
The proof is almost the same as the proof of Proposition \ref{prop: upper bound}
 except that we should replace matrix $A$ and its lower triangular part $\Gamma$ by the Jacobi-preconditioned versions.
 % In fact, the convergence rate of $A$ is related to the matrix $ D_A^{-1/2} \Gamma^T A^{-1} \Gamma D_A^{-1/2} $,
By taking a closer look into the proof, we find that Jacobi-preconditioning is naturally ``embedded'' in C-CD
\footnote{The Jacobi-preconditioning is also embedded in R-CD, but if we pick the coordinates with
probability proportional to $A_{ii}$, the preconditioning effect disappears. }.
This is not surprising since in the update rule \eqref{cyclic CD update} we need to scale the diagonals $A_{ii}$ at each step,
which is similar to Jacobi-preconditioning (but not the same).
Therefore, we can think of $\hat{\kappa}$ as a more appropriate parameter to characterize the complexity of C-CD than
the original condition number $\kappa$.

Proposition \ref{prop: upper bound, JacobiPre} implies the following upper bound
\begin{equation}\label{upper bound of complexity, JacobiPre}
  \min \left\{ \tO\left( n^3 \hat{\kappa}  \right),  \tO\left( n^2  \hat{\kappa} \log^2 n \frac{ \hat{L} }{ \hat{L}_{\text{avg}} } \right) \right\} ,
  \end{equation}
   where $\hat{L}_{\avg} $ is the average of the diagonal entries of $A$, and it equals $1$ since all diagonal entries of matrix $ \hat{A}$  are $1$.
 For the equal-diagonal case, this upper bound reduces to the upper bound \eqref{upper bound of complexity, kappa} since $\hat{A}$ is just a scaled version of $A$.
  % The bound is also almost tight  it reduces to the bound  \eqref{upper bound of complexity, kappa} which is shown to be almost tight.
  Comparing this bound with \eqref{upper bound of complexity, Lmin} which also holds for the non-equal-diagonal case, we find that the $L_{\avg}/L_{\min}$ factor disappears here (since
   this ratio equals $1$ for the Jacobi-preconditioned matrix).
This can be explained as that the factor $L_{\avg}/L_{\min}$ is absorbed into the new condition number $\hat{\kappa}$;
in fact, it is straightforward to prove
  $$   \hat{ \kappa} \leq \kappa \frac{ L_{\max} }{ L_{\min} } , $$
  thus the upper bound $ \tO\left( n^3 \hat{\kappa}  \right) $  immediately implies
  an upper bound $ \tO\left( n^3 \kappa \frac{L_{\max}}{ L_{\min}}  \right) $ that is slightly weaker than
  the first bound in \eqref{upper bound of complexity, Lmin}.
  It is not easy to explicitly compare the second bound of \eqref{upper bound of complexity, JacobiPre}
  and the second bound of \eqref{upper bound of complexity, Lmin}.
 % One way to interpret this disappereance   Due to the reason in the last paragraph,
  % We think the bound \eqref{upper bound of complexity, JacobiPre} is more natural than \eqref{upper bound of complexity, Lmin} for cyclic CD.

 The reason we still present the bound dependent on $\kappa$, instead of   only presenting the bound dependent on $\hat{\kappa}$,
 is because the former bound allows us to compare C-CD with GD.
 With the new bound \eqref{upper bound of complexity, JacobiPre}, a natural question is how to transform it to a bound that only depends on the parameters of the original matrix, such as $\kappa$.
 This is related to the following classical question on Jacobi-preconditioning:
 $$
   \text{What is the relation between the condition number of $A$ and that of the Jacobi-preconditioned matrix? }
 $$
 Intuitively, larger discrepancy in the diagonal entries leads to a larger condition number, thus Jacobi-preconditioning which makes the diagonals equal should reduce the condition number.
 In other words, one may expect that $\hat{\kappa} \leq \kappa $ holds for most of the time, if not always.
 Unfortunately, it is only known that the relation $\hat{\kappa} \leq \kappa $ holds for some special $A$
 (more precisely, when $A$ satisfies Young's property (A)\footnote{If the rows and columns of a $(p + q)$-dim matrix can be rearranged so that the upper $p \times p$ and lower $q \times q$ submatrices are diagonal, then the matrix is said
 to have Property (A) \cite{young1954iterative}. For example, the tridiagonal matrix satisfies Property (A).
 Also note that the question in \cite{forsythe1955best} appears in a different form:
  when is $\hat{A}$ the best conditioned matrix out of all possible diagonally scaled matrix of $A$?
 }), according to Forsythe and Straus \cite{forsythe1955best}.
For our purpose, the exact relation $\hat{\kappa} \leq \kappa $ is not necessary as we are more interested in
the upper bound of $ \hat{\kappa}/\kappa $.
 %;  an extra constant or $\log n$ factor is fine.
There are some simple bounds (see, e.g., \cite[Lemma 3.2, Lemma 3.3]{widlund1971effects}):
$$
  \hat{ \kappa} \leq   \kappa \cdot \min \{ n, \frac{L_{\max}}{L_{\min}} \} .
$$
As a direct corollary, the first bound of \eqref{upper bound of complexity, JacobiPre} implies two upper bounds
 \begin{equation}\label{yet another bound}
  \min \left\{ \tO\left( n^4 \kappa   \right),  \tO\left( n^3 \kappa \frac{L_{\max}}{L_{\min}}   \right) \right\} .
  \end{equation}
We have already seen a variant of the above second bound in \eqref{upper bound of complexity, Lmin}.

 What is more interesting is the first bound in \eqref{yet another bound} $  \tO\left( n^4 \kappa   \right) $, which is $n^2$ times worse than GD, and $n$ times worse than the equal-diagonal case!
 If we want to express the complexity of cyclic CD purely in terms of $\kappa$ for the non-equal-diagonal case,
  $  \tO\left( n^4 \kappa   \right) $ is the best upper bound we have right now.
  There is an $\O(n)$-factor gap between this upper bound and the lower bound $\O(n^3 \kappa)$.
  We believe this gap is artificial and there should be a stronger proof that establishes an upper bound of $\O(n^3 \kappa)$.
  Such a stronger upper bound might be achieved by proving a constant upper bound of $\kappa/ \hat{\kappa } $.
  We pose two open questions:

 \emph{ Open Question 1}: Is there a non-equal-diagonal example that cyclic CD has complexity worse than $  \tO\left( n^3 \kappa   \right) $?
 If yes, what about $  \tO\left( n^4 \kappa \right) $ ?

 \emph{ Open Question 2}:
 Is there a constant upper bound on $\hat{\kappa}/\kappa$, where $\kappa$ and $\hat{\kappa}$ are the condition numbers of $A$ and the Jacobi-preconditioned $\hat{A}$ respectively? If not, what is the best upper bound of $\hat{\kappa}/\kappa$?
 Is there an example that the ratio $\hat{\kappa}/\kappa$ achieves $\O(n)$?

 We stress again that  $  \tO\left( n^3 \hat{\kappa}   \right) $ is a tight bound in general, and
 $  \tO\left( n^3 \kappa   \right) $ is a tight bound when $A$ has equal diagonal entries.
 Thus the above Question 1 is only valid when we consider non-equal-diagonal matrix $A$ and insist on expressing the complexity in terms
 of the condition number of the original matrix.
 In some sense, it is not as essential as the question whether there is an $\O(n^2)$ gap between C-CD and R-CD studied in this paper.
 Nevertheless, it is still a valid question, and becomes more interesting due to its relation to Jacobi-preconditioning.

\subsection{Necessity of Two Types of Bounds}\label{sec: discuss two types of bounds}
Consider the equal-diagonal case (i.e. $L_i = L_1, \ \forall i$) in this subsection.
We will explain the relation between the two types of bounds, one does not involve $\tau =
\lambda_{\max}/\lambda_{\avg} =  L/L_{1}$ and another does.
We argue that it is not easy, if not impossible, to obtain one single tight bound.
We will also suggest slightly stronger bounds that might be the tightest based on the current parameters (again, for the equal-diagonal case).

We denote two bounds related to $\kappa$ as $B_1 =  \cO( n^3 \kappa  )$ and $ B_2 = \cO( n^2 \kappa \tau \log^2 n ) $; the comparison between the bounds related to $\kappa_{\mathrm CD}$ will be similar and thus omitted.
Since we assume $L_1=\dots = L_n $, we have $L \leq \sum_i L_i = n L_{1}$,
and $ 1 \leq \tau = L/L_1 \leq n $. Therefore, in most cases (more precisely, as long as $\tau \in [1, n/\log^2 n$] while the full range of
$\tau$ is $[1,n]$) the bound $B_2 $ is better than $ B_1 $.
However, $B_2$ does not dominate $B_1$ since for our example $B_1$ is tight while $B_2$ is $\O( \log^2 n)$ times worse.
One natural guess is that maybe the best bound is $B_3 = \cO( n^2 \kappa \tau ) $, which is better than both $B_1$ and $B_2$
and also consistent with our example.  Unfortunately, $B_3$ is probably not the right bound since there exists an example such that the $\log^2 n$ factor is unavoidable \cite{oswald1994convergence}.

Now we discuss the result by Oswald \cite{oswald1994convergence}.
The paper  \cite{oswald1994convergence} establishes an upper bound similar to the second bound in \ref{upper bound of complexity, kappa}.
Then the paper constructs an example that ``matches'' the upper bound; more specifically,
in the example both $\kappa $ and $\tau$ are $\O(1)$ while the spectral radius of the iteration matrix of C-CD is $1 - 1/ \O( \log^2 n)  $.
Thus the complexity for this example is at least $\O(n^2 \log^2 n) = \O( n^2 \log^2 n \kappa \tau ) $ which is $\log^2 n$ times larger than $B_3$
\footnote{This statement is not rigorous. It is tricky: the spectral radius of a non-symmetric iteration matrix may not provide the lower bound
	of the convergence rate; extra effort is needed to rigorously build the connection. We will discuss this issue in more details later.
}.  However, this example only ``matches'' the upper bound in a weak sense as the key parameters $\kappa$ and $\tau$ are constants in the example.
In particular, this example has nothing to do with the question whether the extra factor $\tau$ is necessary or not.
It does not exclude the possibility that the worst-case complexity of C-CD were $ \O( n^2 \log^2 n \kappa  ) $
or even $ \O( n^2 \log^2 n \kappa_{\mathrm CD}  ) $ which are very close to the complexity of GD and R-CD respectively.
We think the extra $\tau$ factor is very important for at least two reasons.
First, for most randomly generated matrices the ratio $\tau = \lambda_{\max}/\lambda_{\avg}$ is much larger than $\O(\log n)$.
This can be tested by numerical experiments, and also validated by theoretical results:
for example, for the Wishart random ensamble $A = U^T U$ where the entries of $U$ are standard Gaussian variables,
the ratio $\tau = L/L_1$ is approximately $\O(\sqrt{n} )$ \footnote{
	According to  \cite[Proposition 6.1]{edelman1988eigenvalues}  the maximum eigenvalue is about $ 4 n$
	and the diagonal entries are the lengths of $n$-dimensional random vectors which are $\O(\sqrt{n})$,
	thus after scaling the diagonal entries $\tau = \lambda_{\max}/\lambda_{\avg}$. }.
% Second, if $\tau \leq \O(\log n)$ then $\kappa $ since $\tau$ is usually much larger than $\log n$; in fact, if
Second,  $\tau$ exactly characterizes the theoretical improvement of R-CD over GD.
When $\tau$ is small, the gain of using R-CD is very limited: either the problem is too easy
and GD already performs well, or the problem is so difficult that even R-CD does not help.
Thus the interesting problems for CD-type methods are those with large $\tau$.

Now we know that  the $\log^2n$ factor is necessary for one extreme case $\tau = \O(1)$, and  the $\log^2 n$ factor can be removed for the other extreme case $\tau = n$.
The transition has to happen somewhere in between, and we guess it happens near $\tau \approx \O(\log^2 n)$.
In other words, we guess the ``best'' bound is
$$ B_{\mathrm{conj}} = \O(n^2 \kappa \max\{ \tau, C \log^2 n \}   )  .$$
Although the operator norm of the triangular operator is $\O(\log n)$, we conjecture that when restricted to a certain class of PD matrices ($\tau$ is not too small) the operator norm becomes $\O(1)$.
%\label{thm: lower bound}(Lower bound of CD dependent on $A_{ii}$)
%   For any $x^0 \in R^n $, any $\delta \in (0,1]$, there exists a quadratic function $f(x) = x^TA x + 2b^T x $ such that
%\begin{subequations} \label{desired bound}
%\begin{align}
% \frac{ f(x^k) - f^*}{  f(x^0) - f^*} \geq (1 - \delta) \left( 1 - \frac{ \max_i A_{ii} }{ \beta } \frac{ 2 \pi^2  }{  \kappa }  \right)^{2k + 2 }
% =(1 - \delta) \left( 1 - \left( \frac{ \max_i A_{ii} }{ \beta } \right)^2 \frac{ 2\pi^2 }{ \kappa_{\mathrm CD} }  \right)^{2k + 2 } .
%\end{align}
%\end{subequations}
% % Corollary:  f(x^k) - f^* \geq (1 - \delta) \mu \left( 1 - \frac{ 2\pi^2 }{ n \kappa }  \right)^{2k + 2} \| x^0 - x^* \|^2,   \\ \\
%where  $x^k$ denotes the output of C-CD after $k$ cycles,
% $x^*$ is the minimum of function $f$, $f^* = f(x^*)$, $\kappa = \frac{\lambda_{\max}(A)}{\lambda_{\min}(A)}$ is the condition number of $f$,
% $\kappa_{\mathrm CD} =  \frac{ \max_i A_{ii} }{ \lambda_{\min}(A) } . $
%
\begin{conjecture}
	If $A$ is symmetric PSD with equal diagonal entries and $\tau = \frac{ \lambda_{\max} }{\lambda_{\avg}}  \geq  C_1 \log^2 n $ for some constant $C_1 $, then the lower triangular part (with diagonals) $\Gamma$, defined as $\Gamma_{ij} = A_{ij}, i\leq j$ and $\Gamma_{ij} = 0, i>j$,
	satisfies $\| \Gamma \| \leq  C_2 \| A \| .  $
\end{conjecture}

\subsection{ Precise Comparison of Time Complexity }
% We come back to the equal-diagonal case in this subsection.
In our previous comparison between C-CD and GD/R-CD we have ignored the constants, and we do not state the comparison in a formal result.
Next we will formally compare them and quantify the exact gap in terms of the time complexity. 
In the first result the error is measured in the objective values.
In the second result  the error is measured in iterates, which allows us to add RP-CD into the comparison and get a better bound for R-CD.
% the previous upper and lower bounds
% Theorem \ref{thm: lower bound} provides lower bounds for the convergence rate.
%  As mentioned earlier, the worst-case complexity of C-CD is $O(n)$ times worse than GD and $O(n^2)$ worse than R-CD.
% The difference in the convergence rate can be translated to the difference in the number of iterations to achieve certain accuracy,
% a metric that can be easily computed in practice.
% We report the comparison of the number of iterations in the following result, which is a corollary of Theorem \ref{thm: lower bound}.

The first proposition shows that to achieve any given relative error in objective values,
C-CD takes at least $n/20$ times more operations than GD, and $n^2/40$ times more operations than R-CD. 
The proof of Proposition \ref{prop: real slower, objective} will be given in Appendix \ref{appen: prop objective error compare proof}.
\begin{proposition}\label{prop: real slower, objective} (Compare C-CD with GD, R-CD; objective error)
  % Consider using C-CD, GD and R-CD to solve the quadratic minimization problem
  %\eqref{quadratic min}. %  with $A = A_c$.
  Let $ k_{\mathrm{CCD} }(\epsilon )$, $k_{\mathrm{GD} }(\epsilon ) $ and $ k_{\mathrm{RCD} }(\epsilon ) $
  be the minimum number of epochs \footnote{For a fair comparison, here one epoch of CD or RP-CD means one cycle of all coordinates, and for R-CD one iteration means randomly selecting coordinates for $n $ times. } for C-CD, GD, R-CD
  to achieve (expected) relative error $$ \frac{ E( f(x^k) - f^* ) }{ f(x^0) - f^* } \leq \epsilon $$ for all initial points in $\R^n$ (for C-CD and GD the expectation operator can be ignored).
There exists a quadratic problem such that
   \begin{subequations}
   \begin{align}
     % \lim_{c \rightarrow 1 }
      \frac{ k_{\mathrm{CD} }(\epsilon ) }{ k_{\mathrm{GD} }(\epsilon ) }\geq \frac{n}{2 \pi^2 } \approx \frac{n}{20},  \label{compare GD with CD, coro}  \\
       \frac{ k_{\mathrm{CCD} }(\epsilon ) }{ k_{\mathrm{RCD} }(\epsilon ) }  \geq \frac{n^2 }{ 4 \pi^2 } \approx \frac{n^2 }{40} .  \label{compare CD with RCD, coro}
     \end{align}
     \end{subequations}
\end{proposition}

Remark: It seems that the comparison of C-CD and R-CD is not fair since for R-CD we record the \emph{expected} number of iterations.
Nevertheless, it is easy to prove that to guarantee the same error with probability $1 - \delta$,
 we only need $\log( 1/\delta)$ times more iterations. % , and the extra multiplicative factor $\log( 1/\delta)$ is typically much smaller than $n$.
For simplicity, we just consider the expected number of iterations of R-CD.

In the above Proposition \ref{prop: real slower, objective}, the relative error is defined for the function values;
next, we prove a result in which the relative error is defined for the (expected) iterates.
The proof of Proposition \ref{prop: slower, iterates} is given in Appendix \ref{appen: prop of iterates compare proof}.

\begin{proposition}\label{prop: slower, iterates} (Compare C-CD with GD,R-CD and RP-CD; iterates error)
%  Consider using C-CD, GD and R-CD to solve the quadratic minimization problem  \eqref{quadratic min}. %  with $A = A_c$.
  Let $ K_{\mathrm{CCD} }(\epsilon )$, $ K_{\mathrm{GD} }(\epsilon ) $, $ K_{\mathrm{RCD} }(\epsilon ) $ and $ K_{\mathrm{RPCD} }(\epsilon )$
  be the minimum number of iterations \footnote{Again, for a fair comparison, here one iteration of CD or RP-CD means one cycle of all coordinates,
  and for R-CD one iteration means randomly selecting coordinates for $n $ times. } for C-CD, GD , R-CD and RP-CD
  to achieve (expected) relative error $$ \frac{  \| E(x^k) - x^* \|^2 }{ \|x^0 - x^* \|^2 } \leq \epsilon $$ for all initial points in $\R^n$ (for C-CD and GD the expectation operator can be igonred).
There exists a quadratic problem such that
\begin{subequations}\label{ite, compare CD with others}
\begin{align}
      \frac{ K_{\mathrm{CCD} }(\epsilon ) }{ K_{\mathrm{GD} }(\epsilon ) }     &   \geq \frac{n}{2 \pi^2 } \approx \frac{n}{20},   \label{compare GD with CD, coro}  \\
      \frac{ K_{\mathrm{CCD} }(\epsilon ) }{ K_{\mathrm{RCD} }(\epsilon ) }    &   \geq \frac{n^2 }{ 2 \pi^2 } \approx \frac{n^2}{20} , \label{compare CD with RCD, coro} \\
       \frac{ K_{\mathrm{CCD} }(\epsilon ) }{ K_{\mathrm{RPCD} }(\epsilon ) }   &   \geq \frac{n (n+1) }{ 2 \pi^2 } \approx \frac{n(n+2)}{20}  .
 \end{align}
  \end{subequations}
\end{proposition}

We present the result for two reasons. First, the convergence of iterates is of interest in some scenarios.
Second, %although the computed ratios seem similar to those in Proposition \ref{corollary},
we can obtain stronger bounds. In particular, the ratio we obtained for the squared iterates of R-CD is twice as large as that for the function values of R-CD ($n^2/20$ v.s. $n^2/40$). % and the new bounds seem to be closer to the practical performance.
 Moreover, we are able to add RP-CD % (randomly permuted CD), another popular variant which randomly permutes the coordinates at each cycle)
into comparison for the iterates error.
We do not include RP-CD in Proposition \ref{prop: slower, iterates} since it seems difficult to compute the convergence rate of the objective error for RP-CD.
Despite the advantages, we need to emphasize that the convergence of expected iterate error is a weaker notion of convergence
than the convergence of objective error, because the former does not lead to a high probability convergence rate while the latter does (which is because $f(x^k) - f^* \geq 0$).
If we could bound $E(\| x^k - x^* \|^2)$ instead of $\| E(x^k - x^*) \|^2$, then high probability convergence rate could also be automatically established; but we are unable to bound $E(\| x^k - x^* \|^2)$ for RP-CD either.

%Remark 2: In Proposition \ref{prop: real slower, objective} we define the relative error as the ratio of the function values so as to be consistent with conventional convergence rate results such as \cite[Theorem 2]{nestrov12}. Nevertheless, a similar result holds if we define the relative error
%in terms of iterates or expected iterates, i.e., $ \|E(x^k) \|/\| x^0 \| $ (see Theorem \ref{thm 2: real slower}).

Our theory shows that there exists one example $A = A_c$ such that C-CD takes at least $\frac{n}{  2\pi^2 } \approx \frac{n }{ 20 } $  times more iterations than GD and $\frac{n^2 }{ 2 \pi^2 } \approx \frac{n^2 }{ 20 } $ times more iterations than R-CD to achieve any accuracy $\epsilon$.
While the theory is only established for the case $c$ is very close to $1$ (recall $c$ is the off-diagonal entry),
we will show in simulations that the predicted gaps do really exist for a wide range of $c$.
Note that ``the number of required iterations'' is defined for ``all initial points'' (in other words, ``worst-case'' initial points).
% In practice people usually pick random initial points, which might lead to better performance than the worst initial point;
 We will show in simulations that even for random initialization the gaps observed in practice match those predicted by Proposition
 \ref{prop: slower, iterates}.
 % the differences still exist.
%% Nevertheless, we will briefly analyze why this will happen in th
%% we fix the problem and the choose initial point, one can easily choose initial point which lies in the
%
%%\section{ Convergence Rate of Deterministic Kaczmarz Method}
%%
%%Kaczmarz method \cite{kaczmarz1937angenaherte} is one of the most popular method to solve overdetermined system
%%$$
%%  U x = b ,
%%$$
%%where $U $ is a full row rank $m \times n$ matrix with $m \geq n$.
%
%

\section{Overview of the Proofs }\label{sec: overview of proofs}
% \subsection{Proof Sketch and Analysis}

\subsection{Overview of Proof of Proposition \ref{prop: upper bound}}
We present two types of bounds: in the equal-diagonal case, the first type only depends on $\kappa$ or $\kappa_{\mathrm CD}$, and the second type depends on $\kappa$ and $\tau$. In the non-equal diagonal case, both bounds depend on $L_{\min}$. 
The first type of bounds can be established by the same techniques as in \cite{sun2015improved}, though  \cite{sun2015improved} only considers non-strongly convex case.
% which is $\O(\log n)$-times stronger for the extreme case $L/L_1 = n$.
% ombines the two types of bounds and improve the constant for one type of bounds.
We give a unified proof framework that leads to both types of bounds. %  and our proof does not require Hoffman bound.
Our proof can be divided into two stages. The first stage is to relate the convergence rate with the spectral norm
of a matrix $ \Gamma^{-1}A\Gamma^{-T}$, which can be proved by two different approaches (from different perspectives):
one is from optimization which views C-CD as inexact GD; the other is from linear algebra which studies
the spectral radius of the iteration matrix $I - \Gamma^{-1}A$.
% Interestingly, our proof clearly shows the connection of the two approaches.
Note that $\Gamma^{-1}A$ is non-symmetric, thus the latter method requires an extra symmetrization technique which
relaxes the spectral radius by the spectral norm. As we will see later, such a technique cannot be used in the proof of the lower bound,
and other techniques are needed for that proof.
In the second stage, we estimate $ \| \Gamma^{-1}A\Gamma^{-T} \| $ via two different methods, leading to the two types of bounds.
As discussed in Section \ref{sec: discuss two types of bounds}, each bound is tight in one scenario, thus the two bounds cannot be combined into one single bound.

\subsection{Overview of Proof of of Theorem  \ref{thm: lower bound}}
\subsubsection{Difficulties}
In general, to prove a lower complexity bound, one only needs to construct an example and compute the convergence rate of the example.
% The difficulty usually lies in finding the appropriate example and the rest is just calculation.
However, in our case, computing the convergence rate of the example is not easy due to (at least) two reasons.

First, Gauss-Seidel method can be written as a matrix recursion and its convergence rate is related to the spectral radius of the update matrix.
It turns out that the spectral radius of our example does not have a closed form expression; in fact,
the spectral radius depends on the roots of an $n$-th order equation.
To resolve this issue, we notice that as the constructed matrix tends to singular (i.e. the off-diagonal entries tend to $1$) the $n$-th order equation will become simple; based on this fact, we are able to bound the spectral radius asymptotically (as off-diagonal entries tend to $1$,
but still for fixed $n$).

Second, the update matrix of Gauss-Seidel method is a non-symmetric matrix.
A simple, though usually ignored, fact is that for non-symmetric matrix recursion, the spectral radius of the iteration matrix
is not the lower bound of the convergence factor in the \emph{real} domain.
Note that if we were allowed to pick initial points in the complex domain, then the spectral radius did provide a lower bound of the convergence rate; but here we are only interested in the real initial points.
We have not seen a general method to deal with this issue; fortunately,
the example we constructed happens to exhibit some special structure so that we can provide a lower bound of the convergence rate.
We will discuss this difficulty in more details in Section \ref{sec: why non-symmetric is difficult}.
There is actually one more difficulty caused by the non-symmetry of the iteration matrix:
it is even harder to bound the function error. Fortunately again, we are able to resolve this difficulty due to
another special property of the problem. See more details in Step 3 of the outline in Section \ref{sec: proof outline}.

% There is a subtle difference between the lower bound and the upper bound, resulting in the issue of non-symmetry and thus a more involved proof for the lower bound, even if it is just for a specific example.
The issue of non-symmetry does not appear in the proof of the upper bound in Proposition
\ref{prop: upper bound}
because a symmetrization technique is used. %  (though most existing proofs do not emphasize this).
% The first obstacle is rather unexpected: in the proof of
Assuming $x^*=0$, we need to compute the convergence rate of $ f(x^k) = (x^k)^T A (x^k) = \| y^k \|^2 $,
where $y^k = U x^k $ in which $U$ satisfies $A = U^T U$. It is easy to get the matrix recursion
$ y^{k+1} = (I - U\Gamma^{-1}U ) y^{k} $, thus one needs to bound the spectral radius of $M_f = I - U\Gamma^{-1}U $.
The spectral radius of a non-symmetric matrix is not easy to directly bound, thus in that proof we
instead upper bound the spectral norm $\| M_f \| = \sqrt{ \| M_f^T M_f \| }$, which gives a upper bound of $\rho(M_f)$.
%Although not explicitly stated, the proof of Prop. \ref{prop: RCD rate} lies in an observation that $M_f^T M_f$
%has a simple form after some transformation.
However, the relaxation from $\rho(M_f)$ to $\| M_f\|$ is not reversible; in other words, even if
we prove that for our example $\| M_f\|$ is large, this does not mean $\rho(M_f)$ is large (or C-CD is slow).
Thus we have to consider the original non-symmetric form $ U\Gamma^{-1}U$ or $\Gamma^{-1}A$ for the lower bound.
% the objective error $f(x^k) - f^*$ decays at a rate proportional to the spectral radius
%of $M_f = I - \Gamma^{-1}A L^{-T}$, but this spectral radius is only useful in the proof of the upper bound.
% The reason is that we do \emph{not} have a matrix recursion of $f(x^k) - f^*$; instead,
% we can only prove $ \frac{f(x^{+1}) - f^*}{f(x^k) - f^*} = \frac{ (d^k)^T d^k  }{ d^k M_f d^k  } $, where
% $d^k$ is iteration dependent.

\subsubsection{Why Non-symmetric Iteration Matrix Causes Difficulty}\label{sec: why non-symmetric is difficult}
We discuss why the spectral radius of a non-symmetric iteration matrix does not necessarily lead to a lower bound of the convergence rate
(for real initial points).
Consider the following matrix recursion
\begin{equation}\label{recursion general}
 y^{k+1} = M y^{k} .
 \end{equation}
 %--------------RESERVE, somewhere else--------
 % All iterative algorithms we will consider later can be written as a matrix recursion (or randomized matrix recursion), denoted as
% For simplicity, for randomized algorithm we just consider the expected iteration matrix.
% For the matrix recursion, we may consider either the convergence rate of iterates $\{x^k \}$ or the convergence rate of
% optimal gap $f(x^k) - f*$; in most cases these two sequences converge in similar rate, and we will not distinguish these two in
% the following discussion. we have used earlier
 %----------------------
% Following the convention,
We say a sequence $\{ y^k \}$ converges with convergence rate $\tau$
if $\| y^k \| \leq C \tau^k $, where $C$ is a constant.

A basic result is that if $M$ is symmetric the convergence rate of $\|y^k\|$ is exactly $\rho(M)$.
% However, if $M$ is not symmetric, $\rho(M)$ is not a lower bound of the convergence rate (to the best
% of our knowledge). % To see why, let us recall how to prove the symmetric case.
%To show that the convergence rate is at most $ \rho(M)$, we need to represent an arbitrary vector as a linear combination of the eigenvectors
%of $M$; however, if $M$ is non-symmetric its eigenvectors may not span the whole space.
How to prove this result?
For the lower bound (i.e. the convergence rate is at least $  \rho(M) $), we need to pick the initial point to be the eigenvector of $M$ corresponding to $\rho(M)$.
This proof no longer works for non-symmetric $M$ since its eigenvectors may be complex vectors.
One way to resolve this issue is to pick the real part of the complex eigenvector; %  as we have done in our proof;
however, this approach requires additional assumptions to work.  % does not work in general (at least to us).
More specifically, suppose $ M v = \lambda v$, where $\lambda = \rho(M) = |\lambda| e^{i \theta}$, and pick the initial point $y^0 = \text{Re}(v) = \frac{1}{2}(v + \bar{v})$.
The update \eqref{recursion general} leads to
$$
y^k = M^k y^0 = \frac{1}{2}M^k (v + \bar{v})
  = \frac{1}{2} (\lambda^k v + \bar{\lambda}^k \bar{v}) = \text{Re}( \lambda^k v )
  = |\lambda|^k \text{Re}( e^{i k \theta} v ) .
$$
Suppose $v = (r_1 e^{i \phi_1}, \dots, r_n e^{i \phi_n}  )$,
then
\begin{equation}\label{xk expression for general initial point}
 \| y^k \| = \rho(M)^k \sqrt{  r_1^2 \cos^2(k\theta + \phi_1) + \dots + r_n^2 \cos^2(k\theta + \phi_n)   } .
\end{equation}
For the lower bound, we want to prove
\begin{equation}\label{lower bound}
 \| y^k \| \geq C \rho(M)^k, \ \forall k ,
 \end{equation}
  where $C$ is a constant.
 % and, for a meaningful result, should not be too smaller than $1$.
Without any additional assumption, this is impossible: if $\phi_j = 0, r_j = 1, \forall j $ and $ k \theta = \frac{\pi}{2} + 2 m \pi $ for some integer $m $, then $\| y^k \| = 0$.
Intuitively, when all $\phi_j$'s are close to each other, it is hard to lower bound $\| y^k\|$;
but if all $\phi_j$'s are evenly spread out, then $\| y^k\|$ can be lower bounded.
For our problem, it turns out the phase $\phi_j$ goes to $2j \pi/n$ as $c$ goes to $1$, which is the the nicest case we can expect
(the phases are equally spaced).
%There is an extra nice property: all $ $
In such a nice case, we are able to give a simple lower bound of $\| y^k \|$.

% Some simple assumptions can work; for example,
One might wonder whether it is easy to obtain a lower bound in the general case under mild assumptions.
We consider the simplest case $n = 2$.
 If $r_1 = r_2 > 0$ and $0< |\phi_1 - \phi_2| < \pi/2 $, then $\| y^k \| \geq \rho(M)^k r_1 |\sin( (\phi_1 - \phi_2)/2)| $.
However, if $r_1 \neq r_2$, then even for $n = 2$ we need more assumptions to find a lower bound. Such assumptions
can be a relation between $r_1/r_2$ and $\phi_1 - \phi_2$, which look non-intuitive and seem to be constructed merely for theory.
Moreover, it is hard to express the corresponding bound (e.g. $ r_1 |\sin( (\phi_1 - \phi_2)/2)|$) as
a function of simple parameters of the original problem. From a practical point of view,
the ``constant'' $ r_1 |\sin( (\phi_1 - \phi_2)/2)|$ can be so small that it already meets the practical need.
These issues will become even more complicated when $ n > 2$.
As a conclusion, when the iteration matrix is non-symmetric, it seems difficult to lower bound the convergence rate in general.

% even though in practice we do not need to worry too much about

%---RESERVE: UPPER BOUND of Non-Symmetric----------------------------------
% As a side remark, bounding the
%------------------------------------
% While for general problems it is difficult to prove \eqref{lower bound},

\subsubsection{Proof Outline for Theorem  \ref{thm: lower bound}}\label{sec: proof outline}
The detailed proof is divided into three steps. We will construct an example $\min_x x^T A_c x$ where the coefficient matrix
 has diagonal entries $1$ and off-diagonal entries $c \in(0, 1)$. Obviously $x^* = 0$ is the unique minimum and $f^* = 0$.

 In Step 1, we compute the spectral radius of the iteration matrix asymptotically.
More specifically, we show that the eigenvalues of the iteration matrix are given by $\lambda_j = 1 - q_j^n$, where
$q_j$'s are the roots of the equation $ q^n(1 - c + q) = 1 $.
While the closed form expression of $\lambda_j$ is difficult to compute (in fact, for a special case,
a very complicated closed form of an infinite series is given in \cite{glasser1994quadratic}), we observe that as $c \rightarrow 1 $, $q_j$'s tend
to the $n$-th unit roots.
We then prove that as $c \rightarrow 1$ the spectral radius of the iteration matrix tends to roughly $ 1 - \frac{2\pi^2 }{ n \kappa } $.

In Step 2, we prove that for a certain real initial point $x^0$,
 the relative error $\frac{ \|x^k - x^* \|^2 }{ \| x^0 - x^* \|^2 }$ is lower bounded by $O\left( \left( 1 - \frac{2\pi^2 }{ n \kappa } \right)^{2k} \right)$. In other words, the sequence $ \{ \|x^k - x^* \|^2 \} $ converges at a rate lower bounded by the spectral radius $ 1 - \frac{2\pi^2 }{ n \kappa } $.
  % rather slowly; in fact, its convergence rate
 The initial point we choose is the real part of the eigenvector corresponding to the spectral radius of the iteration matrix $I - \Gamma^{-1}A$.
 A crucial property is that the eigenvector has an expression $(1,q, \dots, q^n)$ where $q$ is an complex eigenvalue of $\Gamma^{-1}A$,
 thus the phases of the initial elements are roughly $2j/\pi, j=1,\dots, n$.
 This property makes the calculation of the relative error $\frac{ \|x^k - x^* \|^2 }{ \| x^0 - x^* \|^2 }$ possible.

In Step 3, we prove that the relative error $\frac{ f(x^k) - f^* }{ f(x^0) - f^* } = \frac{ f(x^k)}{ f(x^0)}$ is also lower bounded
by $\O\left( \left( 1 - \frac{2\pi^2 }{ n \kappa } \right)^{2k} \right)$.
Again, the special structure of the example is crucial for this step.
Unlike GD method where the iteration matrix $I - \frac{1}{\beta}A$ has the same eigenvectors as $A$, the iteration matrix
of CD method $I - \Gamma^{-1}A$ has different eigenvectors from $A$.
As we pick $x^0$ to be the real part of an eigenvector of $\Gamma^{-1}A$,
it is not clear a priori how to bound $f(x^k) = (x^k)^T A x^k $ and $f(x^0) = (x^0)^T A x^0 $.
Of course one can lower bound $f(x^k)$ by $\lambda_{\min}(A) \| x^k \|^2 $ and upper bound $f(x^0) $ by $\lambda_{\max}(A) \| x^0\|^2$
to get a lower bound of $\frac{f(x^k)}{f(x^0)}$, but this will introduce an extra factor $\frac{\lambda_{\min}(A) }{\lambda_{\max}(A) }
= \frac{1 -c}{ 1 -c + cn} \approx \frac{ 1 -c }{ n }   $ which tends to $0$ as $c \rightarrow 1$.
Thus we need to give a tighter bound of either $f(x^k)$ or $f(x^0)$.
We choose to bound $ f(x^0) $ differently: it turns out $f(x^0) =  (x^0)^T A x^0  $ can be upper bounded by $\lambda_{\min}(A) \| x^0\|^2$ plus some negligible term (as $c \rightarrow 1$), which makes $\frac{ f(x^k)  }{ f(x^0) } $ very close to $ \frac{ \| x^k\|^2 }{ \|x^0 \|^2 }  . $
The crucial property here is that for our example, the eigenvector corresponding to the spectral radius of the iteration matrix $I - \Gamma^{-1}A$
is very close to the eigenvector corresponding to the minimum eigenvalue of $A$.
Needless to say, this property does not hold for general matrix $A$.

%  since $x^k$ is just eigenvalues of $ \Gamma^{-1} A$, not

\section{Formal Proof of Theorem \ref{thm: lower bound} }\label{sec: proof of Main Result} 
 This section contains a full proof of Theorem \ref{thm: lower bound} except the proof for some technical lemmas.
 % a lemma that computes the eigenvalues of a matrix.
 
Assume the initial point is up to our choice for now. We will show in the end of the proof how to deal with an arbitrary initial point.

For any constant $c \in (0,1)$, consider minimizing the following quadratic function
\begin{equation}\label{quadratic min}
  \min_{x \in \R^n } f(x) \triangleq x^T A_c x ,
\end{equation}
where $ A_c \in \R^{n \times n }$ is defined as
% , we pick $A = A_c$ as follows:
\begin{equation}\label{Ac def}
 A_c = \begin{bmatrix}
    1          &  c      &  \dots  &  c   \\
    c          &  1      &  \dots  &  c   \\
    \vdots      &  \vdots  &   \ddots &  \vdots \\
    c          &   c     &   \dots &   1
 \end{bmatrix}
\end{equation}
Simple calculation shows that $A_c$ is a positive definite matrix, with one eigenvalue $1-c$ with multiplicity $n-1$ and one
eigenvalue $1-c + cn $ with multiplicity $1$. Thus the condition number of the matrix is
\begin{equation}\label{kappa expresssion}
  \kappa = \frac{1 - c + cn}{ 1 - c}.
\end{equation}
 The optimum of the problem is $x = (0; 0; \dots ; 0) $.
 Solving this problem is also equivalent to solving a linear system of equations $Ax = 0$.

\textbf{Step 1}: Computing the spectral radius of the iteration matrix, asymptotically.
The following lemma shows that the eigenvalues of the matrix $A_c $ are the roots of a polynomial equation.
The proof of Lemma \ref{lemma of equation} is given in Appendix \ref{appen: lemma equation proof}.
\begin{lemma}\label{lemma of equation}
Suppose $A = A_c$ is defined by \eqref{Ac def} and $\Gamma$ is the lower triangular part of $A_c$ (with diagonals), and denote $\hat{c} = 1 - c$.
% Then $Z = \Gamma^{-1}A $ is a non-singular matrix with $n$ eigenvalues given as follows: $1$ is an eigenvalue of $Z$;
  Suppose the $n + 1 $ roots of
  \begin{equation}\label{eig value equation}
   q^{n} (q - 1 + c) = c q
   \end{equation}
  are $q_0 , q_1, \dots, q_{n-1}, q_n $ among which $q_0 = 0, q_n = 1$, then % $ \lambda_k = \frac{ \hat{c} - \hat{c}q_k }{\hat{c}- q_k } = 1 - q_k^n, k=0, 1,\dots, n-1 $,
   \begin{equation}\label{lambda k def}
  \lambda_k = \frac{(1-c)(1-q_k ) }{ 1- c - q_k } = 1 - q_k^n , \ k=0, 1, \dots, n -1
  \end{equation}
are all $n $ eigenvalues of $Z = \Gamma^{-1}A$.
\end{lemma}

  Note that $\lambda_{n} =  1 - q_n^n  = 0$ is not an eigenvalue of $Z$.
Eliminating a factor of $q$ in \eqref{eig value equation}, we have that $q_1, \dots, q_{n-1}, q_n = 1$
are the $n$ roots of the equation $ q^{n-1} (q - 1 + c) = c $. % has $n $ complex roots, denoted as $q_1(c), q_2(c), \dots, q_n(c) = 1 $.
 Intuitively, as $c$ goes to $1$, the equation becomes $q^n = 1$, thus the roots $q_k $ will converge to  an $n$-th root of unity. 
 The formal statement is given below and the proof is given in Appendix \ref{appen: Lemma of cts roots proof}.
 
 \begin{lemma}\label{claim of n roots convergence}
	% Consider the equation  $ q^{n-1} (q - 1 + c) = c $, for some $c < 1$.
 	There exists some $c_0 \in (0,1)$ such that when $c \in  ( c_0, 1 )$ the following holds:
    the equation $  q^{n-1} (q - 1 + c) = c $ has exactly one solution $q_k(c)$ such that
 	$ |  q_k(c) - e^{  i 2\pi k/n } |  \leq \frac{1}{2} \sin \frac{ \pi }{ n}   $ for $k =  1, \dots, n$; moreover, 
 	$$
 	\lim_{c \rightarrow 1} q_k(c) = e^{ i 2\pi k/n }, \; \forall k. 
 	$$
 \end{lemma}

%There is a small technical difference between the above result and what we claimed: in the above result each root  $q_k$ lies in the
%set $ \{ z \mid  z^{-1} \in B( \eta_k, \epsilon_0  )  \} $, while in the statement of Lemma \ref{lemma of roots converge} we claim each root  $q_k$ lies in a ball with the center $ \eta_k^{-1} $. 
% It is easy to bridge this gap since $ \{ z \mid  z^{-1} \in B( \eta_k, \epsilon_0  )  \} $ implies $ | z - \eta_k^{-1} | < \sin(\pi / n)  $. 
% 
% More formally, let $ \epsilon_1 =  \min\{ \sin( \pi/n),  0.49  \}  $. By a similar argument as above, for any $c \in ( 3/(3 + v( \epsilon_1 ) , 1 ) $,  the equation $  q^{n-1} (q - 1 + c) - c $ has exactly one root $q_k(c)$ such that
% $ |1/ q_k(c) - e^{ - i 2\pi k/n } |  < \epsilon_1 $ for $k = 0, 1, \dots, n-1$, and
% $  \lim_{c \rightarrow 1} q_k(c) = e^{ i 2\pi k/n }, \; \forall k.  $
%One can verify that $ |1/ q_k(c) - e^{ - i 2\pi k/n } |  < \epsilon_1 = \min\{ \sin( \pi/n),  1/3 \}$  implies 
%$ | q_k(c) - e^{  i 2\pi k/n } | <  1.5 \sin( \pi /n )  $.

% (another way to see this is to write the optimization problem in the real domain as
% $ v_k(\epsilon) = \min_{ \theta \in [0, 2 \pi ] }  | (\eta_k +  \epsilon e^{i \theta } )^n - 1 |   $ ). 

Suppose $c \in (c_0, 1)$ from now on.
Note that   $ \min_{ 1\leq j < k \leq n} | e^{ i 2\pi j/n } - e^{ i 2\pi k/n }| = |1 - e^{ i 2 \pi /n}| = 2 \sin (\pi/n)$, thus by Lemma \ref{claim of n roots convergence} $q_k(c), k=1, \dots, n$ are distinct roots of the equation $  q^{n-1} (q - 1 + c) = c $. 
For simplicity of notations, we denote $q_k(c)$ as $q_k$, which satisfies
 \begin{equation}\label{q asymptotic}
  \lim_{c \rightarrow 1} q_k =  e^{i 2 k \pi/ n }.
  \end{equation}
Obviously $\lim_{c \rightarrow 1 } \lambda_k = 0, \forall k $.

Next we prove % that the limit of the third term of the product
\begin{equation}\label{Jk expression}
 J_k \triangleq \lim_{c \rightarrow 1 }   \frac{ 1/ \kappa }{ 1 -  |1 - \lambda_k |  } = \frac{1}{ 2 n \sin^2(k\pi/n) }  , \ k = 1, \dots, n - 1.
 \end{equation}

For notational convenience, let $\hat{c} \triangleq 1 - c, \hat{\lambda}_k  = 1 - \lambda_k $.
Then we have
\begin{equation}\label{intermed lambdak}
 \hat{\lambda}_k \overset{\eqref{lambda k def}}{=} 1 - \frac{  \hat{c}(1 - q_k ) }{ \hat{c} - q_k } = \frac{ (\hat{c} - 1) q_k }{  \hat{c} - q_k }
 = \frac{c q_k }{ q_k - \hat{c} } = \frac{ c }{ 1 - \hat{c} /q_k}  .
 \end{equation}
Then
\begin{equation}
\begin{split}
 J_k  \overset{\eqref{kappa expresssion} }{=} \lim_{ \hat{c} \rightarrow 0 }  \frac{ \hat{c} }{ \hat{c} + c n} \frac{1}{ 1 - |\hat{\lambda}_k| }
   =  \frac{ 1 }{ n } \lim_{ \hat{c} \rightarrow 0 }   \frac{\hat{c}} {1 - |\hat{\lambda}_k | }
   \overset{ \eqref{intermed lambdak} }{=}  \frac{ 1 }{ n } \lim_{ \hat{c} \rightarrow 0 }   \frac{\hat{c}}  { 1 - \left| \frac{c }{ 1 - \hat{c}/q_k } \right| }  \\
    = \frac{ 1 }{ n } \lim_{ \hat{c} \rightarrow 0 }   \frac{\hat{c} |1 - \hat{c}/q_k| }  { |1 - \hat{c}/q_k| - c } .
  \end{split}
 \end{equation}
% Assume $ q_k^{-1} = r e^{i \theta} $ (for simplicity we drop the subscript $k$ for $r$ and $\theta$),
% then $\lim_{\hat{c} \rightarrow 0 } r = 1,  \lim_{\hat{c} \rightarrow 0 } \theta = - 2 k \pi/ n  $.
 Since $ \lim_{ \hat{c} \rightarrow 0}  |1 - \hat{c}/q_k|  = 1 $,  % = \lim_{ \hat{c} \rightarrow 0}  | q_k - \hat{c}|/|q_k|
from the above relation we have
 \begin{equation}
\begin{split}
   J_k & = \frac{ 1 }{ n } \lim_{ \hat{c} \rightarrow 0 } \frac{\hat{c} }  { |1 - \hat{c}/q_k| - c }
     =  \frac{1}{n}  \lim_{ \hat{c} \rightarrow 0 }  \frac{ \hat{c} ( |1 - \hat{c}/q_k| + c)  }{  | 1 - \hat{c}/q_k|^2 - c^2 }
     =  \frac{2 }{n}  \lim_{ \hat{c} \rightarrow 0 }  \frac{ \hat{c}  }{  | 1 - \hat{c}/q_k|^2 - c^2 }  \\
     & = \frac{2 }{n}  \lim_{ \hat{c} \rightarrow 0 }  \frac{ \hat{c} }{ 1 + \hat{c}^2 / |q_k|^2 - 2 \hat{c} \text{Re}(1/q_k) - c^2 } \\
     & = \frac{2 }{n}  \lim_{ \hat{c} \rightarrow 0 }   \frac{ \hat{c} }{ \hat{c}( 1+c ) + \hat{c}^2 / |q_k|^2 - 2 \hat{c} \text{Re}(1/q_k)  }  \\
     & = \frac{2 }{n}  \lim_{ \hat{c} \rightarrow 0 }   \frac{ 1 }{  1+c + \hat{c} / |q_k|^2 - 2  \text{Re}(1/q_k)  }
   \end{split}
 \end{equation}
Since $\lim_{\hat{c} \rightarrow 0} |q_k| = 1, \lim_{\hat{c} \rightarrow 0} \text{Re}(1/q_k) = \cos(- 2k\pi/n ) $, the above relation can be further simplified to
  \begin{equation}\nonumber
\begin{split}
   J_k   % = \frac{ 2 }{n }  \lim_{ \hat{c} \rightarrow 0 }  \frac{1}{ 2 - 2 \text{Re}(1/q_k) }  \\
     = \frac{ 2 }{n } \frac{1}{ 1 + 1 + 0 - 2 \cos(- 2k\pi/n ) }
      =   \frac{1}{ 2 n \sin^2(k\pi/n) }  ,
 \end{split}
 \end{equation}
 which proves \eqref{Jk expression}.

\textbf{Step 2: } Bound the relative iterates error.  % the optimal gap for the iterates.

To simplify the notations, let $q = q_1$ and $\lambda = \lambda_1 = \frac{ \hat{c} - \hat{c}q_1 }{\hat{c}-q_1}$ from now on.

According to the proof of Lemma \ref{lemma of equation}, $ \tilde{v} = (\tilde{v}_1; \dots; \tilde{v}_n)  $ is an eigenvector of $Z = \Gamma^{-1} A $  corresponding to $\lambda $, where
 \begin{equation}
   \tilde{v}_j = \frac{c}{\lambda - \hat{c}} q^{j -1}, j=1,\dots, n.
\end{equation}
We scale each entry of $\tilde{v}_j$ by a constant $\frac{\lambda - \hat{c}} {c} q $ to get a new vector  $ v = (v_1,\dots, v_n)$,
where
 \begin{equation}\label{expression of v, again}
    v_j = q^{j },  \ j = 1, \dots, n.
\end{equation}
Obviously $v$ is also an eigenvector of $Z$ corresponding to $\lambda$, i.e. $Z v = \lambda v$.

Now pick the initial point
 $ x^0 = \text{Re}( v ) $.
 Suppose
 \begin{equation}\label{r,theta def}
   q = r e^{ i \theta } ,
 \end{equation}
  where $i = \sqrt{-1} $, $r > 0$ and $\theta \in [0, 2 \pi)$,  then
  \begin{equation}\label{x0 expression}
  x^0_j =   \text{Re} ( r^j  e^{\sqrt{-1} j \theta } ) = r^j \cos( j \theta), \ j =1,\dots, n .
      % ( r^{1-n} \cos(n-1)\theta ; r^{2-n} \cos(n-2)\theta ; \dots ; r^{-1} \cos \theta ;  1
\end{equation}
Since $x^0 = \frac{1}{2} (v + \bar{v} )$, and $v$ and $\bar{v}$ are eigenvectors of $ M = I - \Gamma^{-1}A$ with
eigenvalues $1 - \lambda$ and $1 - \bar{\lambda}$ respectively, we have
$$
  x^k = M^k x^0 = \frac{1}{2}M^k( v + \bar{v} ) =  \frac{1}{2} ( (1 - \lambda)^k v + (1 - \bar{\lambda})^k \bar{v} )
      = \text{Re}( ( 1 - \lambda )^k v  )  \overset{ \eqref{lambda to n and 1 -q  relation} }{=} \text{Re}( q^{kn} v  ).
$$
According to \eqref{expression of v, again}, the $j$-th entry of $x^k$ is
$$
  x^k_j = \text{Re}(q^{kn + j} )  = r^{kn + j} \cos( kn + j)\theta , \ j=1,\dots, n.
    % ( r^{1 + (1-k)n} \cos(n-1)\theta ; r^{2-n} \cos(n-2)\theta ; \dots ; r^{-1} \cos \theta ;  1  )
$$

Note that $ r = |q| \leq 1 $ (otherwise C-CD will diverge, but we know from classical results that C-CD always converges for solving our problem), then we have
\begin{equation}\label{xk bound ineq}
\begin{split}
  \|x^k \|^2 = \sum_{j =1}^n r^{ 2kn + 2 j } \cos^2 [ (kn + j )\theta ]  &
  \geq    r^{2kn + 2n } \sum_j \cos^2  [ (kn + j )\theta ]   \\
  & =  r^{ (2k+2)n } \frac{1}{2} \left( \sum_j \cos( 2 k n \theta + 2 j \theta ) + n \right) .
  \end{split}
\end{equation}

To calculate the sum in the above expression, we will need the following standard equality;
for completeness, the proof of this claim is given in Appendix \ref{appen: claim sum sin proof}.
\begin{claim}\label{claim: sum cos equality}
For any $z, \phi \in \R $, we have
\begin{equation}\label{sum cos, simpler}
  \sum_{ j =1}^n \cos( z +  j \phi ) = \frac{ \sin(n\phi/2) \cos( z + (n+1)\phi/2 ) }{ \sin (\phi/2 )}.
\end{equation}
\end{claim}

Applying \eqref{sum cos, simpler} to the expression in \eqref{xk bound ineq}, we have
\begin{equation}\label{xk true lower bound}
   \|x^k \|^2 \geq   \frac{1}{2} r^{ (2k+2)n } \left( \frac{ \sin(n \theta ) \cos(  2kn \theta + (n+1)\theta ) }{ \sin \theta } + n  \right)
   \geq  \frac{1}{2} r^{ (2k+2)n  } \left(  n - \left| \frac{ \sin(n \theta ) }{ \sin \theta  }  \right| \right) .
\end{equation}
Similar to \eqref{xk bound ineq} (but bound $r^{2j}$ from above by $1$), we have
\begin{equation}\label{x0 bound ineq}
\begin{split}
  \|x^0 \|^2  & = \sum_{j=1}^n r^{ 2j } \cos^2( j \theta )
  \leq    \sum_j \cos^2 ( j  \theta)    \\
   & =  \frac{1}{2} \left( \sum_j \cos(  2j \theta ) + n   \right)
    =  \frac{1}{2} \left( \frac{ \sin(n \theta ) \cos(  (n+1)\theta ) }{ \sin  \theta } + n  \right)
    \leq   \frac{1}{2} \left(  n + \left| \frac{ \sin(n \theta ) }{ \sin \theta  }  \right| \right) .
  \end{split}
\end{equation}
Combining the above two relations, we have
\begin{equation}\label{ratio of xk x0 bound}
   \frac{ \| x^k \|^2  }{ \| x^0\|^2 } \geq r^{2kn + 2n } \frac{  n -  | \sin(n \theta ) / \sin \theta  |  }{  n +  | \sin(n \theta ) / \sin \theta  |  }
    = r^{2kn + 2n } \omega_c,
\end{equation}
where $$ \omega_c \triangleq \frac{ n -  | \sin(n \theta ) / \sin \theta  | }{ n + | \sin(n \theta ) / \sin \theta  | }  . $$
According to \eqref{q asymptotic},  $q_1 = re^{i \theta} $ converges to $ e^{i 2\pi/n} $ as $c \rightarrow 1$,
thus  $ \theta \rightarrow  2 \pi/n $ and $| \sin(n \theta)/ \sin \theta |  \rightarrow 0  $ as $c \rightarrow 1$,
which further implies $\omega_c \rightarrow 1$ as $c \rightarrow 1$.

\textbf{Step 3: } Bound the relative objective error. % optimal gap for the objective values.

Suppose $A = U^T U$, and denote $y^k = U x^k$, then $$ \| y^k \|^2 = (x^k)^T U^T U x^k = (x^k)^T A x^k = f(x^k)  . $$
%It is not easy to derive a lower bound of $\| y^k \|^2$ from the explicit expression of $y^k$.
%We instead bound $\|y^k \|^2$ based on the bound of $\| x^k\|^2$. More specifically,
Note that the minimum eigenvalue of $A$ is $\hat{c} = 1 - c$, thus
\begin{equation}\label{yk bound}
 \| y^k \|^2 = \| U x^k \|^2 = (x^k)^T A x^k \geq \hat{c} \| x^k \|^2
 \overset{ \eqref{xk true lower bound}  }{\geq}  \frac{  \hat{c} }{2} r^{ (2k+2)n  } \left(  n - \left| \frac{ \sin(n \theta ) }{ \sin (\theta )}  \right| \right) .
% \overset{ \eqref{xk true lower bound} }{\geq} \frac{c^2 }{2 \hat{c} } r^{2kn } \left(  n - \left| \frac{ \sin(n \theta ) }{ \sin (\theta )}  \right| \right).
\end{equation}
%Since $\omega_c \rightarrow 1$ as $c \rightarrow 1$, for any $\delta > 0$, there exists $ c_{\mathrm{u}, 1} < 1 $ such that
%\begin{equation}\label{omega c bound}
%\omega_c > 1 - \delta, \ \forall \ c \in ( c_{\mathrm{u}, 1}, 1 ).
%\end{equation}

%For any $\delta > 0$, pick $c \in (\max\{ c_{\mathrm{u}, 1}, c_{\mathrm{u}, 2}  \}, 1 )$,
%% substituting the above relation and \eqref{coeff close to 1} into \eqref{yk ratio closed form},
% we obtain
%$$
% \frac{ f(x^k) - f^* }{ \mu \| x^0 - x^* \|^2 } = \frac{ \| y^k\|^2 }{ \hat{c} \| x^0 \|^2 }
% \overset{ \eqref{yk bound} }{\geq } \frac{ \| x^k\|^2 }{ \| x^0 \|^2 } \overset{\eqref{ratio of xk x0 bound} }{\geq}
%  r^{2kn + 2n } \omega_c \overset{ \eqref{omega c bound},\eqref{r^n bound for some c}}{\geq } (1 - \delta) \mu \left( 1 - \frac{ 2\pi^2 }{ n \kappa }  \right)^{2k+2}.
%$$
%This prove the first half of \eqref{desired lower bound}.
% To prove the second half of \eqref{desired lower bound},

We need to give an upper bound of $\| y^0 \|^2$.
% e need to use the explicit expression of $y^0$.
Denote
\begin{equation}\label{Gamma def}
 \gamma_j = \text{Re}(q^{  j} ) = r^j \cos(j \theta), \quad  S = \sum_{ l =1}^n  \gamma_l .
\end{equation}
 Then the expression of $x^0$ given in \eqref{x0 expression} becomes
 \begin{equation}\label{x0 expression, again}
   x^0 = ( \gamma_1 ; \dots ; \gamma_n  ).
 \end{equation}
Since the $j$-th row of $A$ is $(c,\dots, c, 1, c, \dots, c)$ where $1$ is in the $j$-th position,
we can compute the $j$-th entry of $A x^0$ as
\begin{equation}\label{y expression}
  (Ax^0)_j = c \sum_{ l =1}^n  \gamma_l + (1-c)  \gamma_j
   \overset{\eqref{Gamma def}}{=}  c  S + \hat{c}  \gamma_j .
\end{equation}
%\begin{equation}\label{y expression}
%  y_j^k = c \sum_{ l =1}^n  \frac{c }{\hat{c}} \text{Re}(q^{kn-n+ l} ) + (1-c)  \frac{c }{\hat{c}} \text{Re}(q^{kn-n+ j} )
%   =  \frac{c^2}{\hat{c} }  \sum_{ l =1}^n \text{Re}(q^{kn-n+ l} )  + c \text{Re}(q^{kn-n+ j} ).
%  %      = \frac{c^2}{\hat{c} } \text{Re}(  q^{kn - n + 1} \frac{  q^n - 1 }{ q - 1}  ) + c \text{Re}(q^{kn-n+ j} ) .
%\end{equation}
Then we have
\begin{equation}\label{y0 first bound}
\begin{split}
  \|y^0 \|^2 & = (x^0)^T A x^0
  =\sum_j \gamma_j( c S + \hat{c} \gamma_j )
       =  c S \sum_j \gamma_j + \hat{c} \sum_j \gamma_j^2
  %     \overset{\eqref{x0 expression, again}}{=}
       = c S^2 + \hat{c} \|x^0 \|^2 .
 \end{split}
\end{equation}
%  Previously Wrong expression; the below expression is for ||Ax||^2.
%       \sum_{j=1}^n ( \frac{c^2}{\hat{c} }  \Gamma + c \gamma_j )^2
%             = \sum_j ( \frac{c^4}{\hat{c}^2 }\Gamma^2 + c^2 \gamma_j^2 +  \frac{ 2 c^3 }{ \hat{c}} \gamma_j \Gamma )  \\
%     &        =  \frac{ n c^4}{\hat{c}^2 }  \Gamma^2  + c^2 \sum_j \gamma_j^2  +  \frac{ 2 c^3 }{ \hat{c}} \Gamma^2  \\
%     &        =  \frac{c^3 (nc + 2\hat{c}) }{ \hat{c}^2 } \Gamma^2 + c^2 \sum_j \gamma_j^2   \\
%    &         =  \frac{c^3 (nc + 2\hat{c}) }{ \hat{c}^2 } \Gamma^2  +  \hat{c}^2 \|x^0 \|^2,

We will show that the second term $\hat{c} \|x^0 \|^2$ is the dominant term, which will imply
that $ \frac{ \| y^k \|^2 }{ \| y^0 \|^2 } \approx \frac{\hat{c} \| x^k\|^2 }{ \hat{c} \|x^0 \|^2 }$. To this end,
we need to bound $S^2$. By the definition of $ S $ in \eqref{Gamma def}, we have
$$
  S = \sum_j \gamma_j = \text{Re}(\sum_j q^{j}) = \text{Re}( q \frac{1 - q^n}{ 1- q}  ),
$$
thus
$$
  S^2 \leq \left|  q \frac{1 - q^n}{ 1- q} \right|^2
   = r^{2} \left| \frac{1 - q^n}{ 1- q} \right|^2
   =  r^{2} \left| \frac{\lambda}{ 1 - q }  \right|^2
   = r^{2} \left| \frac{ \hat{c} }{ \hat{c} - q }  \right|^2
   = \hat{c}^2 \left| \frac{ q }{ \hat{c} - q }  \right|^2  .
$$
Substituting the above relation and \eqref{x0 bound ineq} into \eqref{y0 first bound}, we get
\begin{equation}\label{y0 bound final}
\begin{split}
  \| y^0 \|^2
   & \leq c  \hat{c}^2  \left| \frac{ q }{ \hat{c} - q }  \right|^2 +
  \hat{c} \frac{1}{2} \left(  n + \left| \frac{ \sin(n \theta ) }{ \sin (\theta )}  \right| \right)   \\
  &  =  \hat{c}  \alpha_c,
\end{split}
\end{equation}
%  Previously Wrong expression; the below expression is for ||Ax||^2.
%  & \leq \frac{ c^3 (nc + 2\hat{c}) }{ \hat{c}^2 } \hat{c}^2 r^{-2n} \left| \frac{ q }{ \hat{c} - q }  \right|^2 +
%  \hat{c}^2 \frac{c^2 }{2 \hat{c}^2 } r^{ - 2n }  \left(  n + \left| \frac{ \sin(n \theta ) }{ \sin (\theta )}  \right| \right)   \\
%  &  = r^{-2n} \left(  c^3 (nc + 2\hat{c}) \left| \frac{ q }{ \hat{c} - q }  \right|^2
%    + \frac{ c^2 }{2}   \left(  n + \left| \frac{ \sin(n \theta ) }{ \sin (\theta )}  \right| \right)     \right)  \\
%  &  = r^{-2n} \alpha_c,
where in the last equality we introduce the definition
$$
  \alpha_c \triangleq  c \hat{c} \left| \frac{ q }{ \hat{c} - q }  \right|^2
    +   \frac{1}{2} \left(  n + \left| \frac{ \sin(n \theta ) }{ \sin (\theta )}  \right| \right).
$$
As $ c \rightarrow 1 $, we have $\hat{c} = 1 - c \rightarrow 0$, $ \left| \frac{ q }{ \hat{c} - q }  \right| \rightarrow 1  $ and
$ \left| \frac{ \sin(n \theta ) }{ \sin (\theta )}  \right| \rightarrow 0 $, thus
\begin{equation}\label{alpha limit}
 \alpha_c \rightarrow 0 + \frac{ n }{2}  = \frac{n}{2} ,\quad  \text{as } c \rightarrow 1.
 \end{equation}

Combining \eqref{yk bound} and \eqref{y0 bound final}, we get
\begin{equation}\label{yk ratio closed form}
  \frac{\|y^k\|^2}{\| y^0 \|^2 } \geq \frac{ \hat{c} }{2  } r^{(2k+2) n } \left(  n - \left| \frac{ \sin(n \theta ) }{ \sin (\theta )}  \right| \right)
  \frac{1 }{  \hat{c} \alpha_c }
  = \beta_c r^{(2k + 2)n} ,
\end{equation}
where in the last equality we introduce the definition
$$
 \beta_c \triangleq \frac{ 1 }{2  \alpha_c } \left(  n - \left| \frac{ \sin(n \theta ) }{ \sin (\theta )}  \right| \right) .
$$
According to \eqref{alpha limit} and the fact that $\left| \frac{ \sin(n \theta ) }{ \sin (\theta )}  \right| \rightarrow 0  $ as $c \rightarrow 1$, we have
\begin{equation}\label{coeff beta goes to 1}
  \beta_c \rightarrow  \frac{1}{ n } (n - 0 ) = 1,  \quad \text{as } c \rightarrow 1 ,
\end{equation}
which implies that for any $\delta > 0$, there exists $c_{\mathrm{u}, 1} < 1 $ such that
\begin{equation}\label{coeff close to 1}
  \beta_c > 1 - \delta, \ \forall c \in ( c_{\mathrm{u}, 1} , 1 ).
\end{equation}

By the relation between $\lambda_1$ and $q_1$ and the definition of $r$, we have
$$
    |1 - \lambda_1 | \overset{\eqref{lambda k def}}{=}   |q_1|^n  \overset{\eqref{r,theta def}}{=} r^n .
$$
According to \eqref{Jk expression},
we have
$  \lim_{c \rightarrow 1 }   \frac{ 1/ \kappa }{ 1 -  r^n  } =  \lim_{c \rightarrow 1 }   \frac{ 1/ \kappa }{ 1 -  |1 - \lambda_1 |  } =
\frac{1}{ 2n \sin^2(\pi/n) }  > \frac{ n }{ 2 \pi^2 } .
$
Therefore, there exists $ c_{\mathrm{u}, 2} < 1 $ such that
$
  \frac{ 1/ \kappa }{ 1 -  r^n  } > \frac{ n }{ 2 \pi^2 } ,  \ \forall \ c \in (  c_{\mathrm{u}, 2} , 1 ),
$
i.e.
\begin{equation}\label{r^n bound for some c}
  r^n > 1 - \frac{ 2 \pi^2 }{ n \kappa },  \ \forall \ c \in ( c_{\mathrm{u}, 2} , 1 ).
\end{equation}

For any $\delta > 0$, pick $c \in (\max\{ c_{\mathrm{u}, 1}, c_{\mathrm{u}, 2}  \}, 1 )$ and
substituting \eqref{coeff close to 1} and \eqref{r^n bound for some c} into \eqref{yk ratio closed form},
 we obtain
\begin{equation}\label{boudn f over f}
	 \frac{f(x^k) - f^* }{ f(x^0) - f^* } =  \frac{\|y^k\|^2}{\| y^0 \|^2 } \geq (1 - \delta ) \left( 1 - \frac{ 2 \pi^2 }{ n \kappa }  \right)^{2k + 2} .
\end{equation}
This proves   \eqref{desired bound: kappa}.
To prove the bound \eqref{desired bound: kappaCD}, notice that for our example 
\begin{equation}\label{bound of kappa over kappaCD}
	 \frac{\kappa}{\kappa_{\mathrm CD} } = \frac{L}{L_{\avg} } =  \frac{ 1 - c + cn }{ 1 }
	\rightarrow n ,    \text{ as } c \rightarrow 1.
\end{equation}
According to \eqref{r^n bound for some c} and \eqref{bound of kappa over kappaCD}, 
 for $c$ close enough to $1$, we have $  r^n > 1 - \frac{ 2 \pi^2 }{ n^2 \kappa_{\mathrm CD} } $. 
 Substituting this relation and \eqref{coeff close to 1} into \eqref{yk ratio closed form}, we obtain the desired bound \eqref{desired bound: kappaCD} (similar to the calculation done in \eqref{boudn f over f}).

At last, for an arbitrary initial point $x^0$ our results still hold since C-CD is invariant with respect to the simultaneous shift
of the initial point and the space of variables.
More specifically, pick $c \in (0,1)$ such that \eqref{desired bound} holds and let $v$ be the eigenvector of $A_c$ given in \eqref{expression of v, again}.
Consider using C-CD  to solve the problem $$ \min_z (z - x^0 + v)^T A_c (z - x^0 + v) $$
starting from $x^0$. Applying a linear transformation $z = x - x^0 + v$, this algorithm becomes C-CD for
 solving $ \min_x x^TA_c x $ starting from $ v $ (the optimal solution $x^* $ and optimal value $f^*$ will change accordingly).
 Applying the result we have proved for this case, we get the desired result for the case with initial point $x^0$.
% We only need to transform this initial point $x^0$ to the initial point $v$
 \QED

\section{Numerical Experiments}\label{sec: simulation}

In this section, we present numerical experiments of C-CD, R-CD, RP-CD (randomly permuted CD, i.e., use random orders in each cycle) and GD for minimizing quadratic functions.
%Such experiments have been performed (at least partially) many times before, and different comparison results have been reported: 
In the literature, some papers present examples that C-CD performs better than R-CD (e.g.,
\cite{Beck09}), and others present opposite examples (e.g. \cite{shalev2013stochastic}).
% It seems hard to draw a conclusion about the practical performance of C-CD and R-CD/RP-CD.
%The lack of understanding of the comparison leaves an impression that
%the theoretical analysis does not affect how practitioners use the algorithms:
% no matter how good or bad the results are, practitioners still do not know which algorithm is better and thus have to try both.
% could just try both C-CD and R-CD, and pick the better one.
%The good performance of C-CD in many scenarios seems to indicate that the nice complexity results of R-CD is just theoretical
% and does not mean the superiority of R-CD.
% We argue, at least for the problem we are studying, that theoretical analysis does provide guidance in practice.
Nevertheless, instead of simply stating ``sometimes C-CD converges faster, sometimes R-CD converges faster'',
we will demonstrate that the size of off-diagonal entries (relative to diagonal entries) affect the performance of C-CD. We summarize our numerical findings below:
% explore what features of the problem affect the performance of C-CD.
%In other words, we will try to answer the following question:
%\begin{equation}
% \text{In which cases does C-CD perform well/poorly, in practice?  }
%\end{equation}
%% Of course it is very difficult to give a complete answer to this question, but we will show that
%% ``the level of diagonal dominance'' seems to be one major indicator of the performance of C-CD.
%% show that there are some simple features of the problem that largely determine the performance of C-CD.
%A complete answer to this question seems very difficult, but we have some interesting numerical findings as summarized below:
\begin{enumerate}
	\vspace{-0.15cm}
	\item  C-CD is very slow for solving our example \eqref{Ac def}, as predicted by our theory,
	even for random initial points and non-asymptotic $c$ (e.g. $c > 0.5$).
	In addition, the gap between C-CD and GD/R-CD/RP-CD in our simulation matches the theoretical prediction very well. 
	%  in Proposition  \ref{prop: slower, iterates} very well.
	% For some typical accuracy requirements (e.g. $1e-2, 1e-5, 1e-10$) and
	%   For reasonable $c$ (e.g. $c > 0.8$),
	%   starting from random initial points, C-CD does take about $N/20$ more iterations than $GD$ and $N^2/40$ more iterations than R-CD/RP-CD
	%   to achieve certain accuracy.
	\vspace{-0.15cm}
	\item In the equal-diagonal case, the ratio $\tau = \lambda_{\max}/\lambda_{\avg} = L/L_{\avg}$ is an important indicator of the performance of C-CD. For randomly generated $A$, when $A$ has large $\tau$, C-CD converges much slower than R-CD/RP-CD; when $A$ has small $\tau$, C-CD usually converges as fast as (sometimes faster than) R-CD/RP-CD.
	In these random examples, $\tau$ is closely related to ``off-diagonals-over-diagonals-ratio''
	 (the ratio of the average magnitude of the off-diagonal entries
	over that of the diagonal entries), thus the size of the off-diagonal entries can be a simple indicator of the performance of C-CD.
	%  The ``diagonal dominance ratio'' (ratio of the absolute values of diagonal entries and off-diagonal entries) affects the performance of C-CD.
	\vspace{-0.15cm}
	\item Similar to many experiments in earlier works, we also find that C-CD converges much faster than GD in all cases we test. This is opposite to the theory based on worst-case analysis.
	The bizarre discrepancy between theory and practice has motivated our work, but our work cannot explain
	but rather validate this discrepancy, and new types of analysis might be needed.
	
	%  According to the theory, the performance should be related to the ratio $L/L_{\min}$; since $L/L_{\min}$ is closely related
	%  to the ``diagonal dominance ratio''
	
	%----------RESERVE----------
	%  \item Rather surprisingly, for most random models of $A$ we tested, the spectral radius of the iteration matrix of C-CD
	%  is smaller than that of R-CD/RP-CD. This is true even when C-CD converges much slower than R-CD/RP-CD.
	%   We will discuss this issue in more details later.
	%--------------------
	
\end{enumerate}

%-------------------------------------------
% \subsection{Spectrum and Practical Behavior }
%Before presenting the numerical experiments, we first discuss the relation between eigenvalues of the iteration matrix and
%the practical performance. The goal of this subsection is to show how we can explain the practical behavior by computing
%the eigenvalues, and what level of non-rigorous we have to suffer in such explanations.
% This is very important for understanding the practical behavior of iterative algorithms for quadratic minimization.
% Next we discuss a more important issue: the true performance of a matrix recursion depends on the distribution of the eigenvalues.
%-------------------------------------------

%We have explained in Section \ref{sec: why non-symmetric is difficult} why the spectral radius may not be the indicator of the worst-case convergence rate when the iteration matrix is non-symmetric; more precisely, it provides an upper bound, but not necessarily the lower bound.
%% It is difficult to find a general theory for the worst-case convergence rate of non-symmetric matrix recursion.
%Nevertheless, in most practical scenarios we can still roughly regard the spectral radius as an indicator of the worst-case performance.
% we just need to be aware that  the truth might be different (though unlikely).

\subsection{Experiments for the Bad Example}
We first present simulation results for our example \eqref{Ac def}.

Our theoretical results are established for the asymptotic case $c \rightarrow 1$, and we want to test
whether the same holds for fixed $c$. 
Although the value $ \rho(M)$ does not necessarily represent the convergence rate when $M$ is non-symmetric for C-CD (we have only proved $c \rightarrow 1$ case, not for general $c<1$),
we will still use $\rho(M )$ as a plausible indicator. 
We have computed $1 - \rho(M)$ where $M$ is the (expected) iteration matrix of C-CD, R-CD, RP-CD and GD for various values of $c \in (0,1)$.
% e will regard $\rho(M)$ as a rough lower bound of the worst-case convergence rate of C-CD.
% For C-CD, what we have proved is as $c \rightarrow 1, $ $ \rho(M)$ is indeed the lower bound of the convergence rate; we have not proved the same holds for $c<1$ but we conjecture so. 
%------------------------
%
%\footnote{It is possible to rigorously prove that $\rho(M)$ serves as a lower bound for the converge rate for specific $c < 1$. We can compute the phases of each element of the eigenvector of $M$ corresponding to the spectral radius; as long as the phases are not too clustered (i.e. close to each other), we can give a lower bound of \eqref{xk expression for general initial point}. However, even though the simulations clearly show that these phases are not clustered, to rigorously prove such an obvious phenomenon seems to require long and careful computation. }.
%------------------------
In the last three columns, we divide the values $1 - \rho(M) $ of R-CD, RP-CD and GD by the value of C-CD,
and the resulting ratio represents how many times faster they are than C-CD.
% This is because it roughly takes $ \log(1/\epsilon)/(1 - \rho(M))   $ iterations to achieve relative error $\epsilon$ in the worst-case.
In the rows indicated by ``1(theory)'', we use the theoretical values $  n^2/2\pi^2 \approx n^2/20, n(n+1)/2\pi^2\approx n(n+1)/20,  n/2\pi^2 \approx n/20 $
according to Proposition \ref{prop: slower, iterates}.

\begin{table}[!htbp]
	\caption{ { Comparison of C-CD, R-CD, RP-CD and GD for our example $A = A_c$ } }
	\label{table spectral radius}
	\centering
	\begin{tabular}{|c|c|c|c|c|c|c|c|}
		\hline
		\multirow{2}{*}{ c } &
		\multicolumn{ 4}{c|}{ $1 - \rho(M)$, where $M$ is iteration matrix  } &  \multicolumn{ 3}{c|}{ Ratio over C-CD }  \\  %   \multicolumn{2}{c|}{\multirow{2}{*}{Multi-Row and Col}} &
		\cline{2-8}  % Draw line in 2-5 columns
		&  C-CD  &  GD  & R-CD & RP-CD  &  GD  &   R-CD & RP-CD   \\ %  & \multicolumn{2}{c|}{}
		\hline
		\hline
		\multicolumn{8}{|c|}{ n = 20 }  \\
		\hline
		0.5 & 7.6e-1   &  4.8e-1 & 4.0e-1 & 5.2e-1 &   0.63  &   0.53 &  0.68   \\
		\hline
		0.8 & 1.4e-2   &  1.2e-1 &  1.8e-2 & 2.0e-1  &  0.85  & 12.6 & 14.3  \\
		%       \hline
		%  0.95 &   2.7e-3 &  2.6e-3 & 4.9e-2  & 5.2e-2 &  0.98  & 19.2 & 19.7  \\
		\hline
		0.99 & 4.98e-4 &  5.05e-4 &  1e-2 & 1.03e-2  &  1.01  &  20.0 & 20.7   \\
		\hline
		1 (theory) & -- & --  &  -- & --  &  1.01  & 20.2 & 21.2  \\
		\hline \hline
		\multicolumn{8}{|c|}{ n = 100 }  \\
		\hline
		0.5 & 3.8e-3   & 9.9e-3 &  0.39 & 0.50 &   2.6  &  103  &  132 \\
		\hline
		0.8 & 6.1e-4   & 2.5e-3 & 0.18  & 0.20 & 4.08 & 297 & 328 \\
		\hline
		0.99 & 2.0e-5 & 1.01e-4 & 0.01  & 0.01 &  5.02 & 494 &  497  \\
		\hline
		1 (theory) & -- & --  &  -- & --  & 5.07 & 506 &  512 \\
		\hline  \hline
		\multicolumn{8}{|c|}{ n = 1000 }  \\
		\hline
		0.5 & 3.9e-5   & 9.99e-4 &  0.39 & 0.50 &   25.4  &  9999  &  12717 \\
		\hline
		0.8 & 6.2e-6   & 2.5e-4 & 0.18  & 0.20 & 40.5 & 29411 & 32480 \\
		\hline
		0.99 & 2.01e-7 & 1.01e-5 & 0.01  & 0.01 &  50.2 &  49407 & 49704   \\
		\hline
		1 (theory) & -- & --  &  -- & --  & 50.7 & 50600 &  50760 \\
		\hline
	\end{tabular}
\end{table}

Table \ref{table spectral radius} clearly shows that for $c = 0.8$ the gap between C-CD and other methods is already large, and rather
close to the theory value for $c=1$.
In fact,  the gap   between GD and C-CD for $c = 0.8$  is around $80\%$ of the theoretical gap for $c=1$.
When $c = 0.99$, the gap is about $99\%$ of the predicted gap.
These findings indicate that the gap between GD and C-CD can be uniformly expressed as $ c $ times the theoretical gap for $c = 1$;
similarly the gap between R-CD/RP-CD and C-CD can be expressed as $c^2$ times the theoretical gap for $c = 1$.
This phenomenon suggests that the lower bound \eqref{lower bound, 2nd type}
is not only true for $\tau = \lambda_{\max}/\lambda_{\avg} = n$, but also for many other values of $\tau$ (at least for $\tau  \geq 0.5 n$).
Nevertheless, a rigorous validation requires a non-asymptotic analysis for a given $c$, not for $c\rightarrow 1$, which seems not easy.
% represents the convergence rate when $M$ is symmetric (for R-CD, RP-CD, GD), but 
% can be explained by the following conjecture:
%in the inequalities of Theorem \ref{thm: lower bound}, replacing $n$ by $\frac{ \beta  }{ \max_i A_{ii} } $, which equals $1 - c + cn  \approx c n$
%for our example \eqref{Ac def}, leads to the extra $c$-factor for GD and the extra $c^2$-factor for R-CD/RP-CD.

%------------ RESERVE: ALL RANGE ANALYSIS ------------------------------
%\begin{conjecture}\label{thm: lower bound}(Lower bound of CD dependent on $A_{ii}$)
%   For any $x^0 \in \R^n $, any $\delta \in (0,1]$, there exists a quadratic function $f(x) = x^TA x + 2b^T x $ such that
%\begin{subequations} \label{desired bound}
%\begin{align}
% \frac{ f(x^k) - f^*}{  f(x^0) - f^*} \geq (1 - \delta) \left( 1 - \frac{ \max_i A_{ii} }{ \beta } \frac{ 2 \pi^2  }{  \kappa }  \right)^{2k + 2 }
% =(1 - \delta) \left( 1 - \left( \frac{ \max_i A_{ii} }{ \beta } \right)^2 \frac{ 2\pi^2 }{ \kappa_{\mathrm CD} }  \right)^{2k + 2 } .
%\end{align}
%\end{subequations}
% % Corollary:  f(x^k) - f^* \geq (1 - \delta) \mu \left( 1 - \frac{ 2\pi^2 }{ n \kappa }  \right)^{2k + 2} \| x^0 - x^* \|^2,   \\ \\
%where  $x^k$ denotes the output of C-CD after $k$ cycles,
% $x^*$ is the minimum of function $f$, $f^* = f(x^*)$, $\kappa = \frac{\lambda_{\max}(A)}{\lambda_{\min}(A)}$ is the condition number of $f$,
% $\kappa_{\mathrm CD} =  \frac{ \max_i A_{ii} }{ \lambda_{\min}(A) } . $
%\end{conjecture}
%------------------------------------------

\begin{figure}[H]
	\centering
	\subfigure[$10^2$ iterations]{
		\label{Fig.sub.1}
		\includegraphics[width=0.45\textwidth,height=5cm]{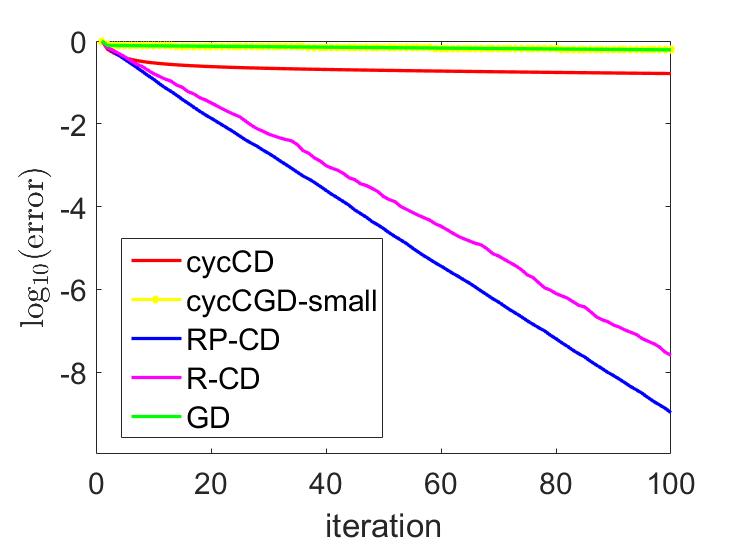}}
	\subfigure[$10^4$ iterations]{
		\label{Fig.sub.2}
		\includegraphics[width=0.45\textwidth,height=5cm]{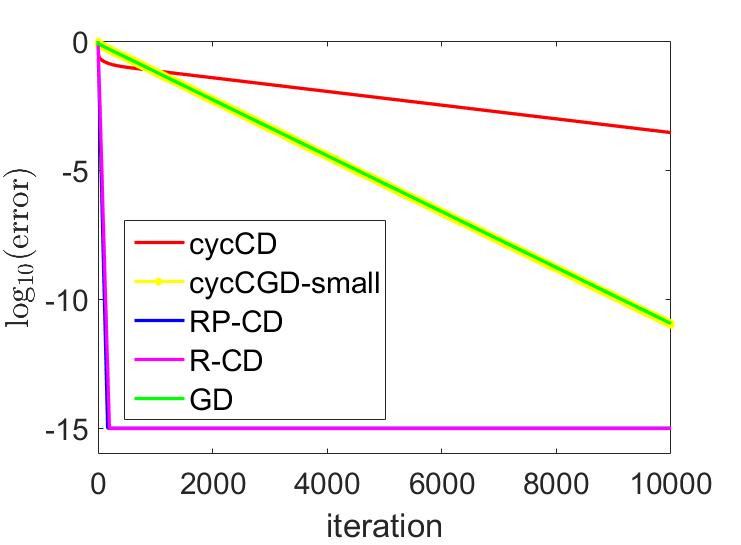}}\hfill %(hfill¿ÉÒÔÈÃtext¸ü½ô´Õ£©
	\caption{Relative error $ \frac{f(x^k) - f^*}{f(x^0)-f^*} $ v.s. iterations, for 5 methods C-CD, cyclic CGD with small stepsize $1/\lambda_{\max}(A)$, randomly permuted CD, randomized CD and GD.
		Minimize $f(x) = x^T A x$,  $n = 100, A = A_c$ with $c = 0.8$.  }
	\label{Fig1ExampleN100}
\end{figure}

Table \ref{table spectral radius} only shows the convergence rate of various methods for the \emph{worst} initial points.
Now we present some simulation results for \emph{random} initialization.
Figure \ref{Fig1ExampleN100} compares the performance of five methods C-CD, cycCGD-small (cyclic CGD with small stepsize $1/\lambda_{\max}$), RP-CD, R-CD and GD, for minimizing $f(x) = x^T A x$, where $n = 100, A = A_c$ with $c = 0.8$.
The left figure shows the first 100 iterations, and the right figure shows $10^4$ iterations.
In the right figure, the large gap predicted by theory clearly exists:
C-CD is about 4 times slower than GD, and GD is about $80$ times slower than R-CD/RP-CD (which means
C-CD is about 320 times slower than R-CD/RP-CD, matching Table \ref{table spectral radius}).
Figure \ref{Fig1ExampleN100} shows that RP-CD is slightly faster than R-CD, which also matches
Table \ref{table spectral radius}.
% One phenomenon not predicted by Table \ref{table spectral radius} is that C-CD is faster than GD in the first 10 iterations;
% is not surprising as Table \ref{table spectral radius} is only related to

\subsection{Experiments for Random Data}
Next, we discuss numerical experiments for randomly generated $A$; for simplicity, we will normalize the diagonal entries of $A$ to be $1$.
Since different random distributions of $A$ will lead to different results, we test many distributions and try to understand for which C-CD performs well/poorly.
To guarantee that $A$ is positive semidefinite, we generate a random matrix $U$ and let $ A = U^T U$.
%(in which case columns of $U$ can be viewed as data points and $A$ can be viewed as empirical covariance matrix if $U$ has zero mean).
We generate the entries of $U$ i.i.d. from a certain random distribution,
% (which implies that the columns of $U$ are independent),
such as $\mathcal{N}(0,1)$ (standard Gaussian distribution), $\text{Unif}[0,1]$ (uniform $[0,1]$ distribution), log-normal distribution, etc.
It turns out for most distributions C-CD is slower than R-CD,
but for standard Gaussian distribution C-CD is better than R-CD.

% At first, we thought that it is because Gaussian distribution is ``nice'' that C-CD performs well.
Inspired by the numerical experiments for the example \eqref{Ac def}, we suspect that the performance of C-CD depends on how large the off-diagonal entries of $A$ are (with fixed diagonal entries). 
%When $U_{ij} \sim \mathcal{N}(0,1)$, $U^T U $ has small off-diagonal entries since $E(A_{ij}) = 0, \ i\neq j$. When $U_{ij} \sim \text{Unif}[0,1]$, $U^T U$ has large off-diagonal entries since  $E( (U^TU)_{ij}) = n/4 $. The difference of these two scenarios lie in the fact that $ \text{Unif}[0,1] $ has non-zero mean while $\mathcal{N}(0,1)$ has zero mean. It is then interesting to compare the case of drawing entries from $ \text{Unif}[0,1] $ with the case of drawing entries from $ \text{Unif}[-0.5,0.5]$.
% turns out for random distributions, a similar phenomenon exists.
% which show that larger off-diagonal entries lead to worse performance of C-CD
To quantify the ``off-diagonals-over-diagonals-ratio'', we define
$$
\chi_{i} =  \frac{ \sum_{j \neq i} |A_{ij}|  }{  A_{ii} } = \sum_{j \neq i} |A_{ij}|, \ i =1,\dots, n,  \quad \chi_{\mathrm{avg}} = \frac{1}{n} \sum_{i} \chi_i,  \quad  \tau = \frac{ \lambda_{\max} }{ \lambda_{\avg} } = \frac{ L }{ L_{\avg} } = L,
$$
where we have used the assumption $A_{ ii } = 1, \forall \; i$ and its consequence $L_{\avg} = 1$.
%Recall that
%$
%  \tau = \frac{ L }{ L_{\avg} } = \frac{ \lambda_{\max} }{ \lambda_{\avg} }.
%$
%we can simply write
%$$
%   \chi_{i} = \sum_{j \neq i} |A_{ij}| , \quad \tau = \lambda_{\max} = L.
%$$
Obviously $ L = \lambda_{\max} \leq 1 + \max_i \chi_i $. In many examples we find $\lambda_{\max}$ to be close to $1 + \chi_{\mathrm{avg}}$, especially when both of them are large.

\begin{figure}
%\begin{tabular}{cc}
\begin{minipage}{0.48\linewidth}
\centerline{\includegraphics[width=\textwidth]{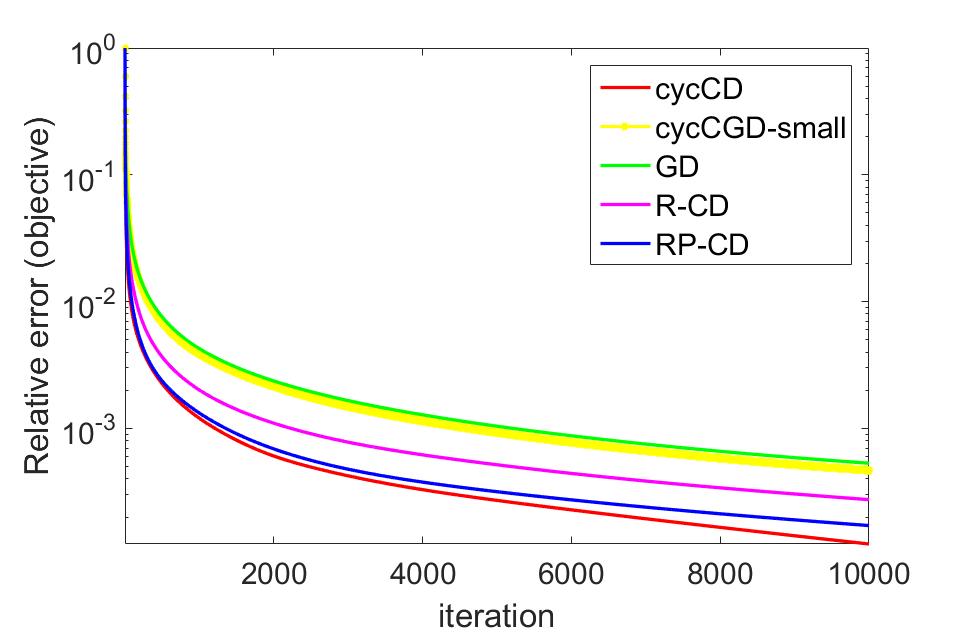}}
\centerline{(a1) Gaussian zero mean. $L \approx  3.8$, $\chi_{\max} \approx 7.9$. }
\end{minipage}
\hfill
\begin{minipage}{.48\linewidth}
\centerline{\includegraphics[width=\textwidth]{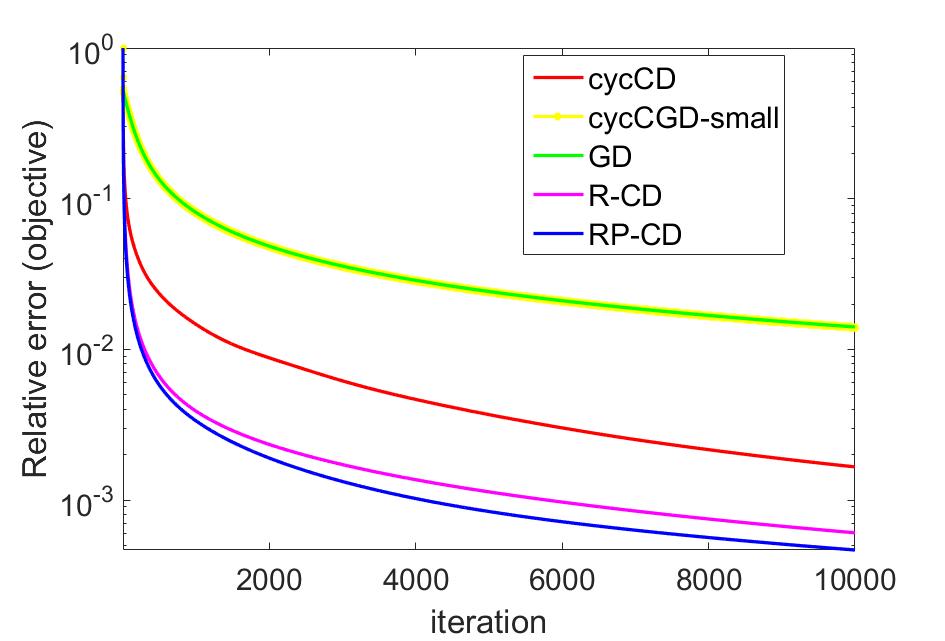}}
\centerline{(a2) Gaussian with mean $2$.  $L \approx  80$, $\chi_{\mathrm{avg}} \approx 79 $}.
\end{minipage}
\vfill
\begin{minipage}{0.48\linewidth}
\centerline{\includegraphics[width=\textwidth]{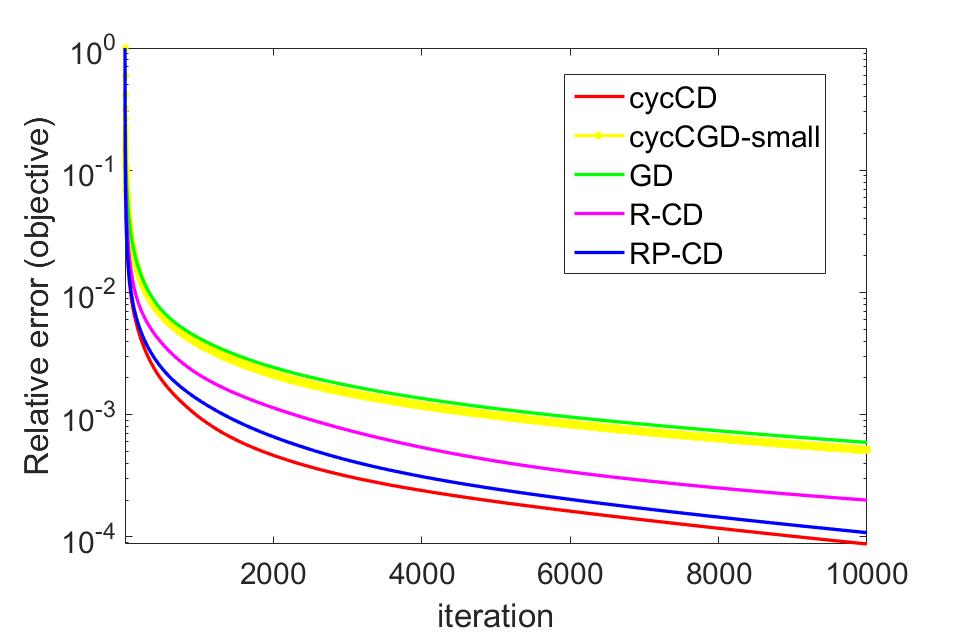}}
\centerline{(b1) Uniform [-0.5, 0.5].  $L \approx  3.8$,  $\chi_{\mathrm{avg}} \approx 7.9 $. }
\end{minipage}
\hfill
\begin{minipage}{0.48\linewidth}
\centerline{\includegraphics[width=\textwidth]{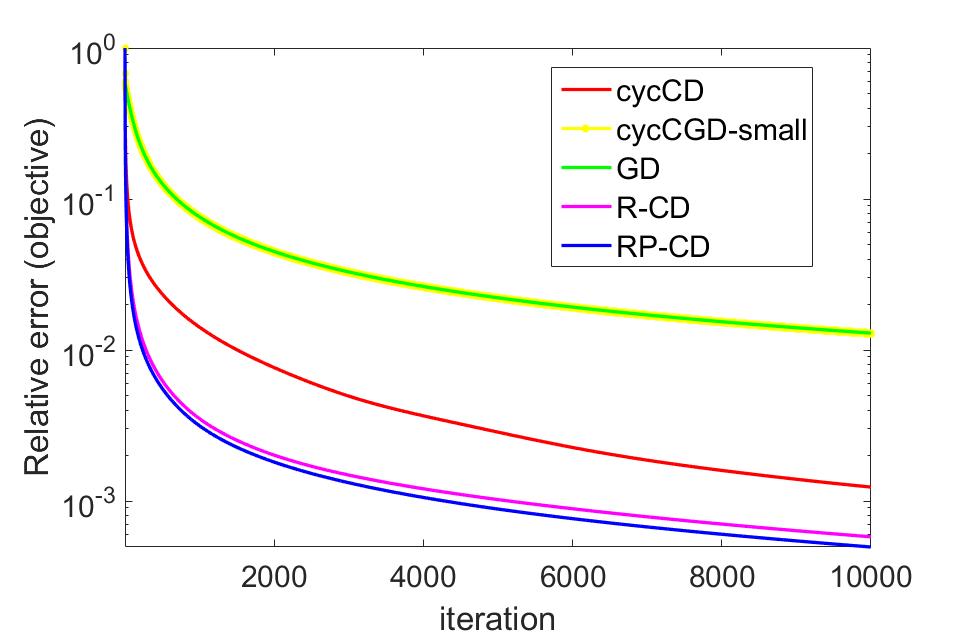}}
\centerline{(b2) Uniform [0,1].  $L \approx 75$, $\chi_{\mathrm{avg}} \approx 74 $. }
\end{minipage}
\begin{minipage}{0.48\linewidth}
\centerline{\includegraphics[width=\textwidth]{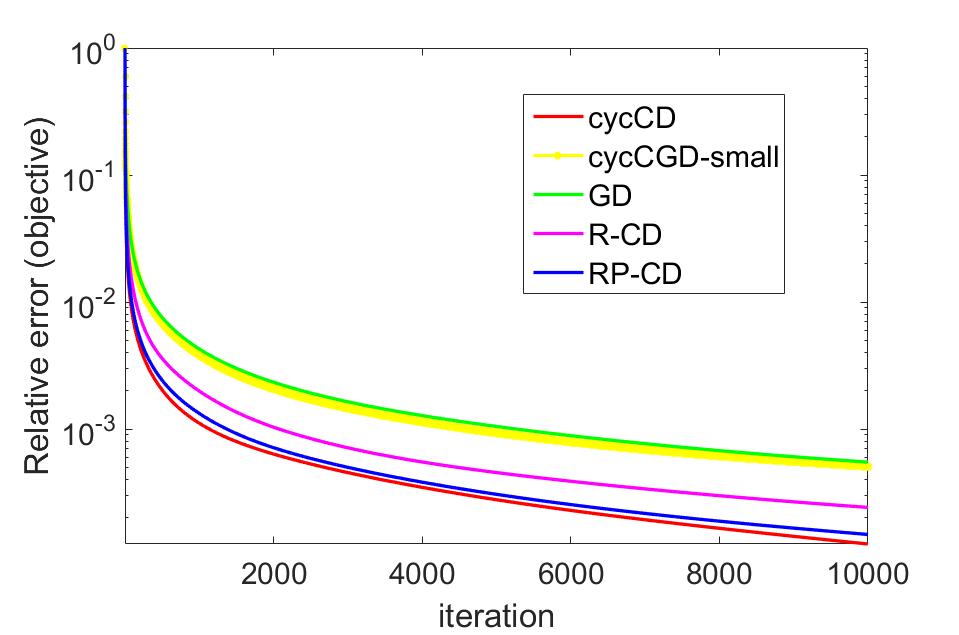}}
\centerline{(c1) Log-normal, with zero mean. $L \approx  3.8$, $\chi_{\mathrm{avg}} \approx 7.6 $. }
\end{minipage}
\hfill
\begin{minipage}{0.48\linewidth}
\centerline{\includegraphics[width=\textwidth]{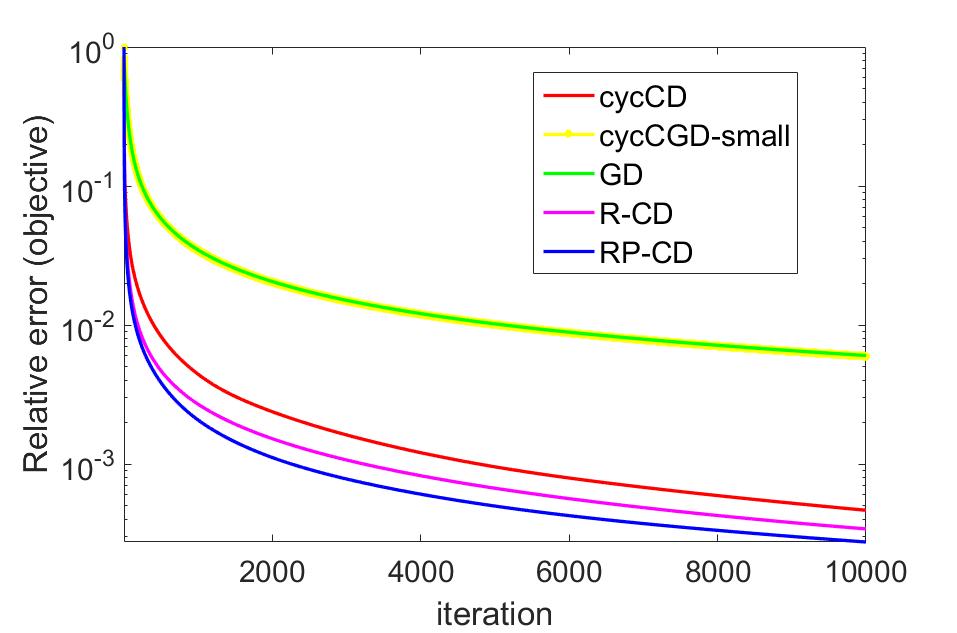}}
\centerline{(c2) Log-normal. $L \approx  42$,  $\chi_{\mathrm{avg}} \approx 41 $.  }
\end{minipage}
%\end{tabular}
\caption{ Comparison of various methods for solving $\min_x x^T A x$, where $A$ is a $ 100 \times 100$ matrix.
$A$ is generated as follows: generate entries of $U$ i.i.d. from a certain distribution and let  $A = U^T U$;
different figures represent different distributions of $U$.
Figures (a1) and (a2): Gaussian distribution with variance $1$; figures (b1) and (b2): uniform distribution;
figure (c1) and (c2): log-normal distribution.
On the left: zero-mean; on the right: non-zero mean.
   }
\label{figureABtest}
\end{figure}

We perform some kind of A/B testing for each distribution: compare the zero-mean case (leading to small off-diagonal entries) with the non-zero mean case (large off-diagonal entries).
We report the simulation results for three distributions Gaussian, uniform and log-normal. 
The simulation results are given in Figure ]\ref{figureABtest}, and the findings from these figures are summarized below.
\begin{enumerate}
	\vspace{-0.15cm}
	\item For all zero-mean cases, C-CD is the fastest; for all non-zero mean cases, C-CD is slower than R-CD/RP-CD.
	% We suspect the performance of C-CD is closely related to $\chi_i$ and $ \chi $.
	This shows that empirically large off-diagonal entries (or large $\tau$) lead to bad performance of C-CD.
	\vspace{-0.15cm}
	\item Different from the example \eqref{Ac def}, C-CD is always much faster than GD in these experiments.
	\vspace{-0.15cm}
	\item Overall, RP-CD is the best algorithm out of the five.
	%  For the zero-mean case, RP-CD is close to C-CD and about twice faster than R-CD; for the non-zero mean case, RP-CD is $\geq 50\%$ faster than R-CD, and several times faster than GD.
	\vspace{-0.15cm}
\end{enumerate}

There are many other ways of generating random $A$. For example, we can multiply $U$ by the square root of a fixed correlation matrix $C$. When $C$ has large off-diagonal entries, the results are similar to those shown on the right column of Figure \ref{figureABtest}.  In statistics, this means that for solving linear regression problems, C-CD is slow when the data have large correlation (\cite{yang2014coordinate} has noticed a related phenomenon).
% We can also consider non-regression setting, letting $U$ be specially structured matrix such as tri-diagonal matrix or Hankel matrix with entries randomly generated. 
% Similar conclusions hold: 
% either when $\tau$ is large we get results similar to the right column of Figure \ref{figureABtest}, or when $\tau $ is small we get the left column.

One interesting question is:
% why is C-CD much faster than GD in these random experiments, while much slower than GD for our example \eqref{Ac def}?
Is randomness crucial in the sense that for any random problem C-CD is faster than GD?
% We do not know the answer, but
It turns out the answer is no.  
We randomly perturb our example, and found that when the perturbation is reasonably small, 
C-CD is still very slow. 
% When the perturbation is large enough, the convergence speed of C-CD will quickly transit from slow to fast. 
%We have investigated a related question:
%is our worst-case example \ref{Ac def} ``stable'', i.e., if we perturb the example randomly, is C-CD still slow?
%It is not hard to imagine that C-CD will still be slow if we randomly perturb our example by a very very small amount; the question is how much perturbation can make C-CD faster than GD.
%It turns out for a reasonable amount of random perturbation C-CD is still slower than GD, and the transition from slow to fast happens rather rapidly.
This also implies that our ``worst-case'' example is robust under small perturbation, which is different from the exponential time example for simplex methods.
Maybe a new type of analysis is needed to explain this phenomenon.

\section{Conclusion}\label{sec: conclulsion}

In this paper, we rigorously establish a $\O(n^2)$ gap between cyclic coordinate descent (C-CD) and randomized coordinate descent (R-CD), when solving quadratic minimization.
More specifically, after presenting an upper bound of $\O(n^4 \kappa_{\text{CD}} \log\frac{1}{\epsilon})$
for C-CD, we prove that this bound is tight in terms of the current parameters.
This is achieved by showing that for a class of examples C-CD does take that many iterations to achieve
accuracy $\epsilon$.
Compared with the complexity of R-CD  $\O(n^2 \kappa_{\text{CD}} \log\frac{1}{\epsilon})$,
our result implies that C-CD can indeed be $O(n^2)$ times slower than R-CD.
When using more parameters such as $\tau = \lambda_{\max}/\lambda_{\avg} $ to characterize the complexity, the complexity of C-CD is approximately $\O(n^2 \tau^2 \kappa_{\text{CD}} \log\frac{1}{\epsilon})$ (up to $\log^2 n$ factor), which is $\O(\tau^2)$ times slower than R-CD.

% Seems UNNECESSARY
%Note that just the numerical experiments for one example is not enough to show
%the gap between the complexity of C-CD  and R-CD. It is not straightforward to interpret the empirical gap.
%In another example (see Example 2 of Section \ref{sec: role of examples}), the gap between R-CD and C-CD is largely because R-CD converges about $15 = \log(1/\epsilon) $ times faster than that predicted by the theoretical bound of R-CD. In addition, that example only shows an $\O(n)$ gap between C-CD and R-CD.
%Even if the example exhibits a gap predicted by theory, it is nontrivial to rigorously prove this gap.
%In our proof, the difficulties include the lack of close-form expressions of the eigenvalues, and the subtle discrepancy between the spectral radius and the lower bound of convergence rate. 

Due to the equivalence of C-CD, Gauss-Seidel method, Kaczmarz method and POCS for solving symmetric PSD linear systems, our result also establishes an $O(n^2)$ gap between the cyclic versions of these methods and their randomized counterparts. 
An interesting finding is that the classical bound of POCS in \cite{smith1977practical} is not better than our bound,
and for the proposed example is infinitely times worse than our bound. 
%The classical convergence rate from POCS literature \cite{smith1977practical}, when specialized to linear systems solving, is dependent on all eigenvalues of the coefficient matrix, instead of the extreme/average eigenvalues  used in our bound. It turns out that the bound in \cite{smith1977practical} is not better than our bound,
%and for the proposed example is infinitely times worse than our bound. 
%It is an interesting open question to obtain a stronger bound for C-CD that depends on all eigenvalues of the matrix.

The simulation partially validates our worst-case analysis. 
% It shows that for the proposed example, there is indeed an $\O(n^2)$ gap between C-CD and R-CD. 
For random coefficient matrices, our numerical experiments show that
the ratio $\tau =  \lambda_{\max}/\lambda_{\avg} $ is closely related to the performance of C-CD. When 
the ratio $\tau$ is large (e.g., in a regression problem with large correlation between the variables), 
C-CD is much slower than R-CD. However, in all random data experiments the gap was never as large as $\O(n^2)$.
More strangely, C-CD is always much faster than GD for random data. 
Thus more theory is needed to explain the worst-case performance of C-CD and typical performance in numerical experiments.

We then discuss some subtle issues and some open questions on the worst-case complexity of C-CD. 
% The parameters used to characterize the complexity play an important role, and can make the complexity comparison rather complicated.  
One subtly arises in the analysis of the non-equal diagonal case, for which we have argued that a more reasonable set of parameters should be based on the Jacobi-preconditioned
version of the original coefficient matrix. 
To perform a comparison with GD and R-CD, we need to explore the relationships between the Jacobi-preconditioned matrix and the original matrix, which is not well understood yet.
This lack of understanding leads to an open question whether an $\O(n^3)$ gap can be established for the non-equal diagonal quadratic case. 
Yet another issue was mentioned in \cite{sun2015improved}: for general convex case (even with equal
per-block Lipschitz constant), a few bounds for C-CD were established but
it is still not known whether he gap between C-CD and R-CD can be $\O(n^3)$. 
It was conjectured in  \cite{sun2015improved} that the current parameters are not enough for characterizing the convergence rate of C-CD for general convex problems. 

A more important open question is whether there is a fundamental gap between deterministic versions of CD
and randomized versions. We have not yet found a deterministic version of CD which can perform as well as R-CD for the proposed example. Either such an example or a proof of the lower bound for \textit{all} deterministic versions of CD would be very interesting. This question is also related to the best complexity of deterministic iterative algorithms for solving symmetric PSD linear systems and positive LP.
A more general version of this question is whether for other algorithms such as ADMM, there is a fundamental gap between all deterministic algorithms and randomized versions.

%Note that this result holds for any fixed cyclic order, not just a particular order. Based on the example, we establish several almost tight complexity bounds of C-CD for quadratic problems. 

\appendix
\newpage

{\huge \textbf{Appendix}}

\section{Proof of Claim \ref{equivalence of Kaczmarz and G-S}}\label{appen: Proof of equivalence of GS and Kaczmarz}

We restate Claim \ref{equivalence of Kaczmarz and G-S} below for readers' convenience:
\begin{claim}\label{restated equivalence of Kaczmarz and G-S}
	Suppose $b \in  \R^{ n \times 1}$, $A = UU^T \in \R^{n \times n} $, where $U \in \R^{n \times n}$ has no zero row. Then  Gauss-Seidel method for solving $A x = b$ is equivalent to Kaczmarz method
	for solving $U y = b$; here, the equivalence means that there is a one-to-one mapping between
	the iterates of the two algorithms. 
\end{claim}

Proof: Suppose $U^T = (u_1, \dots, u_n)$, then $u_j \neq 0, \forall j.$

To solve the linear system $Ax = b$,  the update equation of Gauss-Seidel method can be written as
\begin{equation}\label{cyclic CD update; again}
\begin{split}
& x^{k, j} = x^{k,j-1} - \frac{ A(j,:)x^{k,j-1} - b_j }{ A_{jj} } e_j
 = x^{k,j-1} - \frac{ u_j^T U^T x^{k,j-1} -b_j   }{ \| u_j \|^2 } e_j ,\;, j=1,\dots, n;    \\
 & x^{k+1} = x^{k,n}, \  x^{k+1, 0} = x^{k+1}.
\end{split}
\end{equation}

Let $y^{k,j} = U^T x^{k,j}$ and $y^{k+1} = U^T x^{k+1}, \ y^{k+1, 0} = U^T x^{k+1, 0}$, and
multiply $U^T$ on both sides of the above equations,  we get 
\begin{equation}\label{Kaczmarz update; again}
\begin{split}
& y^{k, j} = y^{k,j-1}  - \frac{ u_j^T y^{k,j-1} -b_j }{ \| u_j \|^2 } u_j, \;, j=1,\dots, n;  \\
& y^{k+1} = y^{k,n}, \  y^{k+1, 0} = y^{k+1}.
\end{split}
\end{equation}

This is exactly the update equation of Kaczmarz method. 
Since $U$ is invertible, define $  x^{k,j} = U^{-T} y^{k,j} $, we can transform the Kaczmarz method to Gauss-Seidel method.     \QED

Remark: The above proof shows that Gauss-Seidel method for any symmetric PSD linear system can be transformed to Kaczmarz method. The other direction is less clear if $U^T$ does not have an inverse. Below we show that in the general case Kaczmarz method is ``almost'' equivalent to Gauss-Seidel method.

To simplify the discussion, we assume $b = 0$ and $ \| u_j\|^2 = 1, \; \forall j. $
The projection onto the hyperplane $\mathcal{H}_j = \{y \mid \langle u_j, y \rangle = b_j \}$ has a simple expression
$I - u_j u_j^T$.  

\textbf{Case 1:} $U $ is square invertible. This means that $u_1, \dots, u_n$ form
a basis of $\R^n$. 
We rewrite the above proof in a more intutive way. 
Any vector $y$ can be represented under the basis $u_1, \dots, u_n$ as
$$
  y = x_1 u_1 + \dots + x_n u_n.
$$
Projecting $y$ onto a hyperplane $\mathcal{H}_1$ is just left multiplying $y$ by $I - u_1 u_1^T$:
\begin{align*}
 (I - u_1 u_1^T) y = (I - u_1 u_1^T ) (x_1 u_1 + \dots + x_n u_n ) 
   = 0 + \sum_{j=2}^n (I - u_1 u_1^T ) x_j u_j   \\
   = \sum_{j=2}^n (u_j - u_1 u_1^T u_j ) x_j
   = - ( u_1^T u_2 x_2 + \dots + u_1^T u_n x_n ) u_1 + x_2 u_2 + \dots + x_n u_n.
\end{align*}
Thus the coordinates $x_2, \dots, x_n$ are unchanged, and the first cordinate $x_1$ is updated to $ - ( u_1^T u_2 x_2 + \dots + u_1^T u_n x_n )$,
which is exactly the optimal solution to  $\min_{x_1}  x^T U U^T x $ with other variables $x_2, \dots, x_n$ fixed.
Under a basis transformation $y = U^T x$, one iterate of Kaczmarz method for updating $y$ is exactly one iterate of Gauss-Seidel method for updating $x$.
Therefore, Gauss-Seidel method is just Kaczmarz method under a different basis. 

\textbf{Case 2:}  Full row-rank linear system, i,e., $ U \in \R^{n \times m}$ with rank $n \leq m$, and the initial point $y^{0,0}$ lies in the span of $u_1, \dots, u_n $.
Note that the row vectors  $u_1, \dots, u_n$ may not span the whole space $\R^m$. 
The  equivalence of G-S method and Kaczmarz method still holds.
 In fact, $y^{0,0} \in \R^m$ can be expressed by $n $ vectors $u_1, \dots, u_n$, and all the iterates stay in the span of $u_1, \dots, u_n$. Thus there is a one-to-one mapping between $y^{k,j}$ and $x^{k,j}$ which is formed by the coordinates of $y^{k,j}$ under the basis $u_1, \dots, u_n$.

\textbf{Case 3:} $U$ is not full row-rank, and the initial point $y^{0,0}$ lies in the span of $u_1, \dots, u_n $.  This includes the overdetermined case $n > m$, as well as
the underdetermined case $n \leq m$ with linearly dependent rows.
There is no one-to-one correspondance between the two methods; nevertheless,
 each sequence of Kaczmarz method corresponds to infinitely many sequences of G-S method.
% Kaczmarz method cannot be interpreted as Gauss-Seidel method by the above procedure. 
This is because $y$ can be represented by $u_1, \dots, u_n$ in multiple ways, i.e., the representation $y = x_1 u_1 + \dots + x_n u_n $ is not unique. 
Fix any representation of the initial point $y^{0,0}$,
the coordinates of $y^{0,0}$ under the spanning set $u_1, \dots, u_n$ can be updated according to the rule described in Case 1, which can be viewed as Gauss-Seidel method.
Thus one representation of $y^{0,0} $ leads to one sequence of Gauss-Seidel iterates.
Different representations of $y^{0,0}$ can lead to different sequences of Gauss-Seidel iterates. 

%One could pick a fixed representation such as $y = x_1 u_1 +  \dots + x_m u_m + 0\cdot u_{m+1} + \dots + 0 \cdot u_n$, and projecting  $y$ onto the hyperplane  $\langle u_j, y \rangle = b_j, j = 1,\dots, m$ can be viewed as updaing the $j$-th coordinate $x_j$; however, when projecting $y$ onto $\langle u_{m+1}, y \rangle = b_{m+1} $, all the coordinates $x_1, \dots, x_m$ may be changed. Thus we cannot interpret Kaczmarz method as Gauss-Seidel method for overdetermined linear system. 

%\begin{equation}\label{cyclic CD update}
%\begin{split}
%& y^{k, j} = \text{Proj}_{\mathcal{H}_j} (y^{k,j-1}) = y^{k,j-1} + \frac{1}{\|u_j \|^2 } (b_j - \langle u_j, x^{k,j-1}  \rangle ) u_j,     \quad j=1,2, \dots, m ; \\
%& y^{k+1} = y^{k,n}, \  y^{k+1, 0} = y^{k+1}.
%\end{split}
%\end{equation}

%\begin{equation}\label{cyclic CD update}
%\begin{split}
%& x^{k, j} = x^{k,j-1} - \frac{ A(j,:)x^{k,j-1} - b_j }{ A_{jj} } e_j,  \quad j=1,2, \dots, n, \\
%& x^{k+1} = x^{k,n}, \  x^{k+1, 0} = x^{k+1}.
%\end{split}
%\end{equation}

%Its normal equation is $U^T U x = U^T \hat{b}$. 
%Define $A = U^T U$, $b = U^T \hat{b}$,
%then the normal equation becomes $Ax = b$. 

\section{Proofs of Upper Bounds}
\subsection{Proof of Proposition \ref{prop: upper bound}}\label{appen: proof of upper bound}
Without loss of generality, we can assume $ b = 0 $.
In fact, minimizing $f(x) = x^T A x - 2 b^T x$ is equivalent to minimizing $ f(x) = (x - x^*)^T A (x - x^*) $
where $x^* = A^{\dag}b $; here we use the fact that $ A x^* = A A^{\dag}b  = b $ when $b \in \mathcal{R}(A)$.
% Since C-CD is invariant with respect to the simultaneous shift of the initial point and the space of variables,
By a linear transformation $z = x-x^*$,
 C-CD for minimizing $ (x - x^*)^T A (x - x^*) $ starting from $x^0$
is equivalent to C-CD for minimizing $ z^T A z $ starting from $z^0 = x^0-x^*$.
Thus we can assume $x^* = 0$, or equivalently, $b = 0$.

The update equation of C-CD now becomes
\begin{equation}\label{update of cyclic CD in proof}
  x^{k+1} = (I - \Gamma^{-1}A) x^k = x^k - d^k,
\end{equation}
where $\Gamma$ is the lower triangular part of $A$ with diagonal entries, i.e., $\Gamma_{ij} = A_{ij}, 1 \leq j \leq i \leq n $, and $d^k = \Gamma^{-1}A x^k $ is the moving direction.
This implies
\begin{equation}\label{d and x relation in C-CD}
\Gamma d^k = A x^k.
\end{equation}

We first assume $A$ is positive definite and will show how to extend to the PSD case in the end.

The proof consists of two main claims.
The first claim relates the convergence rate of C-CD with the spectral radius of a certain matrix.
\begin{claim}\label{claim: C-CD iteration matrix}
  Let $D_A = \text{diag}(A_{11}, \dots, A_{nn})$ be a diagonal matrix with entries $A_{ii}$'s. Then
  \begin{equation}\label{C-CD decrease rate, in spec radius}
    f(x^{k+1}) - f(x^*) \leq \left( 1-  \frac{1}{\| D_A^{-1/2} \Gamma^T A^{-1} \Gamma D_A^{-1/2} \|} \right) (f(x^k) - f(x^*)) .
  \end{equation}
\end{claim}

\emph{First Proof of Claim \ref{claim: C-CD iteration matrix}} (Optimization Perspective):
  Following the proof framework of \cite{sun2015improved}, we bound the descent amount and the cost yet to be minimized (cost-to-go) respectively.
  Suppose $w^0 = x^k, w^n = x^{k+1}$ and $w^1,\dots, w^{n-1}$ are the $n-1$ intermediate iterates.
  Since $w^i$ is obtained by minimizing $f$ over the $i$-th coordinate with other variables fixed,
 it is easy to verify
 \begin{equation}\label{di express as grad}
  d^k_i =\frac{1}{2 A_{ii}} \nabla_i f(w^{i-1}) .
  \end{equation}  % = 2 (e_i^T A w^{i-1}) =
 In the above expression,  $ 2 A_{ii} $ can be viewed as the $i$-th coordinate-wise Lipschitz constant of $\nabla f$ from an optimization perspective.
    We have
   $$ w^1 = w^0 - d^k_1 e_1 = w^0 - \frac{1}{2 A_{11}} \nabla_1 f(w^0), \dots, w^n = w^{n-1} - d^k_n e_n
     = w^{n-1} - \frac{1}{2 A_{nn}} \nabla_n f(w^{n-1} )  , $$  where $e_i$ is the $i$-th unit vector.
  Then
     \begin{equation}\label{one step descent}
  \begin{split}
    f(w^{i-1}) - f(w^{i})  & = (w^{i-1})^T A w^{i-1} - ( w^{i-1} - d^k_i e_i)^T A (w^{i-1} - d^k_i e_i)
      = - d^k_i e_i^T A e_i d^k_i + 2 (w^{i-1})^T A e_i d^k_i  \\
    &  = - A_{ii} (d^k_i)^2 +  \nabla_i f(w^{i-1}) d^k_i
    \overset{ \eqref{di express as grad}}{ =} - A_{ii} (d^k_i)^2 + 2 A_{ii} (d^k_i)^2 = A_{ii} (d^k_i)^2.
    \end{split}
  \end{equation}

  Therefore, the descent amount $ f(x^k) - f(x^{k+1}) $ can be bounded in terms of $d^k$ as
  \begin{equation}\label{sufficient decrease}
   f(x^k) - f(x^{k+1}) = \sum_{i=1}^n f(w^{i-1}) - f(w^{i}) = \sum_{i=1}^n  A_{ii} (d^k_i)^2
    = (d^k)^T D_A d^k.
  \end{equation}

  The cost-to-go estimate is simply
  \begin{equation}\label{cost to go in proof}
   f(x^k) - f(x^*) = f(x^k) = (x^k)^T A x^k
    \overset{\eqref{d and x relation in C-CD}}{=} (d^k)^T \Gamma^T A^{-1} \Gamma d^k   .
    \end{equation}
    Combining with \eqref{sufficient decrease}, we obtain
    \begin{equation}\label{cost-to-go}
     \frac{f(x^k) - f(x^*)}{ f(x^k) - f(x^{k+1})}
      =  \frac{  (d^k)^T \Gamma^T A^{-1} \Gamma d^k }{  (d^k)^T D_A d^k }
      \leq  \| D_A^{-1/2} \Gamma^T A^{-1} \Gamma D_A^{-1/2}     \| ,
     \end{equation}
which implies \eqref{C-CD decrease rate, in spec radius}. % $ f(x^k) - f(x^*) \leq  \| D_A^{-1/2} \Gamma^T A^{-1} \Gamma D_A^{-1/2}     \| $

  \emph{Second Proof of Claim \ref{claim: C-CD iteration matrix}} (Matrix Recursion Perspective):
 One natural idea is to prove $f(x^{k+1}) = M_f f(x^k) $ or $f(x^{k+1}) \leq \| M_f\| f(x^k)$ for a certain matrix $M_f$, based on the update equation of the iterates $x^{k+1} = (I - \Gamma^{-1} A )x^k$. We can write down the expression of $f(x^{k+1})$ in terms of $x^k$ as
  $f(x^{k+1}) = (x^k)^T (I - \Gamma^{-1} A )^T A (I - \Gamma^{-1} A ) x^k  $. 
  However,
  it is not clear how this expression is related to $f(x^k) = (x^k)^T A x^k $.
  A simple trick to resolve this issue is to express everything in terms of $d^k$. More specifically, we have
    \begin{equation}\label{sufficient decrease, another way}
  \begin{split}
    f(x^k) - f(x^{k+1}) = ( x^k)^T A x^k - (x^k - d^k)^T A (x^k - d^k)
     = 2(d^k)^TA x^k -  (d^k)^T A d^k    \\
     = 2 (d^k)^T \Gamma d^k - (d^k)^T A d^k
     = (d^k)^T (\Gamma + \Gamma^T) d^k  - (d^k)^T A d^k
     = (d^k)^T D_A d^k,
     \end{split}
  \end{equation}
  where the last step is because $ \Gamma + \Gamma^T = A + D_A $.
  Equation \eqref{sufficient decrease, another way} is equivalent to \eqref{sufficient decrease} derived earlier using another approach.
  The rest is the same as the first proof.  \QED

Remark: Although the second proof seems simpler, for people who are familiar with optimization the first proof is probably easier to understand:
equation \eqref{one step descent} is just the classical descent lemma (applied to each coordinate), thus \eqref{sufficient decrease} is straightforward to derive.
In the proof of \cite{sun2015improved}, one crucial step is to bound the cost-to-go in terms of $d^k$;
 here for the quadratic case the cost-to-go has a closed form expression given by \eqref{cost-to-go}.
 The second proof is cleaner to write, but it is specifically tailored for the quadratic problem;
 in contrast, the first proof can be extended to non-quadratic problems as done in \cite{sun2015improved} (\eqref{sufficient decrease}
 and \eqref{cost-to-go} will become inequalities).

\begin{claim}\label{claim: bound spectral radius}
  Let $D_A = \text{diag}(A_{11}, \dots, A_{nn})$ be a diagonal matrix with entries $A_{ii}$'s. Then
  \begin{equation}\label{bound G'AinvG}
    \|  \Gamma^T A^{-1} \Gamma \| \leq
     \kappa \cdot \min \left\{   \sum_i L_i, (2 + \frac{1}{\pi} \log n)^2  L    \right\}.
  \end{equation}
\end{claim}

\emph{Proof of Claim \ref{claim: bound spectral radius}:}

 Denote $$ \Gamma_{\mathrm{unit} } = \begin{bmatrix}
   1 &  0 & 0 & \dots & 0  \\
   1 &  1 & 0 & \dots & 0  \\
   \vdots & \vdots & \vdots & \ddots & \vdots  \\
   1 &  1 & 1 & \dots & 0  \\
   1 &  1 & 1 & \dots & 1  \\
     \end{bmatrix}  , $$
     then $ \Gamma = \Gamma_{\mathrm{unit} } \circ A$, where $\circ $ denotes the Hadamard product.
According to the classical result on the operator norm of the triangular truncation operator
\cite[Theorem 1]{angelos1992triangular}, we have
$$
  \| \Gamma  \| = \| \Gamma_{\mathrm{unit}} \circ A \| \leq ( 1 + \frac{1}{\pi} + \frac{1}{\pi} \log n  ) \| A \|
   \leq (2 + \frac{1}{\pi} \log n) \| A \|.
$$
Thus we have
$$ \|  \Gamma^T A^{-1} \Gamma \|
 \leq \| \Gamma^T \Gamma \| \| A^{-1} \|  = \| \Gamma \|^2 \frac{1}{\lambda_{\min}(A) }
  \leq  (2 + \frac{1}{\pi} \log n)^2 \frac{\| A\|^2}{\lambda_{\min}(A)}
  = (2 + \frac{1}{\pi} \log n)^2 \kappa L ,   $$
which proves the second part of \eqref{bound G'AinvG}.

We can bound $ \| \Gamma \|^2 $ in another way (denote $\lambda_i$'s as the eigenvalues of $A$):
\begin{equation}\label{Gamma square bound, another}
\begin{split}
  \| \Gamma \|^2 \leq \|\Gamma \|_F^2 =  \frac{1}{2}( \|A \|_F^2 + \sum_i A_{ii}^2 )
% We could bound the RHS by $\| A\|_F^2$, but this will lose a factor of $\frac{1}{2}$.
% We further bound the RHS while keeping the $\frac{1}{2}$ factor:
% $$  \| \Gamma \|^2
= \frac{1}{2}  \left( \sum_i \lambda_i^2 +  \sum_i A_{ii}^2  \right)   \\
                 \leq \frac{1}{2} \left( (\sum_i \lambda_i)\lambda_{\max} + L_{\max} \sum_i A_{ii}    \right)
                  \overset{\text{(i)}}{=} \frac{1}{2} ( L + L_{\max} )\sum_i L_i \leq L \sum_i L_i.
       \end{split}
\end{equation}
where(i) is because
$ \sum_i \lambda_i = \text{tr}(A) = \sum_i A_{ii} $ and $A_{ii} = L_i $.
Thus
$$
  \|  \Gamma^T A^{-1} \Gamma \|
 \leq  \| \Gamma \|^2 \frac{1}{\lambda_{\min}(A) } \overset{\eqref{Gamma square bound, another}}{\leq} \frac{L}{\lambda_{\min}} \sum_i L_i = \kappa \sum_i L_i.
$$
which proves the first part of \eqref{bound G'AinvG}.
\QED

Finally, according to the fact $ \| D_A^{-1/2} B D_A^{-1/2} \| \leq \frac{1}{\min_i L_i } \| B \| = \frac{1}{L_{\min} } \| B \|  $ for any positive definite matrix $B$, we have
\begin{equation}\nonumber
\begin{split}
 \| D_A^{-1/2} \Gamma^T A^{-1} \Gamma D_A^{-1/2} \|
 \leq \frac{1}{L_{\min} }  \|  \Gamma^T A^{-1} \Gamma \|
 \overset{ \eqref{bound G'AinvG}}{\leq} \frac{1}{L_{\min} }   \kappa \cdot \min \left\{   \sum_i L_i, (2 + \frac{1}{\pi} \log n)^2  L    \right\}
 % \\
% = \kappa \cdot \min \left\{ n \frac{L_{\avg}}{ L_{\min}}, (2 + \frac{1}{\pi} \log n)^2 \frac{L}{L_{\min}}   \right\}.
 \end{split}
\end{equation}
Plugging this inequality into \eqref{C-CD decrease rate, in spec radius} and replacing $\sum_i L_i$ by $n L_{\avg}$, we obtain \eqref{upper bound in kappa}.
%$$
% f(x^{k+1}) - f(x^*) \leq \min \left\{  1 -  \frac{1}{ n \kappa  } \frac{L_{\min}}{L_{\avg}} ,  1 - \frac{L_{\min} }{ L (2 + \log n/ \pi)^2 } \frac{1}{ \kappa  } \right\}  (f(x^k) - f(x^*)) .
%$$
% which is exactly

Now we show how to modify the above proof to the case that $A$ is PSD.
% The intuition is that the update equation \eqref{update of cyclic CD in proof} implies that
% all $\{x^k\}_{k \geq 1}$ lie in the space  $\mathcal{R}(A)$, thus we can view $A$ as a PD matrix when restricted to $\mathcal{R}(A)$. More specifically,
From \eqref{d and x relation in C-CD} we have
$$
  x^k = A^{\dag} \Gamma d^k.
$$
Then \eqref{cost to go in proof} is slightly modified to $ (x^k)^T A x^k
    = (d^k)^T \Gamma^T A^{\dag} \Gamma d^k. $
 We still have \eqref{sufficient decrease}  since its proof does not require $A$ to be positive definite.
 Now we modify \eqref{cost-to-go} to
\begin{equation}\label{PSD case modified ineq}
\begin{split}
  \frac{f(x^k) - f(x^*)}{ f(x^k) - f(x^{k+1})}
  =  \frac{ (d^k)^T \Gamma^T A^{\dag} \Gamma d^k }{  (d^k)^T D_A d^k }
  \leq \frac{ (d^k)^T \Gamma^T \Gamma d^k  \|  A^{\dag}\| }{  (d^k)^T D_A d^k }    \\
  % = \frac{1}{\lambda_{\min}(A) } \frac{ (d^k)^T \Gamma^T \Gamma d^k  }{  (d^k)^T D_A d^k }
  \overset{\text{(i)}} {=}  \frac{1}{\lambda_{\min} } \frac{ (d^k)^T \Gamma^T \Gamma d^k }{  (d^k)^T D_A d^k }
  \leq \frac{1}{\lambda_{\min} } \|D_A^{-1/2} \Gamma^T \Gamma D_A^{-1/2}  \|
  \leq \frac{1}{ \lambda_{\min} L_{\min}  }  \| \Gamma^T \Gamma  \| .
\end{split}
\end{equation}
 where   (i) is because $\|  A^{\dag}\| = 1/\lambda_{\min}$ where $\lambda_{\min}$ is the minimum
 non-zero eigenvalue of $A$.
The rest is almost the same as the proof for the PD case: obtaining the bounds of $  \Gamma^T \Gamma $ as in Claim \ref{claim: bound spectral radius}
and plugging them into \eqref{PSD case modified ineq} immediately leads to \eqref{upper bound in kappa}.

The first bound of result \eqref{upper bound in kappaCD} is a direct corollary of \eqref{upper bound in kappa}
because $ \kappa \leq  n \kappa_{\mathrm CD}$ (which is because $\lambda_{\max}(A) \leq \mathrm{tr}(A) =  n L_{\avg}$).
The second bound of \eqref{upper bound in kappaCD}
 is the same as the second bound of \eqref{upper bound in kappa} because
 $$
    \frac{ \kappa L }{L_{\min}}
    = \frac{ L^2 }{ \lambda_{\min} L_{\min} } = \frac{L^2 }{ L_{\avg} L_{\min} } \frac{ L_{\avg} }{\lambda_{\min}}
     =   \frac{L^2 }{ L_{\avg} L_{\min} }  \kappa_{\mathrm CD}.
 $$
This finishes the proof of Proposition \ref{prop: upper bound}. % \QED
% Applying this relation recursively we obtain

\subsection{Proof of Proposition \ref{prop: upper bound, JacobiPre}}\label{appen: proof of Jacobi Pre upper bound}

This proof is a slight modification of the proof of Proposition \ref{prop: upper bound}.

We first consider the case that $A$ is positive definite.
The insight is to rewrite the relation proved in Claim \ref{claim: C-CD iteration matrix}
\begin{equation}\label{again, C-CD rule}
    f(x^{k+1}) - f(x^*) \leq \left( 1-  \frac{1}{\| D_A^{-1/2} \Gamma^T A^{-1} \Gamma D_A^{-1/2} \|} \right) (f(x^k) - f(x^*))
  \end{equation}
  as
\begin{equation}\label{Jacobi Pre decrease ratio}
   f(x^{k+1}) - f(x^*) \leq \left( 1-  \frac{1}{\| \hat{\Gamma}^T \hat{A}^{-1} \hat{\Gamma} \|} \right) (f(x^k) - f(x^*)) ,
 \end{equation}
  where $ \hat{\Gamma} = D_A^{-1/2} \Gamma D_A^{-1/2}  $ and $ \hat{A} = D_A^{-1/2} A D_A^{-1/2}  $.
  Note that $\hat{\Gamma}$ is still the lower-triangular part (with diagonal entries) of the Jacobi-preconditioned matrix $\hat{A}$.
  The diagonal entries of $\hat{\Gamma}$ and $\hat{A}$ are all $1$, so $\hat{L}_i = 1, \ \forall i$.

  Applying Claim \ref{claim: bound spectral radius} we have
  $$  \|  \hat{ \Gamma}^T \hat{A}^{-1} \hat{ \Gamma } \| \leq
     \hat{ \kappa } \cdot \min \left\{   \sum_i \hat{ L}_i, (2 + \frac{1}{\pi} \log n)^2  \hat{ L}    \right\}
   =  \hat{ \kappa } \cdot  \min \left\{   n , (2 + \frac{1}{\pi} \log n)^2  \hat{ L}    \right\}   . $$
Plugging the above relation into \eqref{Jacobi Pre decrease ratio} we obtain \eqref{Jacobi Pre,upper bound in kappa}.
Similar to Proposition \ref{prop: upper bound}, the second bound \eqref{Jacobi Pre,upper bound in kappaCD}
follows directly from \eqref{Jacobi Pre,upper bound in kappa}.

The case that $ A $ is PSD is can be handled in a similar way to the proof of Proposition \ref{prop: upper bound}.

\section{Supplemental Proofs for Theorem \ref{thm: lower bound}}
\subsection{Proof of Lemma \ref{lemma of equation} }\label{appen: lemma equation proof}
Suppose $\lambda$ is an eigenvalue of $ Z = \Gamma^{-1}A $ and $v = (v_1; v_2; \dots; v_n) \in \C^{n \times 1}$ is the corresponding eigenvector. Then we have
\begin{align}
    & \Gamma^{-1} A v   = \lambda v  \nonumber \\
 \Longrightarrow  \quad &  A v = \lambda L v  \nonumber  \\
 \Longrightarrow  \quad & \begin{cases}   v_1 + c \sum_{j \neq 1} v_j = \lambda v_1 \\
  v_2 + c \sum_{j \neq 2 } v_j = \lambda (c v_1 + v_2 ) \\
   \dots \\
  v_k + c \sum_{j \neq k} v_j  = \lambda (c v_1 + \dots + c v_{k-1} +  v_k), \\
   \dots  \\
  v_n  + c \sum_{j \neq n} v_j = \lambda (c v_1 + \dots + c v_{n-1} +  v_n ).
 \end{cases}   \label{Eig eq. 1st}
\end{align}

Without loss of generality, we can assume
\begin{equation}\label{normalize v}
\sum_{j=1}^n v_j = 1 .
\end{equation}
Let $\hat{c} = 1 - c$. Then \eqref{Eig eq. 1st} becomes
\begin{equation}\label{Eig eq. 2nd}
 \begin{cases}   \hat{c} v_1 + c = \lambda v_1 \\
  \hat{c} v_2 + c = \lambda (c v_1 + v_2 ) \\
   \dots \\
  \hat{c} v_k + c  = \lambda (c v_1 + \dots + c v_{k-1} +  v_k), \\
   \dots  \\
  \hat{c} v_n  + c = \lambda (c v_1 + \dots + c v_{n-1} +  v_n ).
 \end{cases}
\end{equation}
The first equation implies $ v_1 = \frac{c}{ \lambda - \hat{c}}$.
Plugging into the second equation, we get
$$
  v_2 = \frac{c (1 - \lambda v_1 )}{\lambda - \hat{c} } = \frac{c( \lambda - \hat{c} - \lambda c )}{ (\lambda - \hat{c})^2 } = \frac{c \hat{c} (\lambda - 1)}{ ( \lambda - \hat{c})^2} .
$$
Plugging the expression of $v_1, v_2$ into the third equation, we get
$$
  v_3 = \frac{ c(1 - \lambda v_1 - \lambda v_2) }{ \lambda - \hat{c} } = \frac{c (\hat{c})^2 (\lambda - 1)^2 }{ (\lambda - \hat{c})^3 }.
$$
In general, we can prove by induction that
\begin{equation}\label{expression of v}
   v_k = \frac{c (\hat{c})^{k-1} ( \lambda -1 )^{k-1} }{ (\lambda - \hat{c})^k } = \frac{c}{\lambda - \hat{c}} q^{k-1},
\end{equation}
where
\begin{equation}\label{expression of q}
  q = \frac{\hat{c}(\lambda - 1)}{\lambda - \hat{c}}  .
\end{equation}
We can also express $\lambda$ in terms of $q$ as
\begin{equation}\label{expression of lambda}
  \lambda = \frac{ \hat{c} - \hat{c}q }{\hat{c}-q} .
\end{equation}

Note that the expression of $v_k$ given by \eqref{expression of v} satisfies \eqref{Eig eq. 2nd}
for any $\lambda$, but our goal is to compute $\lambda$. To do this,
we need to utilize the normalization assmption \eqref{normalize v}. In particular, we have (when $q \neq 1$)
\begin{align}
 &  1 = \sum_k v_k  = (\sum_{k=1}^n q^{k-1}) \frac{c }{\lambda - \hat{c}}
   = \frac{1 - q^n }{1 - q} \frac{c }{ \lambda - \hat{c} }  \nonumber  \\
  \Longrightarrow \quad &   (1-q ) (\lambda - \hat{c}) = c(1-q^n)   \nonumber  \\
  \overset{ \eqref{expression of q} }{\Longrightarrow} \quad &  c \lambda = c (1 - q^n) \nonumber \\
   \overset{ }{\Longrightarrow} \quad &  q^n = 1 - \lambda \overset{ \eqref{expression of lambda} }{=} 1 -  \frac{ \hat{c} - \hat{c}q }{\hat{c}-q}  \label{lambda to n and 1 -q  relation}    \\
   \overset{ }{\Longrightarrow} \quad & q^n = \frac{qc}{q - \hat{c}}  \nonumber   \\
   \overset{ }{\Longrightarrow} \quad & q^{n} (q - \hat{c}) = c q.   \nonumber
\end{align}

The above procedure is reversible, i.e. suppose $q \neq 1$ is a root of $ q^{n} (q - \hat{c}) = c q $,
then $\lambda = \frac{ \hat{c} - \hat{c}q }{\hat{c}-q}  $ is an eigenvalue of $ Z¡¡$.
Suppose the  $ n + 1 $ roots of $ q^{n} (q - \hat{c}) = c q $ are $q_0 = 0, q_1, \dots, q_{n-1}, q_n = 1$ ($q = 0$ and $q=1$ are always roots),
then $ \lambda_k = \frac{ \hat{c} - \hat{c} q_k }{ \hat{c}- q_k } \overset{ \eqref{lambda to n and 1 -q  relation}}{=} 1 - q_k^n , k= 0, \dots, n - 1$ are the $n$ eigenvalues of $Z$.

\subsection{Proof of Lemma \ref{claim of n roots convergence}}\label{appen: Lemma of cts roots proof}
The roots of a polynomial continuously depend on the coefficients of the polynomial, and thus the roots of a series of polynomials will converge to the roots of tbe limiting polynomial of this family; see \cite[Theorem 4A]{whitney1972complex}. To make our proof self-consistent, 
we will prove Lemma \ref{claim of n roots convergence} by Rouch\'{e}'s theorem in complex analysis. 

When $n=1 $, the only solution of $  q^{n-1} (q - 1 + c) = c $ is $q = 1$, thus the conclusion holds. 
From now on, we assume $n \geq 2$.

Let $p = 1/q$, then the equation $ q^{n-1} (q - 1 + c) = c $ becomes
\begin{align*}
p^{-1} -1 + c = c p^{n-1}  \Longleftrightarrow 
1 - (1-c) p = c p^n   & \Longleftrightarrow    1/c - (1/c- 1) p =  p^n  \\
& \Longleftrightarrow 
p^n -1 + (1/c - 1) (p - 1) = 0.
\end{align*}
This equation can be written as $F(p) + G_c(p) = 0$, where $F$ and $G_c$ are defined as
$$
F (p) =  p^n -1,   \;  G_c(p) = (1/c - 1) (p - 1) . 
$$

\begin{lemma}\label{lemma of roots converge}
	Suppose $n \geq 2$.
	For any $ 0< \epsilon <  \sin ( \pi/n)   $, there exists some
	$\delta > 0  $ such that for any $c \in (1 - \delta, 1 )$,
	$F(p)  + G_c(p) $ has exactly one root $p_k$ in the ball   $ B( e^{- i 2k \pi /n}, \epsilon  )  \triangleq \{ z \mid | z - e^{- i 2k \pi /n} | \leq   \epsilon  \} $, $k = 0, 1, \dots, n - 1$.
\end{lemma}

Clearly, the function $F $ has $n$ roots $ \eta_k \triangleq  e^{- i 2k \pi /n}, k=1, \dots, n $. 
The distance between two adjacent roots are 
$$
| 1 - e^{-2 i \pi /n} |  = 2 \sin ( \pi/n) .    % |1 - \cos (2 \pi/n)   | 
$$
For any $0< \epsilon <   \sin ( \pi/n)   $, consider $n$ balls
$$
B ( \eta_k ,  \epsilon  ) =   \{ z \mid | z - \eta_k | \leq   \epsilon  \}, k= 0, 1, \dots, n-1.
$$
Any two such balls have no intersection since $\epsilon < | \sin ( \pi/n) | = \min_{0 \leq j, k \leq n -1} |\eta_j - \eta_k| $. 

The boundary of the ball $B ( \eta_k ,  \epsilon  )$ is 
$$ \partial B ( \eta_k ,  \epsilon  )  =   \{ z \mid | z - \eta_k | =    \epsilon  \}.  $$

%To apply Roche's theorem, we need to show that $|F(z)| > |G_c(z)|$ on the boundary 
%$ \partial B ( \eta_k ,  \epsilon  )  $. 
Define 
$$ v_k(\epsilon)  \triangleq \inf_{ z \in \partial B ( \eta_k ,  \epsilon  ) }  F(z)
= \min_{ z \in \partial B ( \eta_k ,  \epsilon  ) }  | z^n - 1   |  > 0 . $$
This minimum can be achieved because $v_k(\epsilon) $ is the minimal value of a continuous function on a compact set. It is positive since otherwise there exists some $z \in \partial B ( \eta_k ,  \epsilon  ) $ such that $z^n = 1$ which means $z \in \{ \eta_0, \dots, \eta_{n-1}  \}$. This contradicts the fact that any two balls 
$B ( \eta_j ,  \epsilon  ), B ( \eta_k ,  \epsilon  )$ have no intersection. 

Define
$$
v(\epsilon ) = \min_{ 0 \leq k \leq n-1  }  v_k(\epsilon) > 0. 
$$
For any $ z \in \partial B ( \eta_k ,  \epsilon  ) $, we have
\begin{equation}\label{F1 lower bound}
| F (z) | = | z^n - 1 |   \geq  v( \epsilon ).
\end{equation}
For any $ z \in \partial B ( \eta_k ,  \epsilon  ) $ and any $c > \frac{3 }{ 3 + v(\epsilon)} $, we have
\begin{equation}\label{F2 upper bound}
| G_c (z) |  =   |  (1/c - 1) ( z - 1) | \leq  |1/c -1| ( | \eta_k | + \epsilon + 1   ) \leq 3 |1/c -1| <  v( \epsilon ),
\end{equation}
where the  second inequality is due to $|\eta_k | = 1$ and $\epsilon < \sin ( \pi/n) \leq 1$.

Combining the two bounds \eqref{F1 lower bound} and \eqref{F2 upper bound}, we obtain that 
$$
| F(z) |   > | G_c (z) | , \forall z \in \partial B ( \eta_k ,  \epsilon  ). 
$$
According to Rouch\'{e}'s theorem, $F $ and $F + G_c $ have the same number of zeros inside $B ( \eta_k ,  \epsilon  )  $.
Since $F $ has exactly one root inside $B ( \eta_k ,  \epsilon  )  $ which is $\eta_k$, we obtain that $F + G_c $ has exactly one root $p_k$ inside
$B ( \eta_k ,  \epsilon  )  $. $\Box$

We first let $ \epsilon_0 =  \sin( \pi/n)/ 4 $, which implies $B ( \eta_k ,  \epsilon_0  ), k=0, 1, \dots, n-1$ are $n$ disjoint balls.
For any $c \in ( 3/(3 + v( \epsilon_0 ) , 1 ) $,
Lemma \ref{lemma of roots converge} implies that  
$F (p) + G_c (p) $ has exactly one root inside each ball.
We denote $ p_0(c), p_1(c), \dots, p_{n-1}(c) $ to be the roots of
$F (p) + G_c (p) $ such that $p_k(c) \in  B ( \eta_k ,  \epsilon_0  )  , \forall k $.
Since $F  + G_c $ has exactly $n$ complex roots, thus $p_k(c)$'s are all the roots of $F + G_c $. 
Lemma \ref{lemma of roots converge} implies that for any $\epsilon > 0 $, there exists some $ \delta   $ such that whenever $ c > 1- \delta $, we have
$$
|   p_k( c) - \eta_k | < \epsilon, \; \forall k.
$$ 
This means 
$$
\lim_{c \rightarrow 1} p_k(c) = \eta_k, \; \forall k. 
$$

Since there is a one-to-one mapping between the roots of
$ q^{n-1} (q - 1 + c) - c $ and the roots of 
$ F(p)  + G_c (p) =  p^n -1 + (1/c - 1) (p - 1)   $ by the inverse transformation 
$ p = 1/q$,   we obtain the following result: for any $c \in ( 3/(3 + v(  \epsilon_0  ) , 1 ) $,  the equation
$  q^{n-1} (q - 1 + c) - c $ has exactly one root $q_k(c)$ such that
$ |1/ q_k(c) - e^{ - i 2\pi k/n } |  < \sin( \pi/n)/2 $ for $k = 0, 1, \dots, n-1$; moreover, 
$$
\lim_{c \rightarrow 1} q_k(c) = e^{ i 2\pi k/n }, \; \forall k. 
$$

Since $ |1/ q_k(c) - e^{ - i 2\pi k/n } |  < \sin( \pi/n)/4 $ implies $ | q_k(c) - e^{  i 2\pi k/n } | < \sin(\pi / n)/2   $, 
we obtain the following conclusion:
for any $c \in ( 3/(3 + v(  \epsilon_0  ) , 1 ) $,  the equation
$  q^{n-1} (q - 1 + c) - c $ has exactly one root $q_k(c)$ such that
$ |  q_k(c) - e^{  i 2\pi k/n } |  < \sin( \pi/n)/2 $ for $k = 0, 1, \dots, n-1$; moreover, 
$$
\lim_{c \rightarrow 1} q_k(c) = e^{ i 2\pi k/n }, \; \forall k. 
$$

\subsection{ Proof of Claim \ref{claim: sum cos equality} }\label{appen: claim sum sin proof}
Since $2 \sin(n\phi/2) \cos( x + (n+1)\phi/2 ) =  \sin( z + (n+ 1/2) \phi ) - \sin( z + \phi/2 )$,
the desired equation \eqref{sum cos, simpler} is equivalent to
\begin{equation}\label{sum cos}
  \sum_{j=1}^n \cos( z + j \phi ) = \frac{ \sin( z + (n+ 1/2) \phi ) - \sin( z + \phi/2 ) }{2 \sin (\phi/2 )}.
\end{equation}
We prove \eqref{sum cos} by induction. When $ n = 1$, it holds because
$ \sin( z + 1.5 \phi) - \sin(z + 0.5 \phi) = 2 \sin(\phi/2) \cos( z + \phi) $.
Suppose \eqref{sum cos} holds for $n-1$, i.e.
$$  \sum_{j=1}^{n-1} \cos( z + j \phi ) = \frac{ \sin( z + (n - 1/2) \phi ) - \sin( z + \phi/2 ) }{2 \sin (\phi/2 )}  . $$
Note that $ 2 \cos(z + n \phi) \sin (\phi/2 ) = \sin(z + (n+1/2 \phi)) -  \sin( z + (n - 1/2) \phi ) $,
therefore
\begin{equation}
\begin{split}
 \sum_{j=1}^{n} \cos( z + j \phi ) & = \frac{ \sin( z + (n - 1/2) \phi ) - \sin( z + \phi/2 ) }{2 \sin (\phi/2 )}
+ \cos(z + n \phi)  \\
 & =   \frac{ \sin( z + (n - 1/2) \phi ) - \sin( z + \phi/2 ) + \sin(z + (n+1/2 \phi)) -  \sin( z + (n - 1/2) \phi ) }{2 \sin (\phi/2 )}   \\
 & =  \frac{ \sin( z + (n+ 1/2) \phi ) - \sin( z + \phi/2 ) }{2 \sin (\phi/2 )}.
 \end{split}
 \end{equation}
This completes the induction step, and thus \eqref{sum cos} holds. \QED

\section{Proofs of Propositions on Exact Comparison}
\subsection{Proof of Proposition \ref{prop: real slower, objective}}\label{appen: prop objective error compare proof}
Same as the proof of Theorem  \ref{thm: lower bound}, we pick $A = A_c$ and consider minimizing $f(x) = x^T A_c x$.
Obviously the minimizer $x^* = 0$ and the optimal value $f^* = f(x^*) = 0$.

We first compute $ k_{\mathrm{GD} }(\epsilon ) $.
The update equation of GD is $ x^k =  (I - \frac{1}{\beta} A)x^{k-1} $, where
$\beta = \lambda_{\max}(A) $.
Since $A$ is a symmetric positive definite matrix, we can assume $A = U^T U $, where $U \in \R^{n \times n}$ is non-singular.
Then
\begin{equation}\label{Uxk GD update}
U x^k = U ( I - \frac{1}{\beta} A ) x^{k-1}  = (I - \frac{1}{\beta} U U^T ) Ux^{k-1}  ,
\end{equation}
 The spectral norm of the iteration matrix
 \begin{equation}\label{spectral norm of GD}
 \| I - \frac{1}{\beta} U U^T \| = \| I - \frac{1}{\beta} A \|  = 1 - \frac{1}{\beta} \lambda_{\min}(A) = 1 - 1/ \kappa ,
 \end{equation}
  where $\kappa$ is the condition number of $A$ given by
 \begin{equation}\label{kappa expresssion, again}
 \kappa = \frac{ 1 -c + cn }{ 1-c } .
 \end{equation}
The relation \eqref{Uxk GD update} implies
$$
  f(x^k) = \| Ux^k \|^2 \leq \| I - \frac{1}{\beta} U U^T \|^2 \| Ux^{k-1} \|^2 .
$$
Therefore we have
 \begin{equation}\label{GD compare}
f(x^k) \leq \| I - \frac{1}{\beta} U U^T \|^{2k} f(x^0)  = \left( 1 - \frac{1}{\kappa} \right)^{2k}  f(x^0) .
\end{equation}

% The minimum number of iterations
The minimum number of iterations to achieve $\frac{ (x^k)^H A x^k }{ (x^0)^H A x^0  } \leq \epsilon $ for all initial points $x^0 \in \R^{n \times 1}$ can be upper bounded as
%(to be rigorous, we should round $K_{\mathrm{GD} }(\epsilon ) $ to an integer, but for simplicity we will just ignore this slight difference)
\begin{equation}\label{k GD expression}
  k_{\mathrm{GD} }(\epsilon )  \leq \left \lceil \frac{1}{2} \frac{\ln \epsilon  }{ \| I - A/\beta \| } \right\rceil \leq 
   \frac{1}{2} \frac{\ln \epsilon }{ \ln (1 - 1/\kappa ) } + 1  \triangleq  \tilde{k}_{\mathrm{GD} }(\epsilon ).
\end{equation}

Let $y^k = U x^k$. We will use the same definitions of $q, r, \lambda_i$ as in the proof of Theorem \ref{thm: lower bound}.
According to \eqref{yk ratio closed form}, to obtain a relative error $ \frac{f(x^k)}{f(x^0)} = \frac{\|y^k\|^2}{\| y^0 \|^2 } \leq \epsilon,$ the number of iterations $k = k_{\mathrm{CCD}}(\epsilon)$ should satisfy
$$
  \epsilon \geq \beta_c r^{(2k+2)n}  \overset{\eqref{r,theta def}}{=} \beta_c |q|^{(2k+2)n}  \overset{ \eqref{lambda k def} }{=} \beta_c  |1 - \lambda_1|^{2k + 2} ,
$$
i.e.
\begin{equation}\label{kCD bound}
   k_{\mathrm{CCD}}(\epsilon) \geq \frac{1}{2} \frac{ \ln(1/\epsilon) + \ln(\beta_c) }{ \ln( 1/ |1 - \lambda_1|) } - 1.
\end{equation}

Since $ \lim_{c \rightarrow 1} \ln( 1/ |1 - \lambda_1 | ) = 0 $ and by \eqref{coeff beta goes to 1} $\lim_{c\rightarrow 1} \beta_c = 1$, we have
$$
 \lim_{c \rightarrow 1}  \frac{ \frac{ \ln(1/\epsilon) + \ln(\beta_c) }{ \ln( 1/| 1 - \lambda_1 | ) } - 1 }{  \frac{ \ln(1/\epsilon)  }{ \ln (1/| 1 - \lambda_1) | }  }
  = 1.
$$

\begin{equation}\label{bound of kCD over kGD}
\begin{split}
  \lim_{c \rightarrow 1 } \frac{ k_{\mathrm{CCD} }(\epsilon ) }{ k_{\mathrm{GD} }(\epsilon ) }
  \geq  \lim_{c \rightarrow 1 } \frac{ k_{\mathrm{CCD} }(\epsilon ) }{  \tilde{k}_{\mathrm{GD} }(\epsilon ) }
  & \geq  \lim_{c \rightarrow 1 } \left( \frac{1}{2} \frac{ \ln(1/\epsilon) + \ln(\beta_c) }{ \ln( 1/ |1 - \lambda_1|) } - 1 \right)
      \left( \frac{1}{2} \frac{\ln \epsilon }{ \ln (1 - 1/\kappa ) } + 1 \right)^{-1}   \\
  & = \lim_{c \rightarrow 1 }   \frac{1}{2} \frac{ \ln(1/\epsilon)  }{ \ln( 1/ |1 - \lambda_1|) }
    \left( \frac{1}{2} \frac{\ln \epsilon }{ \ln (1 - 1/\kappa ) }  \right)^{-1} \\
 & =  \lim_{c  \rightarrow 1 }  \frac{ \ln (1 - 1/\kappa )  }{  \ln  |1 - \lambda_1|  }
    =  \lim_{c  \rightarrow 1 }  \frac{ \ln (1 - 1/\kappa ) }{ - 1/ \kappa } \cdot \frac{ -( 1 -  |1 - \lambda_1 |) }{ \ln |1 - \lambda_1| } \cdot \frac{ 1/ \kappa }{ 1 -  |1 - \lambda_1|  }    \\
  &  =  \lim_{c  \rightarrow 1 }   \frac{ 1/ \kappa }{ 1 -  |1 - \lambda_1 |  }
    \overset{\eqref{Jk expression}}{ = }  \frac{1}{ 2 n \sin^2(\pi/n) } > \frac{n}{2 \pi^2 }.
  \end{split}
\end{equation}

The convergence rate of the objective values for R-CD has been given in \cite[Theorem 2]{nestrov12} and \cite[Theorem 3.6]{leventhal2010randomized}.
% but we did not notice a result on the convergence rate of the iterates.
We present the convergence rate of both the iterates and the objective values for R-CD, when solving quadratic problems \eqref{quadratic min}. The proof is quite straightforward and omitted here. Note that the proposition implies $\| E (x^k ) \|^2 $ converges twice as fast as $E( f(x^k)) $,
which explains why in Proposition \ref{prop: slower, iterates} the gap between C-CD and R-CD is twice as large
as that in Proposition \ref{prop: real slower, objective}.
 % given in the appendix.
\begin{proposition}\label{prop: RCD rate}
Consider solving a quadratic minimization problem \eqref{quadratic min} where $A$ is a positive definite matrix
with all diagonal entries being $1$. Suppose R-CD generates a sequence $z^k$ according to \eqref{R-CD} and define $ x^k = z^{kn} $.
Then
\begin{equation}
 \| E (x^k ) \|^2  \leq \left( 1 - \frac{1}{n} \lambda_{\min}  \right)^{2kn} \| x^0\|^2,
\end{equation}  %  2\lambda_{\min} - \lambda_{\min}^2, if changed to E(x^k)^2
and
\begin{equation}
  E( f(x^k)) \leq (1 - \frac{1}{n}\lambda_{\min} )^{kn} f(x^0) ,
\end{equation}
where $\lambda_{\min} $ is the minimum eigenvalue of $A$.
\end{proposition}

%Remark: For simplicity we assume that the diagonal entries of $A$ are all $1$ in the above result.
%When the diagonal entries of $A$ are not all $1$, one may solve a new problem by replacing $A$ with $D^{-1/2}A D^{1/2}$ to make the diagonal entries all $1$,
%where $D$ is a diagonal matrix with entries being the diagonal entries of $A$.
%If one does not want to execute such a preprocessing, the sampling probabilities should depend on the diagonal entries of $A$;
%see, e.g. \cite{nestrov12}.
According to Proposition \ref{prop: RCD rate}, % \cite[Theorem 2]{nestrov12},
$$
  E( f(x^k) ) \leq \left( 1 - 1/\kappa_{\mathrm{CD}}   \right)^{k n} f(x^0) ,
$$
where $ \kappa_{\mathrm{CD} } = \frac{ \max_i A_{ii} }{ \lambda_{\min}(A) }$.
To achieve an error $ \frac{ E( f(x^k) )}{ f(x^0) } \leq \epsilon $, we only need
$$
  \left( 1 - 1/\kappa_{\mathrm{CD}}  \right)^{k n}  \leq \epsilon  \Longleftrightarrow
  k \geq \frac{1}{n} \frac{ \ln \epsilon }{ \ln ( 1- 1/\kappa_{\mathrm{CD} } )  }    .
$$
% For $A = A_c$, we have $\beta_j = 1, \forall i$ and $ \lambda_{\min}(A) = 1 - c  $,
Therefore we have
\begin{equation}\label{kRCD expression}
  k_{\mathrm{RCD}} (\epsilon)  \leq  \frac{1}{n} \frac{ \ln \epsilon }{ \ln ( 1- 1/\kappa_{\mathrm{CD} } ) }  + 1  .
       %              =  n \frac{ \ln \epsilon }{ \ln \left( 1- \frac{ 1 - c }{ n  } \right) }  + 1 .
\end{equation}

Combining the above relation with \eqref{k GD expression}, we have
  \begin{equation}\label{compare GD with RCD}
  \lim_{c \rightarrow 1 } \frac{   \tilde{k}_{\mathrm{GD} }(\epsilon ) }{  k_{\mathrm{RCD} }(\epsilon ) }
  \geq  \lim_{c \rightarrow 1 }  \frac{ n \ln \left( 1-  1/\kappa_{\mathrm{CD} } \right) }{ 2 \ln (1 - 1/\kappa )}
  =  \lim_{c \rightarrow 1 } \frac{ n/\kappa_{\mathrm{CD} } }{ 2/\kappa}
  = \lim_{c \rightarrow 1 } \frac{ n \lambda_{\max}(A) }{ 2 \sum_i \beta_i }
  =  \lim_{c \rightarrow 1 } \frac{ n( 1 - c + c n) }{ 2 n  } = \frac{n}{2}.
     \end{equation}

Combining the above relation with \eqref{bound of kCD over kGD}, we obtain
 \begin{equation}\label{bound of kCD over kRCD}
   \lim_{c \rightarrow 1 } \frac{  k_{\mathrm{CCD} }(\epsilon ) }{  k_{\mathrm{RCD} }(\epsilon ) }
   > \frac{n}{2 \pi^2 } \frac{n}{2} = \frac{n^2 }{ 4 \pi^2 }.
\end{equation}

According to \eqref{bound of kCD over kGD} and \eqref{bound of kCD over kRCD}, there exists $c$ such that
\eqref{compare GD with CD, coro} and \eqref{compare CD with RCD, coro} hold.  \QED

\subsection{Proof of Proposition \ref{prop: slower, iterates}}\label{appen: prop of iterates compare proof}

 Same as the proof of Theorem  \ref{thm: lower bound}, we pick $A = A_c$ and consider minimizing $f(x) = x^T A_c x$.
Obviously the minimizer $x^* = 0$ and the optimal value $f^* = f(x^*) = 0$.

First we consider $ K_{\mathrm{GD} }(\epsilon ) $.
 Since $A$ is a symmetric positive definite matrix, we can assume $A = U^T U $, where $U \in \R^{n \times n}$ is non-singular.
 The update formula of GD is $x^{k+1} = (I - \frac{1}{\beta} A  ) $. The iteration matrix $I - \frac{1}{\beta} A $
  has the same eigenvalues as $I - \frac{1}{\beta}UU^T $,  the iteration matrix of $\{Ux^k \}$ (see \eqref{Uxk GD update}).
  Since both  $I - \frac{1}{\beta} A $ and  $I - \frac{1}{\beta}UU^T $ are symmetric,
   we have the following relation (which means that for GD the squared iterates and the function values converge at the same speed)
\begin{equation}\label{kGD and KGD the same}
 k_{\mathrm GD}(\epsilon) = K_{\mathrm GD}(\epsilon).
\end{equation}

We then consider  $ K_{\mathrm{CCD} }(\epsilon ) $.
 Compare \ref{yk ratio closed form} with \ref{ratio of xk x0 bound} in the proof of Theorem  \ref{thm: lower bound}, we know that
 the bound we obtained for the function values is the same as the bound for the squared iterates.
 Similar to \eqref{kCD bound}, we have
 $$
 K_{\mathrm{CCD}}(\epsilon) \geq \frac{1}{2} \frac{ \ln(1/\epsilon) + \ln(\omega_c) }{ \ln( 1/ |1 - \lambda_1|) } - 1.
 $$
 Similar to \eqref{bound of kCD over kGD} in the proof of Proposition \ref{prop: real slower, objective}, we have
 \begin{equation}\label{asym ite, GD ratio}
 ¡¡\lim_{c \rightarrow 1 } \frac{ K_{\mathrm{CCD} }(\epsilon ) }{ K_{\mathrm{GD} }(\epsilon ) } > \frac{n }{2 \pi^2 }¡¡.
\end{equation}

Next, we consider $ K_{\mathrm{RCD} }(\epsilon ) $.
According to Proposition \ref{prop: RCD rate},
$$
 \| E(x^k) \| \leq (1 - 1/\kappa_{\mathrm CD})^{kn} \| x^0 \|,
$$
which implies
$$
  K_{\mathrm RCD}(\epsilon) \leq \frac{1}{2n} \frac{ \ln(\epsilon)}{ \ln( 1 - 1/\kappa_{\mathrm CD}) } + 1.
$$
Note that the RHS (right-hand side) of the above bound is asymptotically half the RHS of \eqref{kRCD expression}.
Combining with \eqref{compare GD with RCD} and \eqref{kGD and KGD the same}, we have
$$
 \lim_{c\rightarrow 1} \frac{K_{\mathrm GD}(\epsilon) }{ K_{\mathrm RCD}(\epsilon) } \overset{\eqref{kGD and KGD the same}}{=}
  \lim_{c\rightarrow 1} \frac{k_{\mathrm GD}(\epsilon) }{ K_{\mathrm RCD}(\epsilon) }
  = 2 \lim_{c\rightarrow 1} \frac{k_{\mathrm GD}(\epsilon) }{ k_{\mathrm RCD}(\epsilon) } \overset{\eqref{compare GD with RCD}}{ >}  n .
$$
Multiplying this inequality with \eqref{asym ite, GD ratio}, we have
 \begin{equation}\label{asym ite, RCD ratio}
   \lim_{c \rightarrow 1 } \frac{ K_{\mathrm{CCD} }(\epsilon ) }{ K_{\mathrm{RCD} }(\epsilon ) } > \frac{n^2 }{2 \pi^2 }¡¡.
 \end{equation}

 Finally, we compute $ K_{\mathrm{RPCD} }(\epsilon ) $.
 \begin{claim}\label{claim: RP-CD rate for Ac}
   Consider using RP-CD (randomly permuted coordinate descent) to solve the problem $ \min_{x \in \R^n} x^T A_c x$
   with $A_c$ given in \eqref{Ac def}. Suppose the initial point is $x^0$, then we have
    \begin{equation}\label{RP-CD, converge rate for Ac}
     \| E(x^k) \|^2  \leq \left( 1 - (1-c)(1-\gamma) \right)^{2k} \| x^{0}\|^2 ,
      \end{equation}
      where
      \begin{equation}\label{gamma def, first}
      \gamma = \frac{ -n + (1 - (1-c)^n)/c }{ n (n-1) }.
      \end{equation}
 \end{claim}
 The proof of this claim is given in Appendix \ref{appen: RP-CD rate proof}.

 By the definition of $\gamma$ in \eqref{gamma def, first} we have
 \begin{equation}\label{gamma limit}
 \lim_{c \rightarrow 1} \gamma = - 1/n.
 \end{equation}
 Similar to the proof of Proposition \ref{prop: real slower, objective},
 from \eqref{RP-CD, converge rate for Ac} and \eqref{GD compare} we have
 \begin{align}
   \lim_{c \rightarrow 1 } \frac{  k_{\mathrm{GD} }(\epsilon )}{ K_{\mathrm{RPCD} }(\epsilon )  }
     =  \lim_{c \rightarrow 1 }  \frac{ \kappa }{1/[(1-c)(1-\gamma)] }
      \overset{ \eqref{kappa expresssion, again} }{ = } \lim_{c \rightarrow 1 }  \frac{1-c + cn}{ 1-c } (1-c)(1-\gamma) \nonumber  \\
      \overset{ \eqref{gamma def, first} }{= } \lim_{c \rightarrow 1 } ( 1 -c + cn  )( 1- \gamma )
     \overset{ \eqref{gamma limit}}{ = } ( n )( 1+ 1/n ) = n + 1.  \nonumber
 \end{align}
Multiplying this relation with \eqref{asym ite, GD ratio} and use the fact $ K_{\mathrm{GD} }(\epsilon ) =  k_{\mathrm{GD} }(\epsilon )$  we get
\begin{equation}\label{asym ite, RPCD ratio}
   \lim_{c \rightarrow 1 } \frac{ K_{\mathrm{CCD} }(\epsilon ) }{ K_{\mathrm{RPCD} }(\epsilon ) } > \frac{n (n+1) }{2 \pi^2 }¡¡.
 \end{equation}

 According to \eqref{asym ite, GD ratio}, \eqref{asym ite, RCD ratio} and \eqref{asym ite, RPCD ratio}, there exists $c$ such that
all three relations in \eqref{ite, compare CD with others} hold.  %  \QED

 \subsubsection{Proof of Claim \ref{claim: RP-CD rate for Ac}}\label{appen: RP-CD rate proof}

For simplicity, we denote $L = L_{12\dots n}$.
Since $A = A_c$, we have
$$
  L = \begin{bmatrix}
     1 &    0   &  \dots &   0        \\
     c &    1   &  \dots &   0      \\
     \vdots &  \vdots &  \ddots & \vdots    \\
     c   &  c   &  \dots & 1
  \end{bmatrix}
$$
It is easy to get (recall that $\hat{c} \triangleq 1 - c$)
$$
  \Gamma^{-1} = \begin{bmatrix}
     1              &    0             &   0    &    \dots    &      0  &      0        \\
     -c             &    1             &   0    &    \dots    &      0  &      0      \\
      -c \hat{c}    &   -c             &   1    &    \dots    &      0   &      0      \\
       -c \hat{c}^2  &  -c\hat{c}       &  -c    &    \ddots    &     \vdots   &    \vdots    \\
       \vdots       &  \vdots          & \ddots &    \ddots   &     \ddots   &     \vdots   \\
 -c \hat{c}^{n-2}   &-c \hat{c}^{n-3}  & -c \hat{c}^{n-4} & \dots &  -c &   1
  \end{bmatrix}
$$
Since $L_{\sigma}$ can be obtained by permuting the rows and columns of $L$, thus
$L_{\sigma}^{-1}$ can also be obtained by similar permutations based on $L_{\sigma}^{-1}$.
As the matrix $A_c$ has only two distinct values, we know that the expression of $E(L_{\sigma}^{-1})$
(the expectation is taken over the uniform distribution of permutations of $\{1,2,\dots, n \}$)
 must have the following form
 \begin{equation}\label{Linv expect expression}
   E(L_{\sigma}^{-1}) =   \begin{bmatrix}
     1 &    \gamma   &  \dots &   \gamma        \\
     \gamma &    1   &  \dots &   \gamma      \\
     \vdots &  \vdots &  \ddots & \vdots    \\
     \gamma   &  \gamma   &  \dots & 1
  \end{bmatrix},
 \end{equation}
 where $\gamma$ only depends on $c$.
Due to symmetry,  $\gamma$ must be the average of all off-diagonal entries of $\Gamma^{-1}$, i.e.
\begin{equation}\label{gamma def}
  \gamma = \frac{ (n-1)c + (n-2)c\hat{c} + \dots  + \hat{c}^{n-2} }{ n(n-1) } = \frac{ -n + (1 - \hat{c}^n)/c }{ n (n-1) }.
\end{equation}
By the expressions \eqref{Linv expect expression} and \eqref{Ac def}, we have
$$
  E(L_{\sigma}^{-1}) A = \begin{bmatrix}
     \alpha &    \tau   &  \dots &   \tau      \\
      \tau   &    \alpha  &  \dots &     \tau      \\
     \vdots &  \vdots &  \ddots & \vdots    \\
      \tau    &   \tau     &  \dots & \alpha
       \end{bmatrix},
$$
where $ \alpha = 1 + (n-1)c \gamma $, $\tau = c + \gamma + (n-2) c \gamma $.
The minimum eigenvalue of this matrix is
\begin{equation}\label{QA lam min of RCD}
  \lambda_{\min}(E(L_{\sigma}^{-1}) A) = \alpha - \tau = 1 + (n-1)c \gamma - c - \gamma - (n-2) c \gamma
  = (1-c)(1-\gamma).
\end{equation}

According to \eqref{RCD update equation}, we have (note that $\sigma^k$ is independent of $x^k$)
$$
 E(x^{k}) = E(I - L_{\sigma}^{-1}A) E(x^{k-1}),
$$
where the expectation is taken over the uniform distribution of permutations of $\{1,2,\dots, n \}$.
This implies
$$
 \| E(x^k) \| \leq \| I - E( L_{\sigma}^{-1}) A \|  \|E(x^{k-1})\|,
$$
which further implies
\begin{equation}\label{interim RPCD}
   \| E(x^k) \|^2  \leq \| I - E( L_{\sigma}^{-1}) A \|^{2k} \| x^{0}\|^2 .
 \end{equation}
It is easy to verify that $ I - E( L_{\sigma}^{-1}) A $ is a positive semidefinite matrix, thus
 $ \| I - E( L_{\sigma}^{-1}) A  \| = 1 - \lambda_{\min}( E( L_{\sigma}^{-1}) A ) $.
 Plugging this and \eqref{QA lam min of RCD} into \eqref{interim RPCD},
 we obtain the desired inequality \eqref{RP-CD, converge rate for Ac}.

% The expectation in \eqref{Linv expect expression} is taken over all $n!$ permutations.

% \newpage

\bibliographystyle{IEEEbib}
\vspace{0.3cm}
% \bibliography{refs}

{\footnotesize
\bibliography{refBCD}
}

\end{document}